\pdfoutput=1
\documentclass[11pt]{article}

\usepackage[textwidth=468pt,textheight=46\baselineskip]{geometry}
\usepackage[T1]{fontenc}                        %
\usepackage[utf8]{inputenc}                     %
\usepackage{lmodern}                            %
\usepackage{dsfont}                             %
\usepackage{mathrsfs}                           %
\usepackage{tikz}
\usetikzlibrary{cd}
\usepackage{xcolor}
\usepackage{booktabs}           %
\usepackage[inline,shortlabels]{enumitem}                           %
\usepackage{csquotes}                           %
\usepackage{amsmath}
\usepackage{amssymb}                        %
\usepackage{amsthm}                         %
\usepackage{thmtools}                        %
\usepackage{acro}                               %
\usepackage{microtype}                       %
\usepackage{xparse}                             %
\usepackage{etoolbox}                           %
\usepackage{aligned-overset}                    %
\usepackage{mathtools}                      %
\usepackage[comma,numbers,sort&compress]{natbib}            %
\PassOptionsToPackage{hyphens}{url}                             %
\usepackage[pdfstartview=FitH,      %
            breaklinks=true,        %
            bookmarksopen=true,     %
            bookmarksnumbered=true,  %
            pdfencoding=unicode,
            psdextra,
            colorlinks=true,
            linkcolor=black,
            citecolor=blue!75!black,
            hypertexnames=false,
            urlcolor=blue!75!black%
            ]{hyperref}
\DeclareUrlCommand\doi{\def\UrlLeft##1\UrlRight{\color{blue!75!black}\href{http://dx.doi.org/##1}{ DOI\nobreakspace##1}}\urlstyle{rm}}
\usepackage[nameinlink]{cleveref}                           %

\usepackage{graphicx}
\graphicspath{{Images/},{./}}                                    %
\makeatletter
\@ifundefined{thmt@autorefsetup}{}{\@ifundefined{newcounteralias}{}{%
  \renewcommand\thmt@autorefsetup{\@xa\def\csname\thmt@envname autorefname\@xa\endcsname\@xa{\thmt@thmname}}%
}}
\makeatother
\theoremstyle{plain}
\newtheorem{theorem}{Theorem}[section]
\newtheorem{lemma}[theorem]{Lemma}
\newtheorem{corollary}[theorem]{Corollary}
\newtheorem{example}[theorem]{Example}
\theoremstyle{definition}
\newtheorem{definition}[theorem]{Definition}
\newtheorem{remark}[theorem]{Remark}
\allowdisplaybreaks[1]                          %
\newenvironment{delayedproof}[1]
 {\begin{proof}[\MakeLinkTarget*{#1}Proof of \autoref{#1}]}
 {\end{proof}}
\newcommand{\proofref}[1]{\hyperlink{#1}{proof}}
\setlist[1]{labelindent=\parindent} %
\setlist[itemize]{leftmargin=*}
\setlist[enumerate,1]{
  label = \alph*),
  ref = (\alph*)}

\crefrangeformat{section}{#3Sections #1#4--#5#2#6}
\Crefrangeformat{section}{#3Sections #1#4--#5#2#6}
\setlist[itemize]{  
  wide,
}
\makeatletter
\DeclareRobustCommand{\eqref}{\@ifstar{\supp@eqref{\ref*}}{\supp@eqref{\ref}}}
\newcommand*{\supp@eqref}[2]{\textup{\tagform@{#1{#2}}}}
\makeatother

\DeclarePairedDelimiter{\parenthesis}\lparen\rparen
\DeclarePairedDelimiter{\normBase}\lVert\rVert
\DeclarePairedDelimiter\absBase{\lvert}{\rvert}
\DeclarePairedDelimiter{\bracket}{\lbrack}{\rbrack}
\makeatletter
\long\def\NF@ifnext#1#2#3{%
    \let\NF@tok= #1\def\NF@yes{#2}\def\NF@no{#3}\futurelet\NF@next\NF@ifnext@}
\def\NF@ifnext@{\ifx\NF@next\@sptoken\let\NF@do\NF@ifnext@sp\else\ifx\NF@next\NF@tok\let\NF@do\NF@yes\else\let\NF@do\NF@no\fi\fi\NF@do}
\@firstofone{\def\NF@ifnext@sp} {\futurelet\NF@next\NF@ifnext@}
\def\NF@so#1#2{%
    \NF@ifnext*{\@firstoftwo{\NF@so@s{#1}{#2}}}{\NF@so@n{#1}{#2}}}
\def\NF@so@s#1#2{\NF@ifnext[{\NF@so@so{#1}{#2}}{#2{#1*}{#1*}{}}}
\def\NF@so@so#1#2[#3]{#2{#1*}{#1*}{#3}}
\def\NF@so@n#1#2{\NF@ifnext[{\NF@so@no{#1}{#2}}{#2{#1}{#1}{}}}
\def\NF@so@no#1#2[#3]{#2{#1[#3]}{#1}{#3}}
\def\NF@named#1#2#3{\NF@ifnext\bgroup{#1}{}}%
\def\NF@namedbracket#1#2#3{\NF@ifnext\bgroup{#1}{\NF@bracketlast{#1}{#2}{#3}}}
\def\NF@bracketlast#1#2#3{\NF@ifnext[{\NF@bracketlast@o{#1}}{#2{#3}}}%
\def\NF@bracketlast@o#1[#2]{#1{#2}}
\def\NF@before#1#2#3#4#5{%
    \NF@ifnext#1{\NF@before@t{#1}{#2}{#3}}{\NF@ifnext\bgroup{\NF@before@g{#3}}{}}}
\def\NF@before@t#1#2#3#4{\NF@ifnext\bgroup{\NF@before@tg{#1}{#2}{#3}}{}}
\def\NF@before@tg#1#2#3#4{#1{#2#4}\NF@ifnext\bgroup{#3}{}}
\def\NF@before@g#1#2{#1{#2}\NF@ifnext\bgroup{\@firstofone}{}}
\def\NF@beforebracket#1#2#3#4#5{%
    \NF@ifnext#1{\NF@beforebracket@t{#1}{#2}{#3}{#4}{#5}}{\NF@ifnext\bgroup{\NF@beforebracket@g{#3}}{\NF@bracketlast{#3}{#4}{#5}}}}
\def\NF@beforebracket@t#1#2#3#4#5#6{\NF@ifnext\bgroup{\NF@beforebracket@tg{#1}{#2}{#3}{#4}{#5}}{\NF@bracketlast{#3}{#4}{#5}}}
\def\NF@beforebracket@tg#1#2#3#4#5#6{#1{#2#6}\NF@ifnext\bgroup{#3}{\NF@bracketlast{#3}{#4}{#5}}}
\def\NF@beforebracket@g#1#2{#1{#2}\NF@ifnext\bgroup{\@firstofone}{}}
\def\NF@delimiter#1#2#3#4{\NF@ifnext\bgroup{\NF@delimiter@g{#1}{#2}}{#2{\:\cdot\:}#1}}
\def\NF@delimiter@g#1#2#3{#2{\ifblank{#3}{\:\cdot\:}{#3}}#1}
\def\NF@DefineDelimiter@o#1#2[#3]{\NF@DefineDelimiter{#1}{#2}{_{#3}}}
\def\NF@DefineDelimiter#1#2#3{\DeclareRobustCommand{#1}{\NF@so{#2}{\NF@delimiter{#3}}}}
\def\NF@delimiterbefore#1#2#3#4#5{#4#1{#5}}
\newcommand{\NamedFunction}[2]{%
    \DeclareRobustCommand{#1}{\operatorname{#2}\NF@so\parenthesis\NF@named}%
}
\newcommand{\NamedFunctionbefore}[2]{%
    \DeclareRobustCommand{#1}{#2\NF@so\parenthesis{\NF@before_{}}}%
}%
\newcommand{\NamedFunctionbeforespace}[2]{%
    \DeclareRobustCommand{#1}{#2\NF@so\parenthesis{\NF@before_{\!\!}}}%
}%
\newcommand{\NamedFunctionbeforeup}[2]{%
    \DeclareRobustCommand{#1}{#2\NF@so\parenthesis{\NF@before^{}}}%
}%
\newcommand{\NamedFunctionbeforeBracket}[2]{%
    \DeclareRobustCommand{#1}{#2\NF@so\bracket{\NF@beforebracket_{}}}%
}%
\newcommand{\DefineDelimiter}[2]{%
    \NF@ifnext[{\NF@DefineDelimiter@o{#1}{#2}}{\NF@DefineDelimiter{#1}{#2}{}}}%
\newcommand{\DefineDelimiterBefore}[2]{%
    \DeclareRobustCommand{#1}{\NF@so{#2}\NF@delimiterbefore}}%
\makeatother
\newcommand{\EE}{\mathbb{E}}
\newcommand{\FF}{\mathbb{F}}
\newcommand{\NN}{\mathbb{N}}
\newcommand{\PP}{\mathbb{P}}
\newcommand{\RR}{\mathbb{R}}
\newcommand{\AC}{\mathcal{A}}
\newcommand{\BC}{\mathcal{B}}
\newcommand{\DC}{\mathcal{D}}
\newcommand{\MC}{\mathcal{M}}
\newcommand{\NC}{\mathcal{N}}
\newcommand{\PC}{\mathcal{P}}
\newcommand{\QC}{\mathcal{Q}}
\newcommand{\WC}{\mathcal{W}}
\newcommand{\XC}{\mathcal{X}}
\newcommand{\YC}{\mathcal{Y}}
\DeclareMathOperator{\Exists}{\exists}                          %
\providecommand\given{}
\newcommand\SetSymbol[1][]{%
    \nonscript\:#1\vert
    \allowbreak
    \nonscript\:
    \mathopen{}}
\DeclarePairedDelimiterX\set[1]\{\}{%
\renewcommand\given{\SetSymbol[\delimsize]}
#1
}
\providecommand{\setcount}{}
\NamedFunction{\setcount}{\#}%
\newcommand*{\defeq}{\mathrel{\vcenter{\baselineskip0.5ex \lineskiplimit0pt \hbox{\scriptsize.}\hbox{\scriptsize.}}}=}                     %
\renewcommand{\defeq}{\coloneqq}
\newcommand*{\eqdef}{=\mathrel{\vcenter{\baselineskip0.5ex \lineskiplimit0pt \hbox{\scriptsize.}\hbox{\scriptsize.}}}}                           %
\renewcommand{\eqdef}{\eqqcolon}
\newcommand{\R}{\RR}                         %
\newcommand{\Rbar}{\bar{\R}}                               %
\newcommand{\N}{\NN}                         %
\newcommand{\tensor}{\otimes}
\DeclareMathOperator{\Id}{id}
\newcommand{\id}{\Id}
\providecommand{\Functionimage}{}
\DefineDelimiterBefore{\Functionimage}{\parenthesis}
\makeatletter
\DeclareRobustCommand{\preimage}{\NF@so\Functionimage\NF@preimage}%
\def\NF@preimage#1#2#3#4#5{#1{{#4}^{-1}}{#5}}
\makeatother
\NamedFunction{\exp}{exp}
\NamedFunction{\log}{log}
\DeclarePairedDelimiterX{\restrBase}[1]{.}{\rvert}{%
    \kern-\nulldelimiterspace %
    #1 %
    \littletaller %
}
\newcommand{\littletaller}{\mathchoice{\vphantom{\big|}}{}{}{}}
\makeatletter
\DeclareRobustCommand{\restr}{\NF@ifnext*{\@firstoftwo{\NF@restr\@firstoftwo}}{\NF@restr\@secondoftwo}}%
\def\NF@restr#1{\NF@ifnext[{\NF@restr@o#1}{\NF@restr@t{#1{\restrBase}{\restrBase*}}}}
\def\NF@restr@o#1[#2]{\NF@restr@t{#1{\restrBase[#2]}{\restrBase*}}}
\def\NF@restr@t#1{\NF@ifnext_{\NF@restr@u{#1}}{\NF@restr@m{#1}}}
\def\NF@restr@u#1#2{\NF@restr@m{#1}}
\def\NF@restr@m#1#2#3{#1{#2}_{#3}}
\makeatother
\DeclareMathOperator*{\argmin}{arg\,min}
\DeclareMathOperator*{\essliminf}{ess\,lim\,inf}
\DeclareMathOperator*{\essinf}{ess\,inf}
\NamedFunction{\sgn}{sgn}               %
\newcommand{\convexindex}{w}
\newcommand{\weightings}[1]{\Delta^{#1}}             %
\newcommand{\nondec}{
{
\pgfmathsetmacro{\scalefactor}{0.55}
\pgfmathsetmacro{\etemp}{0.86 * \scalefactor}
\pgfmathsetmacro{\ftemp}{0.88 * \scalefactor}
\mathchoice
{{\tikz[baseline=0.18em]{\path (0,0);\draw[->, line width=0.51pt, >={To[length=0.27em, width=0.55em]}] (0, 0.02em) --
 (0.86em, 0.88em);}}}%
{{\tikz[baseline=0.18em]{\path (0,0);\draw[->, line width=0.51pt, >={To[length=0.27em, width=0.55em]}] (0, 0.02em) --
 (0.86em, 0.88em);}}}%
{\!{\tikz[baseline=0.1375em]{\path (0,0);\draw[->, line width=0.40pt,
 >={To[length=0.1485em, width=0.3025em]}] (0, 0.011em) -- (0.473em, 0.484em);}}}%
{\!{\tikz[baseline=0.1375em]{\path (0,0);\draw[->, line width=0.40pt,
 >={To[length=0.1485em, width=0.3025em]}] (0, 0.011em) -- (0.473em, 0.484em);}}}%
}
}
\providecommand{\increasingfunctionclasses}{}
\NamedFunctionbeforeup{\increasingfunctionclasses}{D_\nondec}
\providecommand{\increasingfunctions}{}
\NamedFunction{\increasingfunctions}{\DC_\nondec}
\newcommand{\monotonetransformation}[1]{\prescript{}{*}{#1}}
\newcommand{\integratedgeneralizedinverse}[1]{K_{#1}}
\providecommand{\generalizedinverse}{}
\NamedFunction{\generalizedinverse}{Inv}
\makeatletter
\DeclareRobustCommand{\interior}{\NF@so\parenthesis\NF@interior}
\def\NF@interior#1#2#3#4{#1{#4}^\circ}
\makeatother
\newcommand{\closure}[1]{\overline{#1}}
\DefineDelimiter{\norm}{\normBase}%
\DefineDelimiter{\abs}{\absBase}%
\makeatletter
\def\NF@inner#1{\NF@ifnext[{\NF@inner@o#1}{#1{\scalprod}{\scalprod*}}}
\def\NF@inner@o#1[#2]{#1{\scalprod[#2]}{\scalprod*}}
\makeatother     %
\newcommand{\emptyconvergencearrow}[1]{\xrightarrow{#1}}
\newcommand*{\convN}{n \to \infty}
\newcommand{\nto}{\emptyconvergencearrow{\convN}}
\newcommand{\RadonNikodym}[2]{\frac{\diff#1}{\diff#2}}
\newcommand*\diff{\mathop{}\!\mathrm{d}}
\newcommand{\dd}{\diff}
\newcommand{\X}{\XC}                        %
\newcommand{\Y}{\YC}                        %
\newcommand{\A}{\AC}                        %
\NamedFunctionbeforeup{\Prob}{\PP}%
\renewcommand{\P}{\Prob}
\newcommand{\popmeasure}{\mathfrak{m}}
\newcommand{\randommeasure}{\kappa}
\makeatletter
\def\NF@distributionson[#1]{\NF@ifnext[{\NF@distributionson@ii{#1}}{\setofdistributions_{}(#1)}}
\def\NF@distributionson@ii#1[#2]{\setofdistributions_{#2}(#1)}
\makeatother
\providecommand{\distributionon}{}
\NamedFunctionbeforespace{\distributionon}{\PC}
\newcommand{\Lp}{\mathscr{L}}                       %
\newcommand{\Lpmod}{L}                       %
\NamedFunction{\B}{\BC}                        %
\newcommand{\Lebesgue}{\lambda}
\newcommand{\measuresetname}{\MC}
\providecommand{\measureon}{}
\NamedFunction{\measureon}{\Tilde{\measuresetname}}
\providecommand{\Radonmeasureon}{}
\NamedFunction{\Radonmeasureon}{\measuresetname_R}
\providecommand{\locallyfinitemeasureon}{}
\NamedFunction{\locallyfinitemeasureon}{\measuresetname_{l}}
\providecommand{\boundedlyfinitemeasureon}{}
\NamedFunction{\boundedlyfinitemeasureon}{\measuresetname_{b}}
\providecommand{\sigmafinitemeasureon}{}
\NamedFunction{\sigmafinitemeasureon}{\measuresetname_{\sigma}}
\providecommand{\generatedmeasureon}{}
\NamedFunction{\generatedmeasureon}{\measuresetname_{\CDFtransform}}
\makeatletter
\def\NF@measureson[#1]{\NF@ifnext[{\NF@measureson@ii{#1}}{\MC_{}(#1)}}
\def\NF@measureson@ii#1[#2]{\MC_{#2}(#1)}
\makeatother
\newcommand{\Deq}{\overset{\DistributionD}{=}}
\providecommand{\support}{}
\NamedFunction{\support}{spt}
\newcommand{\pushforward}[2]{#1_\##2}
\newcommand{\ind}{\mathds{1}}               %
\newcommand{\1}{\ind}
\newcommand{\CDFset}{D_{\text{{c.d.f.\@}\nocorr}}}
\newcommand{\QFset}{D_{\text{{q.f.\@}\nocorr}}}
\providecommand{\Quantiletransform}{}
\NamedFunction{\Quantiletransform}{\QC}
\providecommand{\CDFtransform}{}
\NamedFunction{\CDFtransform}{\mathfrak{F}}
\providecommand{\Stieltjestransform}{}
\NamedFunction{\Stieltjestransform}{\CDFtransform^{-1}}
\providecommand{\Transportmeasure}{}
\NamedFunction{\Transportmeasure}{\Quantiletransform^{-1}}
\providecommand{\InverseMeasure}{}
\NamedFunction{\InverseMeasure}{\mathcal{I}}
\providecommand{\DensSelec}{}
\NamedFunction{\DensSelec}{\mathfrak{f}}
\providecommand{\QuantDensSelec}{}
\NamedFunction{\QuantDensSelec}{\mathfrak{q}}
\makeatletter
\DeclareRobustCommand{\uniformmeasure}{\NF@ifnext^{\NF@uniformmeasure}{\restr{\Lebesgue}{(0,1)}}}%
\def\NF@uniformmeasure^#1{\restr{\Lebesgue}{(0,1)^{#1}}}
\makeatother
\NamedFunction{\Normal}{\NC}
\newcommand{\gaussian}{\Normal}
\NamedFunction{\betad}{Beta}%
\newcommand{\Betadist}{\betad}
\newcommand{\meansymbol}{\EE}
\providecommand{\mean}{}
\NamedFunctionbeforeBracket{\mean}{\meansymbol}%
\providecommand{\median}{}
\NamedFunction{\median}{med}
\providecommand{\midrange}{}
\NamedFunction{\midrange}{midrange}
\DeclareMathOperator{\DistributionD}{\DC}
\newcommand{\weakconvergencesymbol}{w}
\newcommand{\vagueconvergencesymbol}{v}
\newcommand{\vagueto}{\emptyconvergencearrow{\vagueconvergencesymbol}}
\newcommand{\coupling}{\gamma}
\NamedFunction{\transportplans}{\Pi}           %
\NamedFunction{\W}{\WC}                        %
\NamedFunctionbefore{\Wp}{\WC_{p}}
\NamedFunctionbefore{\Wq}{\WC_{q}}
\NamedFunctionbefore{\Wtwo}{\WC_{2}}
\newcommand{\Barycenterfunc}[2]{\Lambda_{#1}^{#2}}
\newcommand{\WB}[1]{$#1$-{WB}}
\makeatletter
\def\NF@optimallineartransport[#1]{A^*_{#1}}
\makeatother
\makeatletter
\def\NF@optimalrotation[#1]{U^*_{#1}}
\makeatother
\DeclareAcronym{RHS}{
    short = r.h.s\acdot,
    long = right-hand side,
    first-style = short
}
\DeclareAcronym{LHS}{
    short = l.h.s\acdot,
    long = left-hand side,
    first-style = short
}
\DeclareAcronym{CDF}{
    short = c.d.f\acdot,
    long = cumulative distribution function
}
\DeclareAcronym{QF}{
    short = q.f\acdot,
    long = quantile function
}
\DeclareAcronym{pWB}{
    long = $p$-Wasserstein barycenter,
    short = $p$-WB
}
\DeclareAcronym{2WB}{
    long = $2$-Wasserstein barycenter,
    short = $2$-WB
}
\DeclareAcronym{WB}{
    long = Wasserstein barycenter,
    short = WB,
    first-style = short
}
\DeclareAcronym{ae}{
    short = a.e\acdot,
    long = almost everywhere,
    first-style = short
}
\DeclareAcronym{wrt}{
    short = w.r.t\acdot,
    long = with respect to,
    first-style = short
}
\DeclareAcronym{ie}{
    short = i.e.\accomma,
    long = that is,
    first-style = short
}
\DeclareAcronym{eg}{
    short = e.g\acdot,
    long = for example,
    first-style = short
}
\newcommand{\abbrcomma}{,}
\newcommand*{\accomma}{\aciftrailing{comma}{}{\abbrcomma}}
\makeatletter
\patchcmd{\section}{\bfseries}{\bfseries\boldmath}{}{\PackageWarning{arxiv_draft}{Patching \string\section\space failed}}
\patchcmd{\subsection}{\bfseries}{\bfseries\boldmath}{}{\PackageWarning{arxiv_draft}{Patching \string\subsection\space failed}}
\patchcmd{\subsubsection}{\bfseries}{\bfseries\boldmath}{}{\PackageWarning{arxiv_draft}{Patching \string\subsubsection\space failed}}
\patchcmd{\paragraph}{\bfseries}{\bfseries\boldmath}{}{\PackageWarning{arxiv_draft}{Patching \string\paragraph\space failed}}
\makeatother

\begin{document}
\def\dobib{}

\title{Structure Preservation of \texorpdfstring{$2$}{2}-Wasserstein Barycenters on the Real Line with Applications to Log-Concave Distributions}
\author{Markus Zobel%
    \thanks{%
        Institute for Mathematical Stochastics, %
        University of G\"ottingen, %
        Email: \href{mailto:markus.zobel@uni-goettingen.de}{\nolinkurl{markus.zobel@uni-goettingen.de}}. %
        }%
        \and Axel Munk%
    \thanks{%
        Institute for Mathematical Stochastics, %
        University of G\"ottingen, %
        Email: \href{mailto:munk@math.uni-goettingen.de}{\nolinkurl{munk@math.uni-goettingen.de}}. %
    A.M. acknowledges support from DFG RU 5381 Mathematical Statistics in the Information Age.}%
}
\date{November 2026} %
\maketitle
\begin{abstract}

    The $p$-Wasserstein metric, $p\geq1$, is a popular distance between probability distributions. When the ground space is the real line, these metrics become particularly simple, as they can be equivalently expressed in terms of quantile functions. Consequently, the corresponding Fréchet mean of univariate distributions, denoted as the $p$-Wasserstein barycenter ($p$-WB), can also be expressed directly in terms of quantile functions and is hence straightforward to calculate.
    We study to what extent the $p$-WB preserves the structure of the underlying distributions, focusing on shape constraints on the densities, such as monotonicity, unimodality, or log-concavity.
    Until now, this question has been predominantly treated for $2$-WBs computed from a finite number of absolutely continuous input distributions.
    Here, we consider $p$-WBs of arbitrary collections of distributions. We formulate structure preservation as a generalized convexity property in Wasserstein space, which we denote as $p$-WB convexity.
    As the literature related to this question is scattered, we provide a unifying exposition for the univariate setting.
    Our approach is to find a quantile representation of the structural property of interest and then explore its preservation under Fréchet means on $L^p(0,1)$.
    We investigate $p$-WB convexity for various subfamilies of unimodal distributions. To this end, we introduce the notion of $\phi$-unimodality: an absolutely continuous distribution $\mu$ is $\phi$-unimodal if composition with $\phi$ turns a density into a concave function on the interior of the convex support of $\mu$.
    This notion thus generalizes log-concavity and $s$-concave densities.
    We show that $2$-WB convexity of univariate $\phi$-unimodal distributions depends on the regularity property enforced by $\phi$. 
    The family of univariate log-concave distributions is $2$-WB convex, while the family of univariate distributions with $s$-concave densities is $2$-WB convex if and only if $s\geq-2$. General unimodality, corresponding to $s=-\infty$, is not preserved.
    Our proofs rely on a characterization of $\phi$-unimodal families by a mean-convexity property of their quantile density functions, which also extends the notion of $\phi$-unimodality beyond absolutely continuous distributions.
    To use this characterization, we extend the closed-form expression for the quantile function of the $2$-WB to its quantile density.
    In higher dimensions or for $p\neq2$, we give counterexamples ruling out the $p$-WB convexity of $\phi$-unimodal distributions.
    Finally, we discuss applications of our results to shape-constrained density estimation.

\end{abstract}


\section{Introduction}

Distributional data \cite{kleinDistributionalRegressionData2024,kneipInferenceDensityFamilies2001,kneibRageMeanReview2023,petersenFunctionalDataAnalysis2016}, {e.g.\@}, in the context of functional data \cite{ramsayAppliedFunctionalData2002,ferratyNonparametricFunctionalData2006,wangFunctionalDataAnalysis2016,Panaretos.2020}, play an important role in modern data analysis.
In this context, $q$th power Fréchet means (see, {e.g.\@}, \cite{jaffeFrechetMeansInfinite2026,schotzWeightedPowerFrechet2026}), {i.e.}, barycenters, are of particular interest as they can summarize important structural features (such as modes) of an ensemble of distributions $\popmeasure$ on a given ground (metric) space $(\X,d)$.
They are obtained as minimizers of a Fréchet functional $\FF_{d_{\distributionon},q}$, induced by the choice of an extended metric $d_{\distributionon}:\distributionon(\X)^2\to[0,+\infty]$ on the space of distributions $\distributionon(\X)$ on $(\X,\B(\X))$ and an exponent $q>0$,
\begin{equation*}
    \FF_{d_{\distributionon},q}(\popmeasure,\mu)=\mean_{\randommeasure\sim\popmeasure}[d_\PC^q(\randommeasure,\mu)],\quad \mu\in\distributionon(\X).
\end{equation*}
Here, $\randommeasure:(\Omega,\A,\P)\to(\distributionon(\X),\sigma(\iota_{\B(\X)}))$ is a \emph{random distribution}, a measurable map from an ambient complete probability space $(\Omega,\A,\P)$ to $\distributionon(\X)$ equipped with the $\sigma$-algebra initial to all evaluations $\iota_B$ on Borel sets $B\in\B(\X)$ \cite{kallenbergRandomMeasuresTheory2017}.
\subsection{Wasserstein barycenter convexity}
    Prominent examples of Fréchet means are the $(p,q)$-Wasserstein barycenters, which have received considerable attention recently ranging from applications to theory \cite{lecozBarycenterbasedApproachMultimodel2026,singhModelFusionOptimal2020,zhuangWassersteinKmeansClustering2022,Panaretos.2020,
    lauWassersteinDistributionallyRobust2022,srivastavaScalableBayesBarycenter2018,zemelFrechetMeansProcrustes2019,zhouWassersteinRegressionEmpirical2024,bigotDatadrivenRegularizationWasserstein2019}. Let $(\X,d)$ be a Polish metric space and let $p,q\geq1$. Then
    \begin{equation}
        \label{eqn: def Wasserstein barycenter}
        \Barycenterfunc{q}{\Wp}(\popmeasure)=\argmin_{\mu\in\Wp(\X)}\FF_{W_{\!p},q}(\popmeasure,\mu)=\argmin_{\mu\in\Wp(\X)}\mean_{\randommeasure\sim\popmeasure}[W_{\!p}^q(\randommeasure,\mu)]
    \end{equation}
    denotes the set of all $(p,q)$-Wasserstein barycenters.
    Here, the \emph{$p$-Wasserstein space} $\Wp(\X)$ for $p\geq1$ is the metric space of distributions with finite $p$th metric moment
    \begin{equation*}
        \distributionon_{p}(\X)\defeq\set{\mu\in\distributionon{\X}\given\Exists x_0\in\X, \int d^p(x_0,x)\dd{\mu}(x)<\infty}
    \end{equation*}
    together with the \emph{$p$-Wasserstein metric}, the $p$th root of the optimal transport cost {w.r.t.\@} the $p$th power $d^p$ of the metric,\footnote
    {
        For two distributions $\mu,\nu\in\distributionon(\X)$, the set $\transportplans{\mu,\nu}\subset\distributionon(\X^2)$ denotes the set of all transport plans.
        The marginals of a transport plan reproduce the original distributions, {i.e.}, they satisfy $(\pi_1)_\#\coupling = \mu$, $(\pi_2)_\#\coupling = \nu$. For more information on optimal transport and Wasserstein barycenters, we refer the reader to \cite{Santambrogio.2020,Panaretos.2020}.
    }
    \begin{equation*}
        W_{\!p}^p(\mu,\nu)=\inf_{\coupling\in\transportplans{\mu,\nu}}\int d^p(x,y)\dd{\coupling}(x,y),\quad \mu,\nu\in\distributionon_{p}(\X),\quad p\geq1.
    \end{equation*}
    The Fréchet functional $\FF_{W_{\!p},q}$ is finite if and only if $\popmeasure\in\Wq(\Wp(\X))$. Under mild regularity conditions on the metric space $(\X,d)$, satisfied by $\X=\R^d$ with the Euclidean metric, a minimizer in \eqref{eqn: def Wasserstein barycenter} exists \cite{LeGouic.2017}. For $\popmeasure\in\Wq(\Wp(\X))$ such that the Wasserstein barycenter is unique, {i.e.}, $\Barycenterfunc{q}{\Wp}(\popmeasure)=\set{\mu^*}$, we interpret $\Barycenterfunc{q}{\Wp}$ as a distribution-valued map $\Barycenterfunc{q}{\Wp}(\popmeasure)=\mu^*$.

    Most prominently, minimizing the squared $2$-Wasserstein distance $(p=q=2)$ leads to the {{$2$-Wasserstein barycenter} ({$2$-WB})} \cite{aguehBarycentersWassersteinSpace2011}.
    Generalizations to the {{$p$-Wasserstein barycenter} ({$p$-WB})} ($p=q\geq 1$) have been investigated in \cite{LeGouic.2017,brizziWassersteinBarycenters2025}.
    The term Wasserstein median has been coined by multiple authors, in all cases referring to $q=1$ but not necessarily also to $p=1$ \cite{carlierWassersteinMediansRobustness2024,youWassersteinMedianProbability2025}.
    Further generalizations to Wasserstein barycenters of finite measures based on unbalanced optimal transport are possible \cite{frieseckeBarycentersHellingerKantorovichDistance2021,Heinemann.2023}.
    Wasserstein barycenters of distributions of pairwise distances occur in the context of the Gromov-Wasserstein distance between metric measure spaces \cite{steinkampGromovWassersteinBarycenterSurrogates2026}.
    The most studied tool remains the classical {$2$-WB} \cite{bigotStatisticalDataAnalysis2020,chewiStatisticalOptimalTransport2025}, partly because of its geometry. For example, for $\X=\R$, the $2$-Wasserstein space $\Wtwo(\R)$ is isometric to a subset of the flat Hilbert space $\Lpmod^2{(0,1)}$. In higher dimensions $d\geq2$, $\Wtwo(\R^d)$ has nonnegative curvature \cite{ambrosioGradientFlowsMetric2005}.

    Despite the conceptual appeal of Wasserstein barycenters, the literature on their preservation of statistically relevant features of the input distributions is scattered, even for the {$2$-WB} on $\Wtwo(\R)$.
    Formally, this preservation can be phrased in terms of convexity in the Wasserstein space. We call a family of distributions sharing such a feature {$p$-WB} convex if the feature is preserved upon taking arbitrary {$p$-WBs} \cite{AlvarezEsteban.2018,backhoffStochasticGradientDescent2025}.
    \begin{definition}[{$p$-WB} convex]
        \label{def: Wasserstein p barycenter convexity}
        A subset $A\subset\distributionon(\R^d)$ is called \emph{{$p$-WB} convex} for $p\geq1$ if for all random distributions $\popmeasure\in\mathfrak{M}_{+\infty}\defeq\Wp(\Wp(\R^d))$ \emph{supported on $A$}, {i.e.}, $\overline{\popmeasure}(A\cap\W_p(\R^d))=1$, where $\overline{\popmeasure}$ is the inner measure\footnote{
            Let $(\X,\B,\mu)$ be a measure space. Denote by $\overline{\mu}$ the \emph{inner measure} defined for all $A\subset \X$ as $\overline{\mu}(A)\defeq\sup\set{\mu(B)\given A\supset B\in\B}$. If $\randommeasure:(\Omega,\A,\P)\to(\X,\B)$ is a random variable with distribution $\randommeasure\sim\mu\in\distributionon(\X)$ defined on the complete ambient probability space $(\Omega,\A,\P)$, then $\overline{\mu}(A)=1$ is equivalent to $A$ having a subset $B\in\B$ with $\mu(B)=1$, which in turn implies $\P(\randommeasure\in A)=1$.
        }
        of $\popmeasure$, every {$p$-WB} of $\popmeasure$ is again in $A$.
    \end{definition}
    If we replace $\mathfrak{M}_{+\infty}=\Wp(\Wp(\R^d))$ with different choices of subsets $\mathfrak{M}\subset\Wp(\Wp(\R^d))$ shown in \autoref{tab: different concepts of barycenter convexity}, we recover existing notions of convexity (see \autoref{sec: Relation of different notions of barycenter convexity results}).
    For $p=2$, \emph{$p$-displacement convexity} corresponds to McCann interpolations \cite{mccannConvexityPrincipleInteracting1997},\footnote{
        The McCann interpolations of two distributions $\mu_0,\mu_1\in\W_p(\R^d)$, $p>1$, are given by $\mu_{t}\defeq\pushforward{(t\pi_1+(1-t)\pi_2)}{\coupling_p^*}$, where $t\in(0,1)$ and $\coupling_p^*\in\transportplans{\mu_1,\mu_0}$ is any optimal coupling from the definition of $W_{\!p}$ \cite[Theorem 7.2.2]{ambrosioGradientFlowsMetric2005}.
        For $p=2$, the set of all McCann interpolations evaluated at $\convexindex$ thus coincides with the {$2$-WB} $\Barycenterfunc{2}{\Wtwo}(\popmeasure_{\convexindex})$ for $\popmeasure_{\convexindex}=\convexindex\delta_{\mu_1}+(1-\convexindex)\delta_{\mu_{0}}$. For general $p>1$, we have $\mu_{t}=\Barycenterfunc{p}{\Wp}(\popmeasure_{{\convexindex_t}})$, with $\convexindex_t$ a reweighting depending on $p$ \cite{brizziWassersteinBarycenters2025}.
    }
    meaning that the constant-speed $2$-Wasserstein geodesics between two measures $\mu_0,\mu_1\in A\cap\Wtwo(\R^d)$ remain in $A$. This characterization extends to $p>1$ with constant-speed $p$-Wasserstein geodesics instead.
    The extension of $p$-displacement convexity to {$p$-WBs} of finitely many distributions is called \emph{convexity along $p$-barycenters} and is considered for $p=2$ in \cite{aguehBarycentersWassersteinSpace2011}.

    In the univariate case, the flat geometry of $\Wtwo(\R)$ leads to a general approach to determining the {$2$-WB} convexity of sets using quantile functions.
    For spatial dimensions $d>1$, no such general approach is known, and individual proof strategies are necessary \cite{backhoffStochasticGradientDescent2025,AlvarezEsteban.2018}.
    \begin{table}[ht]
        \centering
        \caption{Different convexity notions based on Wasserstein barycenters differentiated by cardinality conditions on the support of $\popmeasure$.\protect\footnotemark
        }
        \label{tab: different concepts of barycenter convexity}
        \begin{tabular}{lll}\toprule
            \cite{ambrosioGradientFlowsMetric2005}      & $p$-displacement convexity      & $\mathfrak{M}_2=\set{\popmeasure\in\distributionon(\Wp(\R^d))\given \setcount{\support{\popmeasure}}=2}$     \\
            \cite{aguehBarycentersWassersteinSpace2011} & convexity along $p$-barycenters & $\mathfrak{M}_{<+\infty}=\set{\popmeasure\in\distributionon(\Wp(\R^d))\given \setcount{\support{\popmeasure}}<+\infty}$ \\
                                                        & {$p$-WB} convexity              & $\mathfrak{M}_{+\infty}=\Wp(\Wp(\R^d))$                                                                      \\\bottomrule
        \end{tabular}
    \end{table}
    \footnotetext{
            The support of $\mu\in\distributionon{\X}$ is the smallest closed set to which $\mu$ assigns full probability,
            {i.e.}, $\support{\mu}\defeq\bigcap_{A\in S_\mu} A,$ $S_{\mu}=\set{A\subset\X\given\text{$A$ is closed, $\mu(\X\setminus A)=0$}}$.
            }
\subsection{The univariate setting}
    \label{sec: the univariate setting}
    In the univariate case, the {$2$-WB} $\mu^*\defeq\Barycenterfunc{2}{\Wtwo}\,(\popmeasure)$ has a closed-form expression in terms of its {{quantile function} ({q.f.\@})} \cite{bigotCharacterizationBarycentersWasserstein2018,Panaretos.2020}
    \begin{equation}
        \label{eqn: characterization of Qunaitle fucntionas mean}
        Q_{\mu^*}(t)=\mean_{\randommeasure\sim\popmeasure} {Q_{\randommeasure}(t)},\quad t\in(0,1),\popmeasure\in\Wtwo(\Wtwo(\R)).
    \end{equation}
    The \emph{left-continuous {q.f.\@}\nocorr} $Q_\mu$ of $\mu\in\distributionon(\R)$ is a generalized inverse of the \emph{{{cumulative distribution function} ({c.d.f.\@})}\nocorr} $F_\mu(x)\defeq\mu((-\infty,x])$, $x\in\R$. It is given by $Q_\mu(t)\defeq\inf(\preimage{F_\mu}{[t,1]})$ for each $t\in(0,1)$.
    While this is the standard definition of a {q.f.\@}, we will instead use the \emph{right-continuous {q.f.\@}\nocorr} defined by
    \begin{equation*}
        Q_\mu(t)\defeq\sup(\preimage{F_\mu}{[0,t]}), \quad t\in(0,1).
    \end{equation*}

    The closed-form expression in \eqref{eqn: characterization of Qunaitle fucntionas mean} is based on a) the isometry between the Wasserstein space $\Wp(\R)$ and the convex subset of $p$-integrable {q.f.s} in $\Lpmod^p(0,1)$
    \begin{equation*}
        W_{\!p}(\mu,\nu)=\norm{Q_\mu-Q_\nu}_{\Lpmod^p(0,1)},\quad \mu,\nu\in\Wp(\R),
    \end{equation*}
    which holds because the coupling $\coupling^*=\pushforward{(Q_{\mu},Q_{\nu})}{\uniformmeasure}$ is optimal for all $p\geq1$,\footnote{
        Therefore, on $\R$, we can call a set $A$ displacement convex if $\mu_t\defeq\pushforward{(tQ_{\mu_1}+(1-t)Q_{\mu_0})}{\uniformmeasure}\in A$ for all $\mu_0,\mu_1\in A$ without restricting to $\mu_0,\mu_1\in\W_p(\R)\cap A$. Thus, displacement convexity implies $p$-displacement convexity for all $1<p<+\infty$. Further, a set $A$ is displacement convex if the corresponding set of {q.f.s} is convex \cite[Definition 2.4]{matsudaShapeconstrainedDensityEstimation2026}.
    }
    b) the characterization of the Fréchet mean on $\Lpmod^2(0,1)$ as the pointwise average of random functions, and c) the non-decreasingness and continuity of this barycenter map, which ensure that the pointwise average of random {q.f.s} is again a {q.f.\@} \cite[Theorem 3.1, resp., Theorem 3.2.11]{bigotCharacterizationBarycentersWasserstein2018,Panaretos.2020}. The arguments proving \eqref{eqn: characterization of Qunaitle fucntionas mean} in \cite[Proof of Theorem 3.2.11]{panaretosSupplementaryProofsInvitation2019arxiv} can be directly adapted to right-continuous {q.f.s}{}.

    This closed-form expression suggests the following approach to checking the {$2$-WB} convexity of a family of univariate distributions:
    \begin{enumerate}
        \item Characterize the family based on a property of the {q.f}.
        \item Show that this property is preserved for the pointwise mean.
    \end{enumerate}
    Closed forms similar to \eqref{eqn: characterization of Qunaitle fucntionas mean} for $p\neq2$ have so far been reported only for $\setcount{\support{\popmeasure}}<+\infty$ \cite{brizziWassersteinBarycenters2025}, but would give rise to similar strategies for proving {$p$-WB} convexity.
    For now, we restrict ourselves to $p=2$ and only return to $p\neq2$ to establish counterexamples in \autoref{sec: Counterexamples for extensions of the main theorem}. 
\subsection{\texorpdfstring{$2$}{2}-WB convexity of absolute continuity of quantile functions}
    For our main result presented in \autoref{sec: main result}, it is useful to characterize a distribution $\mu\in\distributionon(\R)$ using a quantile density $q_\mu$ instead of the {q.f.\@} $Q_\mu$.
    We call a univariate distribution $\mu\in\distributionon(\R)$ \emph{quantile absolutely continuous} if $Q_\mu$ is locally absolutely continuous.\footnote{We call a function $G:I\to\R$ from an open interval $I$ \emph{locally absolutely continuous} if $\restr{G}{(a,b)}$ is absolutely continuous for all $a<b\in I$.} In that case, a \emph{quantile density} $q_\mu$ is any Radon-Nikodym derivative of $Q_\mu$. We show that the set of quantile absolutely continuous distributions is {$2$-WB} convex and establish a closed-form expression analogous to \eqref{eqn: characterization of Qunaitle fucntionas mean}.
    To this end, we require the existence of a quantile density selection, a jointly measurable function $\QuantDensSelec:\distributionon(\R)\times (0,1)\to\Rbar_{\geq0}$ satisfying $\QuantDensSelec(\mu,\cdot)=q_\mu$ $\uniformmeasure$-{a.e.\@} for any quantile density $q_\mu$ of a quantile absolutely continuous $\mu\in\distributionon(\R)$.
    \begin{theorem}
        \label{prop: Wasserstein 2 Barycenter of quantile absolutely continuous distributions on R have absolutely continuous quantile functions}
        Let $\randommeasure\sim\popmeasure\in\Wtwo(\Wtwo(\R))$ be supported on quantile absolutely continuous distributions, and denote its {$2$-WB} by $\mu^*=\Barycenterfunc{2}{\Wtwo}\,(\popmeasure)$. Then
        \begin{equation*}
            q_{\mu^*}(t)\defeq\mean[\QuantDensSelec(\randommeasure,t)],\quad t\in(0,1),
        \end{equation*}
        defines a valid quantile density for every quantile density selection
        $\QuantDensSelec$.
    \end{theorem}
    \begin{corollary}
        The set of quantile absolutely continuous distributions is {$2$-WB} convex.
    \end{corollary}

    The proof given in \autoref{sec: 2-WB convexity of quantile absolute continuity} is based on an analysis of the interplay between random distributions, random {q.f.s}, and the regularity of the mappings $\mu\mapsto F_\mu$ and $\mu\mapsto Q_\mu$. Following \cite{bobkovOnedimensionalEmpiricalMeasures2019}, we focus on the inverse measure $\mu^{-1}$, which is the Stieltjes measure of $Q_\mu$. The inverse measure of the {$2$-WB} is the intensity measure of the random inverse measures $\randommeasure^{-1}$, and we use results from random measure theory \cite{kallenbergRandomMeasuresTheory2017} to work with this object.

\subsection{Application to unimodality}
    We follow \cite{zobelmarkusQuantileCharacterizationUnivariate2026arxiv} to characterize unimodal distributions in terms of their quantile densities.
    An absolutely continuous distribution $\mu\in\distributionon(\R)$ is \emph{unimodal} if there exists a quasi-concave density $f_\mu$ of $\mu$ \cite{Dharmadhikari.1988}, that is,
    \begin{equation*}
        f_{\mu}(\convexindex x +(1-\convexindex)y)\geq \min(f_\mu(x),f_\mu(y)),\quad\convexindex\in(0,1),x,y\in\R.
    \end{equation*}
    This is equivalent to the existence of a mode $m\in\R$ such that a density $f_\mu$ is non-decreasing on $(-\infty,m)$ and non-increasing on $(m,+\infty)$ \cite{martosNonlinearProgrammingTheory1975}.
    \begin{theorem}[{\cite[Theorem 1.6]{zobelmarkusQuantileCharacterizationUnivariate2026arxiv}}]
        \label{Prop: Characterization of unimodal distributions}
        A measure $\mu\in\distributionon(\R)$ is unimodal if and only if its {q.f.\@} $Q_\mu$ is locally absolutely continuous and there exists a quantile density $q_\mu$ that is quasi-convex, that is,
        \begin{equation*}
            q_\mu(\convexindex x+ (1-\convexindex)y)\leq \max(q_\mu(x),q_\mu(y)),\quad\convexindex,x,y\in(0,1).
        \end{equation*}
    \end{theorem}
    A \emph{quantile mode} of a unimodal $\mu$ is then a point $\beta\in[0,1]$ such that a quantile density $q_\mu$ is non-increasing on $(0,\beta)$ and non-decreasing on $(\beta,1)$.\footnote{For the general definition of not necessarily absolutely continuous unimodal distributions, see the beginning of \autoref{sec: Quantile Characterization of Ms concave distributions}.}
    While this characterization seems promising in light of \autoref{prop: Wasserstein 2 Barycenter of quantile absolutely continuous distributions on R have absolutely continuous quantile functions}, quasi-convexity of the quantile density is not preserved by {$2$-WBs}, and thus unimodality is not {$2$-WB} convex as already claimed in \cite{backhoffStochasticGradientDescent2025}. We present a counterexample in \autoref{example: unimodal distributions are not 2 Wasserstein barycenter convex}.
    \begin{figure}[ht]
        \centering
        \includegraphics[width=\textwidth]{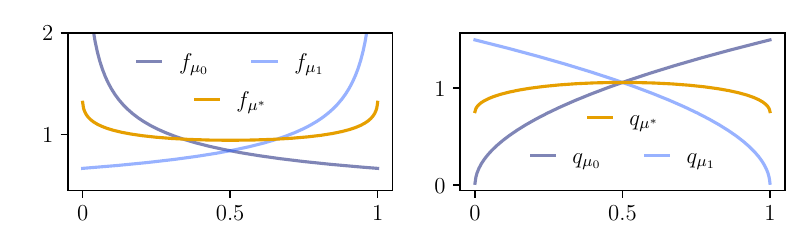}
        \caption{{L.h.s.\@}: While the original densities of the two beta distributions $\mu_0\sim\Betadist[\big]{\frac{2}{3},1}$, $\mu_1\sim\Betadist[\big]{1,\frac{2}{3}}$ (blue) each have one mode, the density of their equal-weight {$2$-WB} $\mu^*$ (orange) is bimodal. {R.h.s.\@}: The quantile densities of $\mu_0,\mu_1$ are quasi-convex as they are monotone, but their mean $q_{\mu^*}=\frac12(q_{\mu_0}+q_{\mu_1})$ is concave and not quasi-convex.}
        \label{fig: unimodal counterexample}
    \end{figure}
    \begin{example}[The set of unimodal distributions is not {$2$-WB} convex]
        \label{example: unimodal distributions are not 2 Wasserstein barycenter convex}
        Take two beta distributions $\mu_0\sim\Betadist[\big]{\frac{2}{3},1}$, $\mu_1\sim\Betadist[\big]{1,\frac{2}{3}}$ and their McCann interpolant $\mu^*=\Barycenterfunc{2}{\Wtwo}(\frac12(\delta_{\mu_0}+\delta_{\mu_1}))$. Both distributions are unimodal because their densities are monotone and hence quasi-concave.
        The density of their {$2$-WB} $\mu^*$ depicted in \autoref{fig: unimodal counterexample} has a global minimum at $\frac12$ and thus fails to be quasi-concave. Therefore, $\mu^*$ is not unimodal.
    \end{example}
    In the following, we consider subfamilies of unimodal distributions defined by concavity-like requirements on the density that are stronger than quasi-concavity.

\subsection{2-WB convexity of quantile \texorpdfstring{$\phi$}{ɸ}-unimodal distributions}
    Generalizing the quantile characterization of unimodality in \autoref{Prop: Characterization of unimodal distributions}, we introduce the notion of a quantile $\phi$-unimodal distribution.
    \begin{definition}
        \label{def: quantile phi-convex distribution}\begin{enumerate}[wide, labelwidth=!, labelindent=0pt,itemindent=!]
            \item
                  \label{enum: def: generator}
                  A strictly increasing and continuous function $\phi:\R_{>0}\to\R$ is called a \emph{generator} and generates a weighted \emph{quasi-arithmetic mean} $M_\phi$ by $M_{\phi,\convexindex}(x,y)=\phi^{-1}(\convexindex\phi(x)+ (1-\convexindex)\phi(y))$ for $(x,y)\in \R_{>0}^2$, $\convexindex\in(0,1)$.
            \item
                  \label{ enum: extension to Rbar}
                  We extend $M_\phi$ to $\Rbar_{\geq0}^2$,\footnote{
                      We write $\Rbar=\R\cup\set{-\infty,+\infty}$ for the usual set of extended real numbers, the order completion of $\R$ with $-\infty\defeq\inf(\R)$, $+\infty\defeq\sup(\R)$. We use the addition rules $\pm\infty+x=\pm\infty$ for all $x\in\R$.
                      Lastly, we denote by $\B(\Rbar)=\sigma(\B(\R)\cup\set{-\infty,+\infty})$ its Borel $\sigma$-algebra.
                  }
                  by setting $\phi(0)=\liminf_{t\searrow0}\phi(t)$ and $\phi(+\infty)=\limsup_{t\nearrow+\infty}\phi(t)$.
                  If $\phi(0)=-\infty$ and $\phi(+\infty)=+\infty$, then $M_\phi(0,+\infty)$ is not defined. Therefore, in that case, $M_{\phi}$ is extended only to $\Rbar_{\geq0}^2\setminus\set{(0,+\infty),(+\infty,0)}$.
            \item \label{enum: M_phi convex}
                  Let $I\subset\R^d$ be convex. A function $g: I \to\Rbar_{\geq0}$ is \emph{$M_\phi$-convex (-concave)} if
                  \begin{equation*}
                      g(\convexindex x+(1-\convexindex) y)\overset{(\geq)}{\leq} M_{\phi,\convexindex}(g(x),g(y)),\quad x,y\in I,\convexindex\in(0,1).
                  \end{equation*}
            \item
                  \label{enum: quantile phi unimodal}
                  A distribution $\mu\in\distributionon(\R)$ is \emph{quantile $\phi$-unimodal} if $\mu$ is quantile absolutely continuous and there exists a quantile density $q_\mu$ that is $M_\phi$-convex.
        \end{enumerate}
    \end{definition}
    As all quasi-arithmetic means are means, they satisfy $\min(x,y)\leq M_{\phi,\convexindex}(x,y)\leq\max(x,y)$ for all $\convexindex\in(0,1)$. Therefore, any $M_\phi$-convex quantile density is quasi-convex \cite{avrielGeneralizedConvexFunctions1974}, and thus any quantile $\phi$-unimodal distribution is unimodal.
    Our main theorem characterizes when the set of quantile $\phi$-unimodal distributions is {$2$-WB} convex.
    \begin{theorem}
        The set of quantile $\phi$-unimodal distributions is {$2$-WB} convex if and only if $M_\phi$ is concave.
    \end{theorem}

    The proof of our main theorem relies on the existence of a measurable selection of the $M_\phi$-convex representatives of the quantile density. This selection enables us to use \autoref{prop: Wasserstein 2 Barycenter of quantile absolutely continuous distributions on R have absolutely continuous quantile functions} in the \proofref{theo: Wasserstein 2 barycenter convexity of M convex funcitons}.
    \begin{lemma}
        \label{lem: existence of phi-convex modification of quantile density}
        There exists a quantile density selection $\QuantDensSelec:\distributionon(\R)\times(0,1)\to\Rbar_{\geq0}$ such that the following holds for every generator $\phi$: a quantile absolutely continuous distribution $\mu\in\distributionon(\R)$ is quantile $\phi$-unimodal if and only if $\QuantDensSelec(\mu,\cdot)$ is $M_\phi$-convex.
    \end{lemma}
    We give a \proofref{lem: existence of phi-convex modification of quantile density} in \autoref{sec: existence quantile density selection}.
\subsection{Characterization of quantile \texorpdfstring{$\phi$}{ɸ}-unimodal distributions}
    \label{sec: main result}
    While the definition of quantile $\phi$-unimodal distributions in \autoref{def: quantile phi-convex distribution} \ref{enum: quantile phi unimodal} is novel, the following notion of distributions with $M_\phi$-concave densities dates back to \cite{tammGCONCAVEFUNCTIONSPROBABILITY1977,avrielGeneralizedConvexFunctions1974}.
    \begin{definition}
        \label{def: phi unimodal}
        An absolutely continuous distribution $\mu\in\distributionon(\R^d)$ is \emph{$\phi$-unimodal} if its support $\support{\mu}$ is convex and there exists a density $f_\mu$ such that $f_\mu$ is $M_\phi$-concave when restricted to the interior of its support.
    \end{definition}
    \begin{example}
        \begin{enumerate}[wide, labelwidth=!, labelindent=0pt,itemindent=!]
            \item
                  An absolutely continuous distribution $\mu$ is called \emph{log-concave} if $\log\circ f_\mu$ is concave for some density $f_\mu$. This is equivalent to
                  \begin{equation*}
                      f_{\mu}(\convexindex x +(1-\convexindex)y)\geq\exp(\convexindex \log\circ f_\mu(x)+ (1-\convexindex) \log\circ f_\mu(y)),\quad x,y\in\R^d,\convexindex\in(0,1).
                  \end{equation*}
                  Additionally, the support of a log-concave distribution is convex, and thus every log-concave distribution is $\log$-unimodal, with $M_{\log}$ being the geometric mean. Furthermore, since the $M_{\log}$-concave density version $f_\mu$ of a distribution $\mu$ can be extended by $0$ outside the interior of its support, $\log\circ f_\mu:\R^d\to\R\cup\set{-\infty}$ is concave, so $\mu$ is log-concave.
            \item
                  More generally, $\phi_s$-unimodal distributions are known as \emph{distributions with $s$-concave densities} \cite{Dharmadhikari.1988}, where $\phi_s(x)=\sgn(s)x^s$ for $s\in\R\setminus\set{0}$, $x\in\R_{>0}$ and $\phi_0=\log$.
            \item
                  All power-law distributions, {i.e.}, absolutely continuous distributions with density $f_\mu(x)\propto\frac{1}{\abs{x}^a}$, $x\in\support{\mu}\subset\R_{\geq0}$, $a>0$, are $\phi_{-\frac1a}$-unimodal, {e.g.\@}, Pareto distributions. Similarly, Student's $t$ distribution with $\nu>0$ degrees of freedom is $\phi_{-\frac1{\nu+1}}$-unimodal.
        \end{enumerate}
    \end{example}
    In the univariate case, we can establish a correspondence between the two notions.
    \begin{theorem}
        \label{thm: equivilance between M concave and quantile M convex quasi arithmetic}
        Fix $c=\frac12$. Every generator $\phi$ as in \autoref{def: quantile phi-convex distribution} \ref{enum: def: generator} induces a co-generator $\monotonetransformation{\phi}$ defined by
        \begin{equation*}
            \monotonetransformation{\phi}(x)\defeq -\int_{\phi(\frac1c)}^{\phi(\frac1x)}\phi^{-1}\diff\Lebesgue +\phi(c),\quad x\in\Rbar_{\geq0},
        \end{equation*}
        such that every absolutely continuous $\mu\in\distributionon(\R)$ is quantile $\monotonetransformation{\phi}$-unimodal if and only if $\mu$ is $\phi$-unimodal.
    \end{theorem}
    An outline of the proof is given in \autoref{Proof of quantile characterization of ɸ}, and the general \proofref{thm: equivilance between M concave and quantile M convex quasi arithmetic} can be found in \autoref{sec: Proof of general generalized inverse pair theorem on equivalenve of derivatives}.
    We calculate in \autoref{example: monotone transform power means} that $M_{\monotonetransformation{\phi_s}}=M_{\phi_{r}}$ for $r=-1-s$, $s\in\R$.
    This unifies several known characterizations of $\phi_s$-unimodal distributions.
    An absolutely continuous distribution is log-concave with a log-concave density $f_\mu$ if and only if $f_\mu\circ Q_\mu$ is concave \cite[Proposition A.3]{Bobkov.1996}.
    This condition is equivalent to quantile $\phi_{-1}$-unimodality of $\mu$, in line with $\monotonetransformation{\phi_0}=\phi_{-1}$, since $q_\mu\defeq\frac{1}{f_\mu\circ Q_\mu}$ is a valid quantile density in this case and $\phi_{-1}\circ q_\mu=-f_\mu\circ Q_\mu$.
    Under smoothness assumptions on the density, it was shown in \cite[Lemma E.1]{malamutEntropicOptimalTransport2025} that having $s$-concave densities is equivalent to quantile $\phi_r$-unimodality for $r=-s-1$, $s\in\R$.

\begin{figure}[ht]
        \centering
        \includegraphics[width=\textwidth]{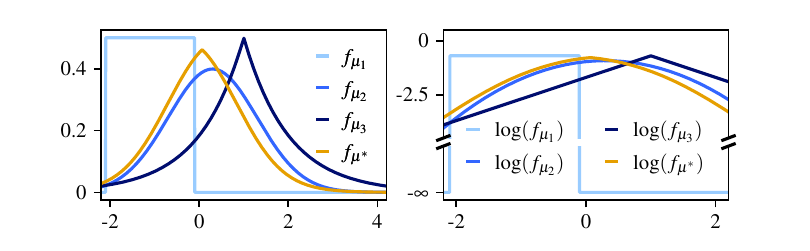}
        \caption{The density of the equal-weight {$2$-WB} $\mu^*$ (orange) of three log-concave (uniform, Gaussian, Laplace) distributions (blue) on the {l.h.s.\@} is again log-concave, as the {r.h.s.\@} shows.}
        \label{fig: log-concave plot}
    \end{figure}

    In contrast to $\phi$-unimodality, quantile $\phi$-unimodality does not assume absolute continuity of the measure $\mu$. Thus, we use the latter to extend the former and reformulate our main result.
    \begin{definition}
        \label{def: M concave general}
        Let $\phi$ be a generator. A distribution $\mu\in\distributionon(\R)$ is \emph{$\phi$-unimodal} if $\mu$ is quantile $\monotonetransformation{\phi}$-unimodal.
    \end{definition}
    \begin{theorem}
        \label{theo: Wasserstein 2 barycenter convexity of M convex funcitons}
        The set of univariate $\phi$-unimodal distributions is {$2$-WB} convex if and only if $M_{\monotonetransformation{\phi}}$ is concave.
    \end{theorem}
    The proof is given in \autoref{sec: Main theorem}. 
    \begin{corollary}
        \label{theo: Wasserstein 2 barycenter convexity of log-concave distributions}
        The set of log-concave distributions is {$2$-WB} convex.
    \end{corollary}
    \begin{proof}
        Log-concavity is quantile $\phi_{-1}$-unimodality, and the harmonic mean $M_{\phi_{-1}}$ is concave.
    \end{proof}
    \autoref{fig: log-concave plot} visualizes \autoref{theo: Wasserstein 2 barycenter convexity of log-concave distributions} for three absolutely continuous log-concave distributions. We include non-absolutely continuous distributions as we employ the {q.f.\@} perspective for our proofs. Restricting to absolutely continuous distributions afterwards comes at little additional cost.
    
    \begin{corollary}
        \label{cor: Wasserstein 2 barycenter convexity of absolutely continuous M convex functions}
        The set of absolutely continuous $\phi$-unimodal distributions is {$2$-WB} convex if and only if $M_{\monotonetransformation{\phi}}$ is concave.
    \end{corollary}
    \begin{proof}
        Indeed, the set of all absolutely continuous distributions is {$2$-WB} convex \cite{maAbsoluteContinuityWasserstein2025}.\footnote{They prove the stronger statement that if $\popmeasure\in\Wtwo(\Wtwo(\R))$ gives any mass to absolutely continuous distributions, $\popmeasure(\text{$\mu\ll\Lebesgue$})>0$, then the {$2$-WB} is absolutely continuous.} Thus, the intersection of all absolutely continuous and all $\phi$-unimodal distributions is {$2$-WB} convex if the latter is. The counterexample in the proof of \autoref{theo: Wasserstein 2 barycenter convexity of M convex funcitons} uses absolutely continuous distributions, so it also establishes the \enquote{only if} direction here.
    \end{proof}
\subsection{Prior work}
    The $2$-displacement convexity of the set of univariate log-concave distributions was a folklore result before being proved in \cite{matsudaShapeconstrainedDensityEstimation2026,malamutEntropicOptimalTransport2025}.
    In \cite[Proposition 10.8]{malamutEntropicOptimalTransport2025}, the author also proves the $2$-displacement convexity of absolutely continuous distributions with $s$-concave densities if and only if $s\geq-2$, relying on smoothness assumptions on the densities. The regime $s\geq-2$ is precisely the regime in which $M_{\phi_r}$ is concave, as it is equivalent to $-s-1=r\leq1$ \cite{chudziakWeightedQuasiarithmeticMeans2019}. This covers log-concavity, corresponding to $s=0$.
    Thus, in \autoref{theo: Wasserstein 2 barycenter convexity of M convex funcitons}, we unify and generalize these results by encompassing the larger class of $\phi$-unimodal families and showing {$2$-WB} convexity instead of $2$-displacement convexity. For log-concavity, the extension to {$2$-WB} convexity can also be shown directly from $2$-displacement convexity.
    In \autoref{sec: Relation of different notions of barycenter convexity results}, we explain how to extend $2$-displacement convexity to {$2$-WB} convexity under the additional assumption that the family is closed in $\Wtwo(\R)$. Without this additional assumption, the two are not equivalent in general (see \autoref{example: difference dispalcement convexity and barycenter convexity} for $p=2$). Such an extension is possible in the log-concave case.
\subsection{Other results on \texorpdfstring{{$2$-WB}\nocorr}{2-WB} convexity of families of unimodal distributions}
    Other examples of {$2$-WB} convex families of unimodal distributions are known in the literature. In \cite{backhoffStochasticGradientDescent2025}, the {$2$-WB} convexity of symmetric unimodal distributions has been shown.
    The proof does not use the quantile density approach developed here but rather works directly with the known closed-form expression for the {q.f.\@} of the {$2$-WB} in \eqref{eqn: characterization of Qunaitle fucntionas mean}.\footnote{They claim that their results also hold for $p\neq2$, but their argument is not valid because it assumes that the characterization of the {q.f.\@} as the pointwise mean \eqref{eqn: characterization of Qunaitle fucntionas mean} also holds for $p\neq 2$. \autoref{example: Counterexample Wasserstein p barycenter convexity for p neq 2} falsifies their claim.}
    They directly characterize all symmetric distributions and symmetric unimodal distributions via properties of $Q_\mu$ and show that the defining equalities and inequalities can be extended to the {q.f.\@} of the {$2$-WB}.
    From their proof one can directly obtain the {$2$-WB} convexity of unimodal distributions with a shared quantile mode $\beta\in[0,1]$, {i.e.}, of all distributions satisfying $Q_\mu$ being concave on $(0,\beta]$ and convex on $[\beta,1)$. In particular, the sets of all absolutely continuous distributions with non-increasing and of all with non-decreasing densities are each {$2$-WB} convex (compare also \cite[Proposition 3.2]{matsudaShapeconstrainedDensityEstimation2026}), as they correspond to a fixed quantile mode $\beta=0$ and $\beta=1$ under absolute continuity, respectively.
    \autoref{theo: Wasserstein 2 barycenter convexity of M convex funcitons} does not apply to unimodal distributions with a fixed quantile mode, but its proof can be modified to do so, as one again has to show that an average of convex functions is convex.
\subsection{Statistical applications}
    \autoref{cor: Wasserstein 2 barycenter convexity of absolutely continuous M convex functions} has a notable consequence for 
    shape-constrained nonparametric density estimation, see, {e.g.\@}, \cite{dumbgenMaximumLikelihoodEstimation2009} for log-concave, \cite{koenkerQuasiconcaveDensityEstimation2010} for $s$-concave, and \cite{sereginNonparametricEstimationMultivariate2010} for $M_\phi$-concave density estimation.
    Our result specifies when the {$2$-WBs} of such estimators share these properties.
    For example, for univariate log-concave distributions, no further estimation or projection onto the family of log-concave distributions is needed when applying the {$2$-WB} to average arbitrarily many log-concave density estimates.
    More generally, the same holds true for $s$- and $M_\phi$-concave densities with $s\geq-2$ and $M_{\monotonetransformation{\phi}}$ concave, respectively (see \autoref{cor: Wasserstein 2 barycenter convexity of absolutely continuous M convex functions}).
    In these cases, the {$2$-WB} automatically preserves the shape constraint and can be easily computed via its quantile representation in \eqref{eqn: characterization of Qunaitle fucntionas mean}. Hence, no additional, potentially computationally expensive steps are needed to enforce the shape constraint for the barycenter.
    This principle can be transferred to other shape constraints. Based on the preceding discussion, it holds for monotone densities or unimodal distributions with fixed quantile mode.

    A second application to shape-constrained density estimation has been put forward in \cite{matsudaShapeconstrainedDensityEstimation2026}. Instead of using maximum likelihood ideas for the estimator, they investigate a $W_2$-projection-based estimator. The empirical distribution $\mu_n$ of a sample $X_1,\dots,X_n\sim\mu\in\distributionon(\R)$ is projected onto the nearest distribution in 
    the desired set of shape-constrained distributions {w.r.t.\@} the $W_2$ distance.
    The shape constraints considered are non-increasing densities on $\R_{\geq0}$ as well as log-concave densities. A crucial step to establish the uniqueness of the projection estimator is to show the displacement convexity of these families. Thus, based on our results this estimation scheme can be extended to the {$2$-WB} convex $\phi$-unimodal families.
    When the target families are additionally $W_2$ closed, the projection is well-defined.
    As the considered families with shape-constrained densities are not $W_2$ closed, their closure, which contains, {e.g.\@}, Dirac distributions, is used instead.
    To yield a density estimator, one has to further show that the projection is absolutely continuous. This works for log-concave distributions as long as the sample does not concentrate on one point \cite[Proposition 4.6]{matsudaShapeconstrainedDensityEstimation2026}.
    We speculate that this result holds for $\phi$-unimodal distributions whose generators satisfy $\phi(0)=-\infty$, just like $\phi_0$ in the log-concave case.

\subsection{Structure of the paper}

We proceed as follows.
In \Cref{sec: non-decreasing functions mappings,sec: Random quantile functions and inverse measures,sec: Random inverse measures of quantile function and inverse measure of 2-WB}, we introduce the technical framework for the one-dimensional setting to establish the connection between random measures and random {q.f.s}{}.
We then characterize the absolute continuity of a {q.f.\@} in \autoref{sec: Quantile absolute continuity and quantile density selections}.
This allows us to prove \autoref{prop: Wasserstein 2 Barycenter of quantile absolutely continuous distributions on R have absolutely continuous quantile functions} on the {$2$-WB} convexity of quantile absolutely continuous distributions in \autoref{sec: 2-WB convexity of quantile absolute continuity}.

Afterwards, we return to $\phi$-unimodality and $M_\phi$-convexity in \Cref{sec: Quasi-arithmetic means,sec: Quasi-arithmetic mean-convexity}. In \autoref{Proof of quantile characterization of ɸ}, we outline the proof of the quantile characterization of $\phi$-unimodal distributions stated in \autoref{thm: equivilance between M concave and quantile M convex quasi arithmetic}.

\autoref{sec: Main theorem} contains the proof of our main result \autoref{theo: Wasserstein 2 barycenter convexity of M convex funcitons}.
We first discuss the existence of an $M_\phi$-convex quantile density selection in \autoref{sec: existence quantile density selection} and the concavity of quasi-arithmetic means in \autoref{sec: convexity of quasi-aritmetic means}.
With these results at hand, we give the proof in \autoref{sec: Proof of main theorem} and present the counterexamples ruling out extensions of \autoref{theo: Wasserstein 2 barycenter convexity of M convex funcitons} to $p\neq2$ or higher dimensions in \autoref{sec: Counterexamples for extensions of the main theorem}.

We finish with a discussion of the relations between the barycenter convexity properties mentioned here in \autoref{sec: Relation of different notions of barycenter convexity results} and the implications of our quantile characterization of generalized $\phi$-unimodality in \autoref{sec: discussion generalization of phi concavity}.

Proofs omitted from the main text can be found in the \hyperref[sec supp]{Appendix}.

\section[Quantile absolute continuity is \texorpdfstring{$2$}{2}-Wasserstein barycenter convexity]{Quantile absolute continuity is \texorpdfstring{$2$}{2}-Wasserstein barycenter\\ convexity}
 \label{sec: quantile absolute continuity is Wasserstein 2 barycenter convex}
 
We begin by providing a thorough discussion of the relation between distributions and {q.f.s}{}. Various technical details are provided in \autoref{sec app: Non-decreasing functions}.
\subsection{Equivalent characterizations of univariate measures}
    \label{sec: non-decreasing functions mappings}
    Recall the setting in \autoref{sec: the univariate setting}. For each distribution $\mu\in\distributionon(\R)$, we define its {c.d.f.\@} as $F_\mu(x)=\mu((-\infty,x])$ and its right-continuous {q.f.\@} as $Q_\mu(t)=\sup\preimage{F_\mu}{[0,t]}$, $t\in(0,1)$. Denote by $\CDFset(\R)$ the set of all {c.d.f.s} and by $\QFset(0,1)$ the set of all {q.f.s}{} of univariate distributions.
    This gives rise to the following bijective maps:
    \begin{align*}
         & \text{\emph{{c.d.f.\@} transform}}                                       & \CDFtransform                                               & :\distributionon(\R)\to\CDFset(\R),             & \CDFtransform(\mu)         & = F_\mu                                       \\
         & \text{with inverse \cite{klenkeProbabilityTheoryComprehensive2008}}
         & \Stieltjestransform                                                     & :\CDFset(\R)\to\distributionon(\R), & \Stieltjestransform(F_\mu)((a,b])                                              & =F_\mu(b)-F_\mu(a),                                                                \\
         & \text{\emph{{q.f.\@} transform}}                                        & \Quantiletransform                                          & :\distributionon(\R)\to\QFset(0,1),                 & \Quantiletransform(\mu)    & =Q_\mu                                        \\
         & \text{with inverse \cite{bobkovOnedimensionalEmpiricalMeasures2019}}   & \Transportmeasure                                           & :\QFset(0,1)\to\distributionon(\R),                 & \Transportmeasure(Q_\mu) & =\pushforward{{Q_\mu}}{\uniformmeasure},            \\
         & \text{\emph{generalized inverse}}                                      & \generalizedinverse                                         & :\QFset(0,1)\to\CDFset(\R), & \generalizedinverse(Q_\mu) & =F_\mu                                        \\
         & \text{with inverse \cite{delafortelleStudyGeneralizedInverses2015}}    & \generalizedinverse                                         & :\CDFset(\R)\to\QFset(0,1), & \generalizedinverse(F_\mu) & =Q_\mu.
        \intertext{
            We can further associate distributions on $\R$ with \emph{locally finite measures} on $(0,1)$, denoted by $\locallyfinitemeasureon(0,1)$, {i.e.}, measures $\mu$ such that $\mu([a,b])<+\infty$ for all $[a,b]\subset(0,1)$.
            To this end, we consider the Stieltjes measure associated with the {q.f.\@} $Q_\mu$ and define
        }
         & \text{via}                                                             & \Stieltjestransform                                          & :\QFset(0,1)\to\locallyfinitemeasureon(0,1),        & \Stieltjestransform(Q_\mu)((a,b]) & = Q_\mu(b)-Q_\mu(a)                                   \\
         & \text{or \cite[Prop. A.4]{bobkovOnedimensionalEmpiricalMeasures2019}} & \Transportmeasure                                           & :\CDFset(\R)\to\locallyfinitemeasureon(0,1),    & \Transportmeasure(F_\mu)       & = \restr{(\pushforward{{F_\mu}}{\Lebesgue})}{(0,1)} \\
         & \text{the \emph{inverse measure}}                                      & \InverseMeasure                                             & :\distributionon(\R)\to\locallyfinitemeasureon(0,1),                    & \InverseMeasure(\mu)       & =\mu^{-1}=\Stieltjestransform(Q_\mu).
    \end{align*}
    The latter mappings are not injective, as all translations of a distribution are mapped to the same inverse measure.
    Nevertheless, the mappings are surjective.
    All maps are continuous {w.r.t.\@} weak convergence of distributions, convergence at continuity points on $\QFset(0,1)$ and $\CDFset(\R)$, and vague convergence on $\locallyfinitemeasureon(0,1)$. Proofs are given in \autoref{sec app: Non-decreasing functions}.
    The resulting diagram commutes.
    \begin{figure}[!ht]
        \centering
        \hspace{-1cm}
        
\begin{tikzcd}[
        every label/.append style = {font = \normalsize},
        every cell/.append style = {font = \normalsize},
        column sep=large,
        row sep=large,
        cells={nodes={execute at begin node=\displaystyle}}
    ]
    \distributionon{\R}
    \arrow[rrr,twoheadrightarrow,"\InverseMeasure"]
    \arrow[d,leftharpoondown,"\CDFtransform"]
    &&& \locallyfinitemeasureon(0,1)
    \\
    \CDFset(\R)
    \arrow[rrr,leftharpoondown,"\generalizedinverse"]
    \arrow[u,leftharpoondown,"\CDFtransform^{-1}"]
    \arrow[urrr,twoheadrightarrow,"\Quantiletransform^{-1}",pos=.7]
    &&&
    \QFset(0,1)
    \arrow[ulll,rightharpoonup,"\Quantiletransform",pos=.15,swap,crossing over]
    \arrow[ulll,leftharpoondown,"\Quantiletransform^{-1}",pos=.75]
    \arrow[lll,leftharpoondown,"\generalizedinverse"]
    \arrow[u,twoheadrightarrow,"\CDFtransform^{-1}"]
\end{tikzcd}

        \hspace{-1cm}
        \caption{Commutative diagram between different representations of univariate distributions. Here, $\rightleftharpoons$ denotes bijective mappings and their inverses, while $\twoheadrightarrow$ denotes surjective maps.}
        \label{fig: Commutative diagram of univariate measures}
    \end{figure}
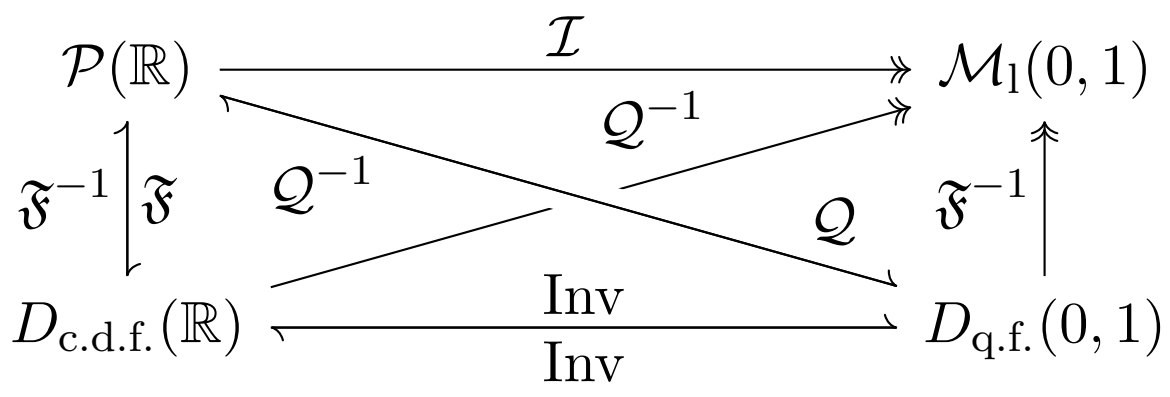

\subsection{Random quantile functions and inverse measures}
    \label{sec: Random quantile functions and inverse measures}
    
Our approach to describing random measures follows \cite{kallenbergRandomMeasuresTheory2017}.
We briefly review the required material.
Recall the ambient complete probability space $(\Omega,\A,\P)$.
\begin{definition}[Compare {\cite[Lemma 1.14]{kallenbergRandomMeasuresTheory2017}}]
    Let $\Y$ be a nonempty set and let $(\X,\B)$ be a function space $\X\subset\Rbar^{\Y}$
    equipped with the initial $\sigma$-algebra of evaluations at elements of $\Y$, {i.e.}, $\B=\sigma(\iota_{\Y})$ with $\iota_{\Y}\defeq\set{\iota_y:\X\to\Rbar\given y\in\Y}$.
    We call a measurable function $\randommeasure:(\Omega,\A,\P)\to(\X,\B)$ a \emph{random \enquote{element of $\X$}}. Equivalently, it is a kernel (from $\Omega$) to $\X$, {i.e.}, a function $\randommeasure:\Omega\times\Y\to\Rbar$ with the properties that $\randommeasure(\omega,\cdot)\defeq\randommeasure(\omega)\in\X$ and that $\randommeasure(\cdot,y)\defeq\iota_y\circ\randommeasure$ is $\A,\B(\Rbar)$-measurable.
    Then $\popmeasure\defeq\pushforward{\randommeasure}{\P}\in\distributionon(\X)\subset[0,1]^{\sigma(\iota_\Y)}$.
\end{definition}
\begin{table}[h]
    \centering
    \caption{Random distributions on $(\R,\B(\R))$ and their associated quantiles and measures.}
    \label{tab: Example random objects}
    \begin{tabular}{lrcl}
        \toprule random objects             & $\X$                           & $\subset$ & $\Rbar^{\Y}$              \\\midrule
        random distributions on             & $\distributionon(\R)$          & $\subset$ & $[0,1]^{\B(\R)}$          \\
        random {c.d.f.s} on                & $\CDFset(\R)$                  & $\subset$ & $[0,1]^\R$                \\
        random {q.f.s} on                  & $\QFset(0,1)$                  & $\subset$ & $\R^{(0,1)}$              \\
        random (locally finite) measures on & $\locallyfinitemeasureon(0,1)$ & $\subset$ & $\Rbar_{\geq0}^{\B(0,1)}$ \\
        \bottomrule
    \end{tabular}
\end{table}
In \autoref{tab: Example random objects}, we summarize random objects associated with random distributions that we use in the following.
Equipping these spaces $\X$ with the associated evaluation $\sigma$-algebras, all the maps in \autoref{fig: Commutative diagram of univariate measures} are measurable (see \autoref{lem: quantile transform homeomorphism Stieltjes transform continuous} in \autoref{sec: Continuity and measurability of the relations between measures and non-decreasing functions}), so we can switch between representations of random distributions.
\begin{lemma}
    \label{lem: equivalence of kernels}
    Let $\randommeasure:(\Omega,\A)\to\distributionon(\R)$. Then $\randommeasure$ is a random distribution if and only if $Q_{\randommeasure}\defeq\Quantiletransform\circ \randommeasure$ is a random {q.f}. In both cases, $\randommeasure^{-1}\defeq\InverseMeasure\circ\randommeasure$ is a random measure on $(0,1)$.
\end{lemma}
The \proofref{lem: equivalence of kernels} is in \autoref{sec: Continuity and measurability of the relations between measures and non-decreasing functions}.
We call $\randommeasure\sim\popmeasure\in\distributionon(\Wtwo(\R))$ a \emph{random finite-variance distribution} if $\randommeasure(\omega)\in\Wtwo(\R)$ for all $\omega\in\Omega$.
The isometry between $\Wtwo(\R)$ and $\QFset^2(0,1)\defeq\Lpmod^2(0,1)\cap\QFset(0,1)$ extends \autoref{lem: equivalence of kernels} to this setting. A random distribution $\randommeasure\sim\popmeasure\in\distributionon(\distributionon(\R))$ is a random finite-variance distribution if and only if $Q_{\randommeasure}(\omega)$ is square-integrable for each $\omega\in\Omega$.
On $\Wtwo(\R)$, the initial $\sigma$-algebra coincides with the metric Borel $\sigma$-algebra.
\begin{theorem}[{\cite[Lemma 2.4.2]{Panaretos.2020}}]
    \label{thm: measurable structure on Wasserstein space}
    It holds that $\B(\Wtwo(\R))=\sigma(\iota_{\B(\R)})$ on $\Wtwo(\R)$.
\end{theorem}
We give a {q.f.\@}-based proof of this special case of \cite[Lemma 2.4.2]{Panaretos.2020}, treating the general cases of $\Wp(\R)$ and $\QFset^{p}(0,1)$ in \autoref{proof: thm: measurable structure on Wasserstein space}.
\subsection{Random inverse measures and the inverse measure of the \texorpdfstring{{$2$-WB}\nocorr}{2-WB}}
    \label{sec: Random inverse measures of quantile function and inverse measure of 2-WB}
    We can now turn to random distributions with a well-defined {$2$-WB}, $\randommeasure\sim\popmeasure\in\Wtwo(\Wtwo(\R))$.
    From \cite[Theorem 3.2.11]{Panaretos.2020}, we directly see that $\randommeasure\sim\popmeasure\in\Wtwo(\Wtwo(\R))$ makes $Q_\randommeasure$ an \emph{integrable kernel}, {i.e.}, $Q_\randommeasure(t)=\iota_t\circ Q_\randommeasure$ is integrable for each $t\in(0,1)$, as $\Quantiletransform_{\mu^*}(t)=\mean[Q_\randommeasure(t)]$ for the {$2$-WB} $\mu^*=\Barycenterfunc{2}{\Wtwo}\,(\popmeasure)$.
    In this case, even $\mean[Q_\randommeasure(\cdot)]\in\QFset^2(0,1)$.
    The same holds for the random inverse measure. A random (locally finite) measure $\rho$ is \emph{integrable} if $\mean[\rho](B)\defeq\mean[\iota_B\circ\rho]=\mean[\rho(B)]$ defines a locally finite measure.
    \begin{theorem}
        \label{lem: random W_2 distribution gives an integrable inverse measure}
        Let $\randommeasure\sim\popmeasure\in\Wtwo(\Wtwo(\R))$. Then $\randommeasure^{-1}$ is an integrable random measure on $(0,1)$ and $\mean{\randommeasure^{-1}}=(\Barycenterfunc{2}{\Wtwo}\,(\popmeasure))^{-1}$.
    \end{theorem}
    \begin{proof}
        For $\randommeasure\sim\popmeasure\in\Wtwo(\Wtwo(\R))$, we have by \autoref{lem: equivalence of kernels} that $\randommeasure^{-1}$ is a random measure on $(0,1)$. We can directly calculate that
        \begin{align*}
            \mean[\randommeasure^{-1}]((a,b])
            \defeq{} & \mean{\randommeasure^{-1}((a,b])}
            =\mean{Q_\randommeasure(b)-Q_\randommeasure(a)}
            =\mean{Q_\randommeasure}(b)-\mean{Q_\randommeasure}(a)                                                                                               \\
            ={}      & Q_{\Barycenterfunc{2}{\Wtwo}\,(\popmeasure)}(b)-Q_{\Barycenterfunc{2}{\Wtwo}\,(\popmeasure)}(a)=(\Barycenterfunc{2}{\Wtwo}\,(\popmeasure))^{-1}((a,b])<\infty
        \end{align*}
        for all $0<a<b<1$. Thus, $\mean_{\randommeasure\sim\popmeasure}[\randommeasure^{-1}]$ is locally finite, as any compact interval is contained in such a half-open interval.
        Therefore, $\randommeasure^{-1}$ is integrable and $\mean[\randommeasure^{-1}]=(\Barycenterfunc{2}{\Wtwo}\,(\popmeasure))^{-1}$.
    \end{proof}

\subsection{Absolute continuity and density selections}
    \label{sec: Quantile absolute continuity and quantile density selections}
    
    A {c.d.f.\@} $F_\mu$ is absolutely continuous if and only if $\mu\ll\Lebesgue$, {i.e.}, $\mu$ is absolutely continuous {w.r.t.\@} the Lebesgue measure on $\R$.
    Similarly, a {q.f.\@} $Q_\mu$ is locally absolutely continuous if and only if $\mu^{-1}\ll\uniformmeasure$, in which case we call $\mu$ \emph{quantile absolutely continuous}.
    We then call any Radon-Nikodym derivative $f_{\mu}\defeq\RadonNikodym{F_\mu}{\Lebesgue}\defeq\RadonNikodym{\mu}{\Lebesgue}$ a \emph{density} of $\mu$ and any Radon-Nikodym derivative $q_\mu\defeq\RadonNikodym{Q_\mu}{\uniformmeasure}\defeq\RadonNikodym{\mu^{-1}}{\uniformmeasure}$ a \emph{quantile density} of $\mu$, both being unique almost everywhere with respect to $\Lebesgue$ and $\uniformmeasure$, respectively.
    Analogously to $\Quantiletransform$ and $\CDFtransform$, we can define a (quantile) density selection for (quantile) absolutely continuous distributions.
    Similarly, we call a random measure $\rho\sim\popmeasure\in\distributionon(\locallyfinitemeasureon(I))$ \emph{absolutely continuous} if $\P(\rho\ll\restr{\Lebesgue}{I})=1$ and a random distribution $\randommeasure$ \emph{quantile absolutely continuous} if $\randommeasure^{-1}$ is absolutely continuous.\footnote{The set of absolutely continuous distributions or absolutely continuous locally finite measures is measurable \cite{kallenbergRandomMeasuresTheory2017}.}
    To define the associated random quantile density, we need a quantile density selection.
    \begin{lemma}
        \label{lem: existence of quantile density selections}
        There exists a \emph{quantile density selection} $\QuantDensSelec:\distributionon(\R)\times (0,1)\to\Rbar_{\geq0}$ that satisfies
        \begin{enumerate}
            \item $\QuantDensSelec$ is $\sigma(\iota_{\B(\R)})\otimes\B((0,1)),\B(\Rbar_{\geq0})$-measurable and
            \item $\QuantDensSelec(\mu,\cdot)=\RadonNikodym{\mu^{-1}}{\restr{\Lebesgue}{(0,1)}}$ $\uniformmeasure$-{a.e.\@} for all quantile absolutely continuous $\mu\in\distributionon(\R)$.
        \end{enumerate}
    \end{lemma}
    The \proofref{lem: existence of quantile density selections}, based on \cite{kallenbergRandomMeasuresTheory2017}, is given in \autoref{sec: appendix density selection}.
    Then $\QuantDensSelec(\randommeasure,\cdot)$ is a corresponding random density.
    \begin{remark}
        A quantile density selection can be understood as a special kernel from $(\distributionon(\R),\sigma(\iota_{\B(\R)}))$ to $\Rbar_{\geq0}^{(0,1)}$. Joint measurability implies, in particular, that the selection $\QuantDensSelec(\cdot,t)$ is measurable for each $t\in (0,1)$.
    \end{remark}
    \subsection{\texorpdfstring{$2$}{2}-WB convexity of quantile absolute continuity}
    \label{sec: 2-WB convexity of quantile absolute continuity}
    The joint measurability of $\QuantDensSelec$ enables the proof of \autoref{prop: Wasserstein 2 Barycenter of quantile absolutely continuous distributions on R have absolutely continuous quantile functions}, which establishes the {$2$-WB} convexity of the set of quantile absolutely continuous distributions and a closed-form expression for the quantile density of the {$2$-WB}.
    \begin{delayedproof}{prop: Wasserstein 2 Barycenter of quantile absolutely continuous distributions on R have absolutely continuous quantile functions}
        For a random measure $\randommeasure\sim\popmeasure\in\Wtwo(\Wtwo(\R))$, we have by \autoref{lem: random W_2 distribution gives an integrable inverse measure} that $\randommeasure^{-1}$ is an integrable random measure and $\mean[\randommeasure^{-1}]=(\Barycenterfunc{2}{\Wtwo}\,(\popmeasure))^{-1}\eqdef(\mu^*)^{-1}$. If, further, $\popmeasure$ is supported on quantile absolutely continuous distributions, {i.e.}, $\P(\randommeasure^{-1}\ll\uniformmeasure)=1$, then we can apply \cite[2.12]{kallenbergRandomMeasuresTheory2017}, which establishes the absolute continuity of the expectation of an absolutely continuous random measure, to obtain that $(\mu^*)^{-1}\ll\uniformmeasure$, so $Q_{\mu^*}$ is locally absolutely continuous.
        By \autoref{lem: existence of quantile density selections}, there exists a jointly measurable quantile density selection $\QuantDensSelec$. Then $q_{\mu^*}\defeq \mean[\QuantDensSelec(\randommeasure,\cdot)]$ is measurable since $\QuantDensSelec$ is also nonnegative, and we can apply Tonelli's theorem to calculate, for arbitrary $A\in\B(0,1)$, that
        \begin{align*}
            (\mu^{*})^{-1}(A)
             &
            =\int_{\Omega} (\randommeasure^{-1}(\omega))(A)\diff \P(\omega)
            =\int_{\Omega} \int_{A}\QuantDensSelec(\randommeasure(\omega),t)\diff\uniformmeasure(t)\diff\P(\omega)    \\
             & =\int_{A}\int_{\Omega} \QuantDensSelec(\randommeasure(\omega),t)\diff\P(\omega)\diff\uniformmeasure(t)
            =\int_{A} q_{\mu^*}(t)\diff\uniformmeasure(t).
        \end{align*}
        Therefore, $q_{\mu^*}=\RadonNikodym{(\mu^*)^{-1}}{\uniformmeasure}$ $\uniformmeasure$-{a.e.\@}, and $q_{\mu^*}$ is a quantile density of $\mu^*$.
    \end{delayedproof}

\section{Different unimodality concepts}
\label{sec: Quantile Characterization of Ms concave distributions}

The general definition of univariate unimodal distributions also encompasses non-absolutely continuous distributions.
A distribution $\mu\in\distributionon(\R)$ is \emph{unimodal} if there exists a mode $m\in\R$ such that its {c.d.f.\@} $F_\mu$ is convex on $(-\infty,m]$ and concave on $[m,\infty)$ \cite[Definition 1.1]{Dharmadhikari.1988}. If we use quasi-concave densities as our definition, we restrict ourselves to the non-weakly closed set of absolutely continuous distributions. Such a restriction is not necessary from the quantile perspective. A distribution $\mu\in\distributionon(\R)$ is unimodal if and only if there exists a \emph{quantile mode} $\beta\in[0,1]$ such that its {q.f.\@} $Q_\mu$ is concave on $(0,\beta]$ and convex on $[\beta,1)$ \cite[Theorem 1.6]{zobelmarkusQuantileCharacterizationUnivariate2026arxiv}.
This condition is strong enough to ensure the local absolute continuity of the {q.f.\@} of a unimodal distribution, and thus of the {$2$-WB} of unimodal distributions by \autoref{prop: Wasserstein 2 Barycenter of quantile absolutely continuous distributions on R have absolutely continuous quantile functions}, but not strong enough to also enforce the quasi-convexity of the quantile density of the {$2$-WB} as in \autoref{example: unimodal distributions are not 2 Wasserstein barycenter convex}.

The restrictions on both the density and the quantile density can be tightened by enforcing stronger concavity or convexity conditions, respectively.
\subsection{Mean-convexity}
    \label{sec: Quasi-arithmetic means}
    A function $M:I^2\to I$, $I\subset\R$ an open, nonvoid interval, is a \emph{mean} if it satisfies the inequalities $\min(x,y)\leq M(x,y)\leq\max(x,y)$ for all $x,y\in I$ and a \emph{weighted mean} if $M:I^2\times(0,1)\to I$ is a mean for each $\convexindex\in(0,1)$, {i.e.}, $M_\convexindex=M(\cdot,\cdot,\convexindex)$ is a mean.
    For two weighted means $M,N$ on intervals $I,J$, respectively, \cite[2.6]{niculescuConvexFunctionsTheir2006} calls a function $f:I\to J$ \emph{$(M,N)$-convex (-concave)} if
    $
        f(M_{\convexindex}(x,y)){\leq}(\geq) N_{\convexindex}(f(x),f(y))$, for all $x,y\in I$, $\convexindex\in(0,1)
    $.
    With this, we recover quasi-convexity (-concavity) as $(M,N)$-convexity (-concavity) with the choices of $M=A$, the arithmetic mean, $A_{\convexindex}(x,y)=\convexindex x+(1-\convexindex) y$, $\convexindex\in(0,1)$, and $N=\max$ ($N=\min$), where we interpret $\min$ and $\max$ as weighted means constant in the weight, {i.e.}, $\min_{\convexindex}=\min$. By the definition of a mean $M$, any $(A,M)$-convex (-concave) function is quasi-convex (-concave).
\subsection{Quasi-arithmetic mean-convexity}
    \label{sec: Quasi-arithmetic mean-convexity}
    Recall the definition of generators $\phi:\R_{>0}\to\R$ in \autoref{def: quantile phi-convex distribution} \ref{enum: def: generator} and the corresponding quasi-arithmetic means on $\R_{>0}$, given by $M_\phi:\R_{>0}^2\times(0,1)\to\R_{>0}$
    \begin{equation}
        \label{def: M_phi mean}
        M_{\phi,\convexindex}(\theta)=M_{\phi}(\theta,\convexindex)=\phi^{-1}(\convexindex\phi(\theta_1)+(1-\convexindex)\phi(\theta_2)),\quad \theta\in\R^2_{>0},\convexindex\in(0,1).
    \end{equation}
    To define $(A,M_\phi)$-convex (-concave) densities, we need to treat the boundary values $0$ and $+\infty$. Thus, we define the set
    \begin{equation*}
        C_\phi\defeq\set[\Big]{\theta\in\Rbar^2_{\geq0}\given \text{$\lim_{\vartheta\to \theta,\vartheta\in\R_{>0}^2} M_{\phi,\frac12}(\vartheta)$ exists in $\Rbar_{\geq0}$}}
    \end{equation*}
    and define $M_{\phi,\convexindex}(\theta)\defeq\lim_{\vartheta\to \theta,\vartheta\in\R_{>0}^2} M_{\phi,\convexindex}(\vartheta)$ for each $\convexindex\in(0,1)$, $\theta\in C_\phi$.
    With $\phi(0)\defeq\lim_{x\searrow0} \phi(x)$ and $\phi(+\infty)\defeq\lim_{x\nearrow+\infty}\phi(x)$, we have that \eqref{def: M_phi mean} also holds for all $\theta\in C_\phi$, and thus $M_{\phi,\convexindex}:C_\phi\to\Rbar_{\geq0}$ is continuous for all $\convexindex\in(0,1)$. We have $C_\phi=\Rbar_{\geq0}^2$ unless both $\phi(0)=-\infty$ and $\phi(+\infty)=+\infty$, in which case $C_\phi=\Rbar_{\geq0}^2\setminus\set{(0,+\infty),(+\infty,0)}$. Thus, in all cases, $C_\phi$ is convex.
    \begin{definition}[Replaces \autoref{def: quantile phi-convex distribution} \ref{enum: M_phi convex}]
        Let $\phi$ be a generator. A function $g:I\to\Rbar_{\geq0}$ is \emph{$M_\phi$-well-defined} if $(g(x),g(y))\in C_\phi$ for all $x,y\in I$, and \emph{$M_\phi$-convex (-concave)} if it is $M_\phi$-well-defined and $g$ is $(A,M_\phi)$-convex (-concave).
    \end{definition}
    The refined definition of $M_\phi$-convex and $M_\phi$-concave functions does not change the resulting definitions of $\phi$-unimodal and quantile $\phi$-unimodal distributions as the well-definedness was implied in the former definitions.
    \begin{lemma}[Aczél-correspondence {\cite[2.6.1. Lemma]{niculescuConvexFunctionsTheir2006}}]
        \label{lem: aczel correspondence}
        Let $I$ be an open, nonvoid interval and $\phi$ a generator. A non-constant function $g:I\to\Rbar_{\geq0}$ is $M_{\phi}$-convex (-concave) if and only if $\phi\circ g:I\to\Rbar$ is convex (concave).
    \end{lemma}
    Analogously, we call a function $g$ $M_\phi$-affine if $\phi\circ g$ is affine. The \proofref{lem: aczel correspondence} is in \autoref{sec: Proof of aczel correspondence}. 
    \begin{remark}
        \label{rem: Aczel correspondence}
        \begin{enumerate}
            \item
                  The restriction to non-constant functions is necessary. The constant quantile density $q_\mu\equiv0$ corresponding to a Dirac or degenerate distribution $\delta_{m}$, $m\in\R$, is $M_\phi$-convex for every generator $\phi$, just like any constant function. On the other hand, if $\phi(0)=-\infty$, then $\phi\circ q_\mu \equiv-\infty$ is not convex under the standard definition of convexity, which excludes the value $-\infty$ from the codomain.
            \item
                  The notion of $M_\phi$-convexity here is not to be confused with $M$-convexity in \cite{matkowskiConvexFunctionsRespect1997}. While we mean $(A,M_\phi)$-convexity, the latter means $(M,M)$-convexity.
            \item
                  \label{enum: rem: Equivalence class of quasi-arithmetic means}
                  Quasi-arithmetic means for different generators $\phi,\varphi$ are the same exactly when there are real numbers $a,b\in\R$, $a>0$, such that $\phi=a\varphi+b$ \cite[Remark 2]{chudziakWeightedQuasiarithmeticMeans2019}. In the literature \cite{hardyInequalities1952}, strictly decreasing continuous functions are also taken as generators. Then $a<0$ would be possible. For example, both $\frac{1}{x}$ and $-\frac{1}{x}$ generate the harmonic mean, but only $-\frac{1}{x}$ is increasing. For a strictly decreasing generator $\varphi$, \autoref{lem: aczel correspondence} would equate $M_{\varphi}$-convexity of $f$ with the concavity of $\varphi\circ f$, and $M_{\varphi}$-concavity with the convexity of the composition, respectively.
            \item
                  Functions $f$ for which there exists a generator $\phi$ such that $\phi\circ f$ is convex (concave) have been studied independently. They are called convex (concave) transformable or convexifiable functions and date back to \cite{avrielGeneralizedConvexFunctions1974,fenchelConvexConesSets1953}. See also \cite{horstConvexificationNonlinearProgramming1984,avrielGeneralizedConcavity2010}.
        \end{enumerate}
    \end{remark}

\subsection{Proof of quantile characterization of \texorpdfstring{$\phi$}{ɸ}-unimodality}
    \label{Proof of quantile characterization of ɸ}
    
We will now prove \autoref{thm: equivilance between M concave and quantile M convex quasi arithmetic} by establishing the relationship between $\phi$-unimodality and quantile $\phi$-unimodality via co-generators.
\begin{definition}
    \label{def: monoton transformation}
    Let $\phi:\R_{>0}\to\R$ be a generator with inverse $\phi^{-1}$, and fix $c=\frac12$. Define the \emph{co-generator} as
    \begin{equation*}
        \monotonetransformation{\phi}(x)
        \defeq -\int_{\phi(\frac{1}{c})}^{\phi(\frac{1}{x})}\phi^{-1}\dd\Lebesgue +\phi(c)=\integratedgeneralizedinverse{\phi}^{\phi(\frac1c),\phi(c)}\circ\phi\circ\frac{1}{\Id}(x),\quad x\in\R_{>0},
    \end{equation*}
    where $\integratedgeneralizedinverse{\phi}^{y,z}(x)
        \defeq-\int_{y}^x\phi^{-1}\diff\Lebesgue +z$ for all $x\in\Functionimage{\phi}{\R_{>0}}$ and constants $y\in\Functionimage{\phi}{\R_{>0}}$, $z\in\R$.
\end{definition}
The function $\integratedgeneralizedinverse{\phi}^{y,z}$ is a locally absolutely continuous, strictly decreasing antiderivative of $-\phi^{-1}$. The resulting function $\monotonetransformation{\phi}$ is again a generator. It is continuous as the composition of continuous functions. The composition of strictly monotone functions is again strictly monotone. It is strictly increasing if an even number of the functions in the composition are strictly decreasing and strictly decreasing if an odd number of these functions are strictly decreasing. Therefore, $\monotonetransformation{\phi}$ is strictly increasing.
The particular choice of the constants $y,z$ ensures that $\monotonetransformation{\phi}(c)=\phi(c)$ and that the mapping is an involution. The same would hold for different choices of $c$, but we only need one such mapping. In the following, we will abbreviate $\integratedgeneralizedinverse{\phi}\defeq\integratedgeneralizedinverse{\phi}^{\phi(\frac1c),\phi(c)}$.
\begin{lemma}
    \label{lem: Bijectivity of Integrated inverse Transform}
    Let $\phi$ be a generator. Then $\monotonetransformation{(\monotonetransformation{\phi})}=\phi$.
\end{lemma}
The \proofref{lem: Bijectivity of Integrated inverse Transform} is in \autoref{proof: lem: Bijectivity of Integrated inverse Transform}.
The definition of $\monotonetransformation{(\cdot)}$ is based on the following equation
\begin{equation}
    \label{eqn: connection between the two phi transforms}
    \monotonetransformation{\phi}\circ q_\mu
    = \integratedgeneralizedinverse{\phi}\circ \phi\circ\frac{1}{\Id}\parenthesis*{\frac1{f_\mu\circ Q_\mu}}
    = \integratedgeneralizedinverse{\phi}\circ (\phi\circ f_\mu)\circ Q_\mu.
\end{equation}
The equality in \eqref{eqn: connection between the two phi transforms} holds as soon as $f_\mu$ and $q_\mu$ are chosen to satisfy $q_\mu=\frac{1}{f_\mu\circ Q_\mu}$, which is always possible if both $F_\mu$ and $Q_\mu$ are locally absolutely continuous \cite[Proposition A.18]{bobkovOnedimensionalEmpiricalMeasures2019}.
This allows us to prove the relation between $M_\phi$-concave densities and $M_{\monotonetransformation{\phi}}$-convex quantile densities stated in \autoref{thm: equivilance between M concave and quantile M convex quasi arithmetic} by using the Aczél correspondence (\autoref{lem: aczel correspondence}).
\begin{proof}[Proof (sketch) of \autoref{thm: equivilance between M concave and quantile M convex quasi arithmetic}]
    Assume for simplicity that all derivatives exist and are positive and that $\support(\mu)=\R$. By the absolute continuity of $\mu$, the {c.d.f.\@} $F_\mu=Q_\mu^{-1}$ is the inverse of the {q.f}. Denote by $f_\mu\defeq F_\mu'$, $q_\mu=Q_\mu'$ the classical local derivatives of $F_\mu$ and $Q_\mu$, respectively. We have $f_\mu=\frac{1}{q_\mu\circ F_\mu}$ and $q_\mu=\frac{1}{f_\mu\circ Q_\mu}$. Therefore, \eqref{eqn: connection between the two phi transforms} holds, and we can differentiate it to obtain
    \begin{equation*}
        (\monotonetransformation{\phi}\circ q_\mu)'
        =(\integratedgeneralizedinverse{\phi}\circ (\phi\circ f_\mu)\circ Q_\mu)'
        =(-\phi^{-1}\circ\phi\circ f_\mu\circ Q_\mu)\cdot(\phi\circ f_\mu)'\circ Q_\mu\cdot q_\mu
        =-(\phi\circ f_\mu)'\circ Q_\mu.
    \end{equation*}
    Therefore, $(\phi\circ f_\mu)'$ is non-increasing if and only if $(\monotonetransformation{\phi}\circ q_\mu)'$ is non-decreasing, and thus $\phi\circ f_\mu$ is concave if and only if $\monotonetransformation{\phi}\circ q_\mu$ is convex.
    For the general \proofref{thm: equivilance between M concave and quantile M convex quasi arithmetic} of \autoref{thm: equivilance between M concave and quantile M convex quasi arithmetic}, see \autoref{sec: Proof of general generalized inverse pair theorem on equivalenve of derivatives}.
\end{proof}
We can also interpret $\monotonetransformation{(\cdot)}$ as a map between quasi-arithmetic means, since equivalence classes of generators (see \autoref{rem: Aczel correspondence} \ref{enum: rem: Equivalence class of quasi-arithmetic means}) are mapped to equivalence classes of generators, {i.e.},
\begin{equation*}
    \monotonetransformation{(a\phi+b)}=a\monotonetransformation{\phi} + b,\quad a\in\R_{>0},b\in\R.
\end{equation*}
\begin{example}[{Application to power means}]
    \label{example: monotone transform power means}
    Denote by $M_s\defeq M_{\phi_s}$ the power or Hölder mean for $s\in\R$. Then the associated mean is again a power mean, as
    \begin{align*}
        \monotonetransformation{\phi_s}(x)
         & =-\int_{\phi_s(\frac1c)}^{\phi_s(\frac1x)}\abs{z}^{\frac1s}\diff\Lebesgue(z)+\phi_s(c)=-\sgn(s)\frac{s}{s+1}\parenthesis*{\abs*{\phi_s\parenthesis*{\frac1x}}^{\frac{s+1}{s}} -\abs*{\phi_s\parenthesis*{\frac1c}}^{\frac{s+1}{s}}}+ \phi_s(c)              \\
         & \eqdef\frac{\abs{s}}{\abs{s+1}}\sgn(-s-1)x^{-s-1}+b(s,c)\eqdef a(s)\phi_{-s-1}(x)+b(s,c),\quad x\in\R_{>0}.
    \end{align*}
    Therefore, $\monotonetransformation{M_s}\defeq M_{\monotonetransformation{\phi_s}}=M_{r}$ for $r=-s-1$ and $s\notin\set{-1,0}$.
    In the log-concave case $(s=0)$, we have
    \begin{equation*}
        \monotonetransformation{\phi_{0}}(x)
        =-\exp*{\log*{\frac{1}{x}}}+\frac{1}{c}+\log(c)
        \eqdef a(0)\phi_{-1}(x)+b(0,c),\quad x\in\R_{>0}.
    \end{equation*}
    Similarly, for the harmonic case $(s=-1)$, we have
    \begin{equation*}
        \monotonetransformation{\phi_{-1}}(x)
        =\log*{\abs*{-\frac{1}{\frac1x}}}-\log*{c}-\frac{1}{c}
        \eqdef a(-1)\phi_{0}(x)+b(-1,c),\quad x\in\R_{>0}.
    \end{equation*}
\end{example}
The calculations in \autoref{example: monotone transform power means} show that \autoref{thm: equivilance between M concave and quantile M convex quasi arithmetic} generalizes the quantile characterization of distributions with $s$-concave densities in \cite[Lemma E.1]{malamutEntropicOptimalTransport2025}, as we do not require any smoothness assumptions on $f_\mu$.

The set of quantile $\monotonetransformation{\phi}$-unimodal distributions includes all absolutely continuous $\phi$-unimodal distributions and thus serves as a natural generalization of $\phi$-unimodal distributions to non-absolutely continuous distributions. We postpone the discussion of the independent interest of this generalization to \autoref{sec: discussion generalization of phi concavity} and here only characterize generalized $\phi$-unimodal distributions from \autoref{def: M concave general}.
\begin{theorem}
    \label{thm: characterization of not absolutely continuous M concave distribution}
    Let $\phi$ be a generator.
    A distribution $\mu\in\distributionon(\R)$ is $\phi$-unimodal if and only if there exists a decomposition $\mu=\convexindex\mu_{a}+(1-\convexindex)\delta_{m}$ for some $\convexindex\in[0,1]$, $m\in\support{\mu_a}$, $\mu_a\ll\Lebesgue$, and there exists a density $f_a$ of $\mu_a$ such that
    \begin{enumerate}
        \item $\convexindex f_a$ is $M_\phi$-concave on $\interior{\support{\mu_a}}$ with mode $m$ and
        \item $\lim\limits_{\substack{x\to m,x\in\interior{\support{\mu_a}}}}f_a(x)=+\infty$ if $\convexindex\in(0,1)$.
    \end{enumerate}
\end{theorem}
For a \proofref{thm: characterization of not absolutely continuous M concave distribution} see \autoref{sec: proof of chararacterization of phi unimodality}.

\section{Main theorem}
\label{sec: Main theorem}
 
We can now return to proving the {$2$-WB} convexity of $\phi$-unimodal distributions. To apply \autoref{prop: Wasserstein 2 Barycenter of quantile absolutely continuous distributions on R have absolutely continuous quantile functions}, we first establish the existence of a measurable selection of  $M_\phi$-convex quantile densities in \autoref{lem: existence of phi-convex modification of quantile density}. Afterwards, we recall the convexity of quasi-arithmetic means. We then prove \autoref{theo: Wasserstein 2 barycenter convexity of M convex funcitons} and finally discuss counterexamples for $d\neq1$ and $p\neq2$.

\subsection{Measurable selection of \texorpdfstring{$M_\phi$}{Mᵩ}-convex quantile densities}
    \label{sec: existence quantile density selection}
    
We establish the existence of a quantile density selection that selects lower semicontinuous versions of quantile densities whenever possible in \autoref{sec: proof of density selection for quantile densities}.
\begin{lemma}
 \label{lem: existence of lower semicontinuous modification of quantile density}
    There exists a quantile density selection $\tilde{\QuantDensSelec}:\distributionon(\R)\times(0,1)\to\Rbar_{\geq0}$ such that, for each $\mu\in\distributionon(\R)$, the function $\tilde{\QuantDensSelec}(\mu,\cdot):(0,1)\to\Rbar_{\geq0}$ coincides with its own greatest lower semicontinuous {a.e.\@} minorant if and only if there exists a lower semicontinuous version of $\tilde{\QuantDensSelec}(\mu,\cdot)$.
\end{lemma}
This selection already selects $M_\phi$-convex quantile densities whenever possible, {i.e.}, it suffices for \autoref{lem: existence of phi-convex modification of quantile density}.
\begin{delayedproof}{lem: existence of phi-convex modification of quantile density}
    We establish the equivalence for the selection $\tilde{\QuantDensSelec}$ from \autoref{lem: existence of lower semicontinuous modification of quantile density}.
    \begin{itemize}
        \item[\enquote{$\Rightarrow$}:]
        Let $q_\mu$ be the $M_\phi$-convex quantile density of $\mu$. Then $q_\mu$ is continuous: Either $\phi\circ q_\mu$ is convex, and thus $q_\mu=\phi^{-1}\circ(\phi\circ q_\mu)$ is continuous, or $q_\mu\equiv0$. Therefore, $q_\mu$ is, in particular, lower semicontinuous. As $\mu$ is quantile absolutely continuous, we have $q_\mu=\tilde{\QuantDensSelec}(\mu,\cdot)$ $\uniformmeasure$-{a.e}. Thus, $\tilde{\QuantDensSelec}(\mu,\cdot)$ has a lower semicontinuous version. By \autoref{lem: existence of lower semicontinuous modification of quantile density}, we therefore have $q_\mu=\tilde{\QuantDensSelec}(\mu,\cdot)$. Thus, $\tilde{\QuantDensSelec}(\mu,\cdot)$ is $M_\phi$-convex if $\mu$ is quantile $\phi$-unimodal.
        \item[\enquote{$\Leftarrow$}:] For every quantile absolutely continuous distribution $\mu$, we have that $\tilde{\QuantDensSelec}(\mu,\cdot)$ is a quantile density of $\mu$. Thus, if $\tilde{\QuantDensSelec}(\mu,\cdot)$ is $M_\phi$-convex, then $\mu$ is quantile $\phi$-unimodal.\qedhere
    \end{itemize}
\end{delayedproof}

 \subsection{Convexity of quasi-arithmetic means}
 \label{sec: convexity of quasi-aritmetic means}

The following characterizations of the concavity of quasi-arithmetic means are from \cite{chudziakWeightedQuasiarithmeticMeans2019}.
\begin{definition}
    We call a quasi-arithmetic mean $M_\phi$ \emph{convex (concave)} if, for each $\convexindex\in(0,1)$, the weighted mean $M_{\phi,\convexindex}$ is convex (concave) as a function from $\R_{>0}^2$ to $\R_{>0}$.
\end{definition}
\begin{lemma}[{\cite[Theorem 1]{chudziakWeightedQuasiarithmeticMeans2019}}]
    \label{lem Convexity of mean}
    A quasi-arithmetic mean $M_\phi$ is convex (concave) if and only if $M_{\phi,\frac12}$ is Jensen convex (concave), that is,
    \begin{equation*}
        M_{\phi,\frac12}\parenthesis*{\frac{\theta+\vartheta}{2}}\overset{(\geq)}{\leq}\frac12(M_{\phi,\frac12}(\theta) + M_{\phi,\frac12}(\vartheta)),\quad \theta,\vartheta\in\R_{>0}^2.
    \end{equation*}
\end{lemma}
\begin{proof}
    The fixed choice $\convexindex=\frac12$ is actually arbitrary. By \cite[Theorem 1]{chudziakWeightedQuasiarithmeticMeans2019}, the convexity (concavity) of $M_{\phi,\convexindex}$ for one $\convexindex\in(0,1)$ is equivalent to its convexity (concavity) for all $\convexindex\in(0,1)$. Furthermore, as quasi-arithmetic means are locally bounded, Jensen convexity (concavity) already implies convexity (concavity) \cite[Remark 1]{chudziakWeightedQuasiarithmeticMeans2019}.
\end{proof}
While convexity of $M_\phi$ implies convexity of $\phi$ \cite[Theorem 3]{chudziakWeightedQuasiarithmeticMeans2019}, the reverse does not hold \cite[Example 2]{chudziakWeightedQuasiarithmeticMeans2019}.
If we extend $M_\phi$ to $C_\phi$, then this extended function is still convex (concave) on $C_\phi$ if $M_\phi$ was convex (concave) on $\R_{>0}$, as $C_\phi$ is convex and $M_{\phi}$ is extended continuously. We call a measurable function $h:\Omega\to\Rbar_{\geq0}^2$ \emph{$M_\phi$-well-defined} if $\Functionimage{h}{\Omega}\subset C_\phi$.
\begin{theorem}[Generalized Jensen inequality]
    \label{theo:generalized Jensen inequality}
    Let $(\Omega,\A,\P)$ be a probability space, let $h:\Omega\to\Rbar^2_{\geq0}$ be a measurable and $M_\phi$-well-defined function such that $\int h \diff\P\in C_\phi\subset\Rbar_{\geq0}^2$, and let $M_\phi$ be a convex (concave) quasi-arithmetic mean on $\R_{>0}$. Then
    \begin{equation*}                                                                M_{\phi,\convexindex}\parenthesis*{\int_\Omega h\diff\P}
        \overset{(\geq)}{\leq}\int_\Omega M_{\phi,\convexindex}\parenthesis*{h}\diff\P,\quad\convexindex\in(0,1).
    \end{equation*}
\end{theorem}
The \proofref{theo:generalized Jensen inequality} of this natural result is included in \autoref{sec: proof of Jensens inequality} for completeness.

\subsection{Proof of \autoref{theo: Wasserstein 2 barycenter convexity of M convex funcitons}}
\label{sec: Proof of main theorem}
    
Equipped with the appropriate quantile density selection and the notion of concavity of quasi-arithmetic means, we now prove our main theorem.
The proof of the reverse implication uses a construction similar to the one in the proof of \cite[Proposition 10.8]{malamutEntropicOptimalTransport2025}.
\begin{delayedproof}{theo: Wasserstein 2 barycenter convexity of M convex funcitons}
    \noindent
    \begin{itemize}    
        \item[\enquote{$\Leftarrow$}:]
              Let $\randommeasure\sim\popmeasure\in\Wtwo(\Wtwo(\R))$ be supported on $\phi$-unimodal distributions, where $M_{\monotonetransformation{\phi}}$ is concave, and denote by $\mu^*$ the {$2$-WB} of $\popmeasure$. By \autoref{prop: Wasserstein 2 Barycenter of quantile absolutely continuous distributions on R have absolutely continuous quantile functions}, we have that $q_{\mu^*}\defeq\mean[\QuantDensSelec(\randommeasure,\cdot)]$ is a valid quantile density of $\mu^*$ for $\QuantDensSelec$ from \autoref{lem: existence of phi-convex modification of quantile density}.
              That is, $q_\randommeasure=\QuantDensSelec(\randommeasure,\cdot)$ is $\P$-{a.e.\@} $M_{\monotonetransformation{\phi}}$-convex.
              To apply the generalized Jensen inequality \autoref{theo:generalized Jensen inequality}, we need to check the $\P$-{a.e.\@} $M_{\monotonetransformation{\phi}}$-well-definedness of $(\QuantDensSelec(\randommeasure,x),\QuantDensSelec(\randommeasure,y)):\Omega\to\Rbar_{\geq0}^2$ and that $(q_{\mu^*}(x),
                  q_{\mu^*}(y))\in C_{\monotonetransformation{\phi}}$ for all $x,y\in(0,1)$.
              \begin{itemize}
                  \item[\enquote{$\monotonetransformation{\phi}(0)\neq-\infty$}:]
                        If $\monotonetransformation{\phi}(0)\neq-\infty$, then, as $C_{\monotonetransformation{\phi}}=\Rbar^2_{\geq0}$, this holds directly.
                  \item[\enquote{$\monotonetransformation{\phi}(0)=-\infty$}:] If $\monotonetransformation{\phi}(0)=-\infty$, then we know that $\P$-almost surely, $\QuantDensSelec(\randommeasure,x)=0$ for one $x\in(0,1)$ if and only if $\QuantDensSelec(\randommeasure,\cdot)\equiv 0$. By \autoref{lem: aczel correspondence}, for $\monotonetransformation{\phi}(0)=-\infty$, the only $M_{\monotonetransformation{\phi}}$-convex function attaining $0$ is the constant $0$, as otherwise the composition with $\monotonetransformation{\phi}$ would be convex and attain the value $-\infty$. Therefore, $(\QuantDensSelec(\randommeasure,x),\QuantDensSelec(\randommeasure,y)):\Omega\to\Rbar_{\geq0}^2$ is $\P$-{a.e.\@} $M_{\monotonetransformation{\phi}}$-well-defined.
                        If $q_{\mu^*}(x)=0$ for some $x$, then $\QuantDensSelec(\randommeasure,x)=0$ $\P$-{a.e.\@}, so again $\QuantDensSelec(\randommeasure,\cdot)\equiv0$ $\P$-{a.e.\@} and $q_{\mu^*}\equiv0$. Therefore, either $q_{\mu^*}\equiv0$ or $q_{\mu^*}>0$, and thus $(q_{\mu^*}(x),
                            q_{\mu^*}(y))\in C_{\monotonetransformation{\phi}}$ for all $x,y\in(0,1)$.
              \end{itemize}
              Thus, we can apply \autoref{theo:generalized Jensen inequality} to $\QuantDensSelec(\randommeasure,x)\times\QuantDensSelec(\randommeasure,y):\Omega\to C_{\monotonetransformation{\phi}}$ and $M_{\monotonetransformation{\phi},\convexindex}$ for all $x,y,\convexindex\in(0,1)$ independently of the value of $\monotonetransformation{\phi}(0)$.
              Take arbitrary $x,y,\convexindex\in(0,1)$. Then
              \begin{align*}
                  q_{\mu^*}(\convexindex x+ (1-\convexindex)y)
                   & = \mean[\QuantDensSelec(\randommeasure,\convexindex x+ (1-\convexindex)y)]
                  \leq \mean*[M_{\monotonetransformation{\phi},\convexindex}\parenthesis*{\begin{pmatrix}
                                                                                                      q_{\randommeasure}(x) \\
                                                                                                      q_{\randommeasure}(y)
                                                                                                  \end{pmatrix}}]                                                                           \\
                   & \leq M_{\monotonetransformation{\phi},\convexindex}\parenthesis*{\mean*[\begin{pmatrix}
                                                                                                         q_{\randommeasure}(x) \\
                                                                                                         q_{\randommeasure}(y)
                                                                                                     \end{pmatrix}]} = M_{\monotonetransformation{\phi},\convexindex}(q_{\mu^*}(x),q_{\mu^*}(y)).
              \end{align*}
        \item[\enquote{$\Rightarrow$}:]
              Let $\monotonetransformation{\phi}:\R_{>0}\to\R$ be a generator of $M_{\monotonetransformation{\phi}}$, and suppose that $M_{\monotonetransformation{\phi}}$ is not concave. By \autoref{lem Convexity of mean}, we can even assume that $M_{\monotonetransformation{\phi}}$ is not Jensen concave.
              Thus, there exist $\theta=(\theta_1,\theta_2),\vartheta=(\vartheta_1,\vartheta_2)\in\R_{>0}^2$ such that
              \begin{equation*}
                  M_{\monotonetransformation{\phi},\frac12}\parenthesis*{\frac{\theta+\vartheta}{2}}<\frac{M_{\monotonetransformation{\phi},\frac12}(\theta)+M_{\monotonetransformation{\phi},\frac12}(\vartheta)}{2}.
              \end{equation*}
              Now, for some $a_1<a_2$ in $(0,1)$, sufficiently close to $0$ and $1$, respectively, there exist two affine functions $X,Y:(0,1)\to\Functionimage{\monotonetransformation{\phi}}{\R_{>0}}$ such that $X(a_j)=\monotonetransformation{\phi}(\theta_j)$ and $Y(a_j)=\monotonetransformation{\phi}(\vartheta_j)$ for $j\in\set{1,2}$.
              The mappings $\monotonetransformation{\phi}^{-1}\circ X,\monotonetransformation{\phi}^{-1}\circ Y$ are positive real $M_{\monotonetransformation{\phi}}$-affine functions on $(0,1)$. We use them to define two distributions with $M_{\monotonetransformation{\phi}}$-convex quantile densities by setting $Q_{\mu_j}(x)\defeq\int_{\frac12}^xq_{\mu_j}\diff\Lebesgue$ for $q_{\mu_1}\defeq\monotonetransformation{\phi}^{-1}\circ X$ and $q_{\mu_2}\defeq\monotonetransformation{\phi}^{-1}\circ Y$ and $j\in \set{1,2}$.
              By construction, we have $q_{\mu_1}(a_j)=\theta_j$ and $q_{\mu_2}(a_j)=\vartheta_j$ for $j\in\set{1,2}$.
              The distribution $\popmeasure=\frac{1}{2}(\delta_{\mu_1}+\delta_{\mu_2})$ is supported on $\phi$-unimodal distributions, and since it is a discrete distribution of distributions with bounded support, we further have $\popmeasure\in\Wtwo(\Wtwo(\R))$. Denote again $\mu^*=\Barycenterfunc{2}{\Wtwo}\,(\popmeasure)$. By the continuity of $q_{\mu_j}$, we have $q_{\mu_j}=\QuantDensSelec(\mu_j,\cdot)$. Thus, $q_{\mu^*}\defeq\frac{1}{2}(q_{\mu_1}+q_{\mu_2})$ is a valid quantile density of $\mu^*$ and is continuous. Therefore, $q_{\mu^*}=\QuantDensSelec(\mu^*,\cdot)$. However, since
              \begin{align*}
                  \QuantDensSelec*{\mu^*,\frac{a_1+a_2}{2}}
                   & = \frac{1}{2}\parenthesis*{q_{\mu_1}\parenthesis*{\frac{a_1+a_2}{2}}+q_{\mu_2}\parenthesis*{\frac{a_1+a_2}{2}}}                                 \\
                   & =\frac{1}{2}(M_{\monotonetransformation{\phi},\frac12}(q_{\mu_1}(a_1),q_{\mu_1}(a_2))+M_{\monotonetransformation{\phi},\frac12}(q_{\mu_2}(a_1),q_{\mu_2}(a_2)))
                  =\frac{1}{2}(M_{\monotonetransformation{\phi},\frac12}(\theta)+M_{\monotonetransformation{\phi},\frac12}(\vartheta))                                                            \\
                   & >M_{\monotonetransformation{\phi},\frac12}\parenthesis*{\frac{\theta+\vartheta}{2}}
                  =M_{\monotonetransformation{\phi},\frac12}\parenthesis*{\frac{q_{\mu_1}(a_1)+q_{\mu_2}(a_1)}{2},\frac{q_{\mu_1}(a_2)+q_{\mu_2}(a_2)}2}                     \\
                   & =M_{\monotonetransformation{\phi},\frac12}(q_{\mu^*}(a_1),q_{\mu^*}(a_2)),
              \end{align*}
              we have that $\QuantDensSelec{\mu^*,\cdot}$ is not $M_{\monotonetransformation{\phi}}$-convex. As $\mu^*$ is quantile absolutely continuous by \autoref{prop: Wasserstein 2 Barycenter of quantile absolutely continuous distributions on R have absolutely continuous quantile functions}, $\mu^*$ is not quantile $\monotonetransformation{\phi}$-unimodal by \autoref{lem: existence of phi-convex modification of quantile density}, and hence, by \autoref{def: M concave general}, not $\phi$-unimodal.\qedhere
    \end{itemize}
\end{delayedproof}

\section[Counterexamples for extensions to higher dimensions and arbitrary \texorpdfstring{$p$}{p}]{Counterexamples for extensions to higher dimensions and \\arbitrary \texorpdfstring{$p$}{p}}
\label{sec: Counterexamples for extensions of the main theorem}
    
Both counterexamples have a common structure. For all generators $\phi$, the set of \emph{$\infty$-unimodal} distributions, {i.e.}, uniform distributions on convex sets, is a subset of the set of all $\phi$-unimodal distributions. These $\infty$-unimodal distributions are obtained by taking $s\to\infty$ in the definition of $\phi_s$-unimodal distributions. In the univariate setting, they are quantile absolutely continuous distributions with an $M_{-\infty}=\min$-convex, and hence constant, quantile density. Any constant quantile density is also $M_{\monotonetransformation{\phi}}$-convex for all generators $\phi$, and any constant density is $M_{\phi}$-concave on the interior of its support.
We give examples showing that a Wasserstein barycenter of uniform distributions may fail to be unimodal. This then rules out Wasserstein barycenter convexity of the set of $\phi$-unimodal distributions for all generators $\phi$ in each setting.
\subsection{The set of \texorpdfstring{$\phi$}{ɸ}-unimodal distributions is not \texorpdfstring{{$p$-WB}\nocorr}{p-WB} convex for \texorpdfstring{$p\neq2$}{p≠2}}
\begin{figure}[ht]
        \centering
        \includegraphics{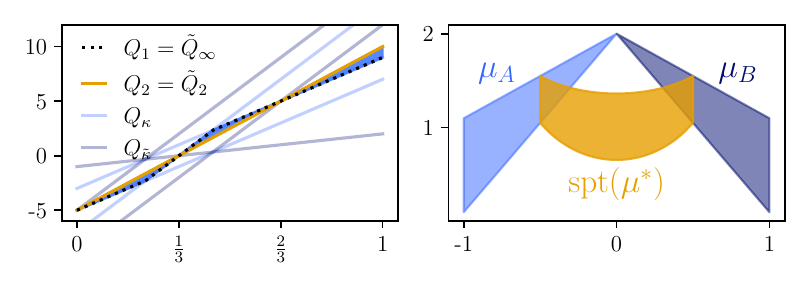}
        \caption{{L.h.s.\@}: The {q.f.s} of the two random distributions $\randommeasure\sim\popmeasure$ (light blue, short-dashed) and $\tilde{\randommeasure}\sim\tilde{\popmeasure}$ (yellow, long-dashed) from \autoref{example: Counterexample Wasserstein p barycenter convexity for p neq 2} both generate the same \WB{2} $\mu_2^*=\tilde{\mu}_2^*$ (orange, solid). Furthermore, the enveloping {q.f.s} $Q_1$ and $\tilde{Q}_{\infty}$ (black, dotted) coincide. By the monotonicity of the barycenter $\Lambda_p^\R$ of three points, we have that the {q.f.s} of the indicated \WB{p} of $\popmeasure$ and $\tilde{\popmeasure}$ must lie in the strict envelope (blue) of $Q_1$ and $Q_{2}$ except at the crossing points $\frac13$ and $\frac23$. However, no {q.f.\@} of a unimodal or symmetric distribution fits into that region. Moreover, $Q_1$ does not correspond to a unimodal or a symmetric distribution. {R.h.s.\@} \cite[Compare Fig. 6.]{santambrogioConvexitySupportDisplacement2016}: The uniform distributions $\mu_A,\mu_B$ with visualized supports $A$ (light blue) and $B$ (dark blue), respectively, are unimodal, but the support of their equal-weight barycenter $\mu^*$ is non-convex. Therefore, its density $f_{\mu^*}$ is not quasi-concave.}
        \label{fig: counterexample unimodal preservation for p neq 2 and d>1}
    \end{figure}

    \begin{example}
        [Counterexample to the extension of \autoref{theo: Wasserstein 2 barycenter convexity of M convex funcitons} for $p\neq2$]
        \label{example: Counterexample Wasserstein p barycenter convexity for p neq 2}
        We consider the equal-weight mixture distributions $\popmeasure=\frac13\sum_{j=1}^3\delta_{\mu_j}$ and $\tilde{\popmeasure}=\frac13\sum_{j=1}^3\delta_{\tilde{\mu}_j}$, each with three continuous uniform mixture components, as given in \autoref{tab: uniform counterexmaples p neq2}. Define the {q.f.s}
        \begin{equation*}
            \begin{aligned}
                Q_1(t)\defeq{}&{}\median(Q_{\mu_1}(t),Q_{\mu_2}(t), Q_{\mu_3}(t))=\midrange(Q_{\tilde{\mu}_1}(t),Q_{\tilde{\mu}_2}(t),Q_{\tilde{\mu}_3}(t)),\\
            Q_2(t)\defeq{}&{}\frac13(Q_{\mu_1}(t)+Q_{\mu_2}(t)+Q_{\mu_3}(t))=\frac13(Q_{\tilde{\mu}_1}(t)+Q_{\tilde{\mu}_2}(t)+Q_{\tilde{\mu}_3}(t)),
            \end{aligned}\quad t\in(0,1).
        \end{equation*}
        They are the {q.f.s} of the $1$-WB and {$2$-WB}, respectively, of $\popmeasure$ by \cite[Proposition 4.1.]{carlierWassersteinMediansRobustness2024} and \cite[Theorem 5.1]{brizziWassersteinBarycenters2025}, respectively.
        We detail in \autoref{sec: Details counterexample} that for all $1< p<+\infty$, the {q.f.s} of the {$p$-WBs} of $\popmeasure$ (for $1< p<2$) and $\tilde{\popmeasure}$ (for $2<p<+\infty$) lie strictly inside the region bounded by $Q_1$ and $Q_{2}$, except at the crossing points $\frac13,\frac23$. As can be seen from the {l.h.s.\@} of \autoref{fig: counterexample unimodal preservation for p neq 2 and d>1}, no function lying strictly inside this region can be first concave and then convex. In particular, the {$p$-WBs} cannot be unimodal \cite[Theorem 1.6.]{zobelmarkusQuantileCharacterizationUnivariate2026arxiv}.
        For $p=1$, we have that $\Barycenterfunc{1}{\W_1}(\popmeasure)$ is not unimodal, as $Q_1$ fails to be first concave and then convex.
    \end{example}
    \begin{table}[!h]
            \centering
            \caption{Mixture components in \autoref{example: Counterexample Wasserstein p barycenter convexity for p neq 2}.}
            \label{tab: uniform counterexmaples p neq2}
            \begin{tabular}{llll}
                \toprule                  & $j=1$   & $j=2$   & $j=3$   \\
                \midrule$\support{\mu_j}$ & [-7,14] & [-5,7]  & [-3,9]  \\
                $\support{\tilde{\mu_j}}$ & [-1,2]  & [-9,12] & [-5,16] \\
                \bottomrule
            \end{tabular}
        \end{table}
        Therefore, the set of $\phi$-unimodal distributions is not {$p$-WB} convex for $p\neq2$ for any generator $\phi$.
    \begin{remark}
        This counterexample rules out the claim in \cite{backhoffStochasticGradientDescent2025} that univariate symmetric distributions are {$p$-WB} convex, as all continuous uniform distributions are symmetric. 
        We call a distribution $\mu$ \emph{symmetric} if there exists $m\in\R$ such that $X-m\Deq-(X-m)$ for $X\sim\mu$. This is equivalent to $F_{\mu}(m+x)=1-F_{\mu}(m-x)$ at all continuity points $m+x\in\R$ of $F_\mu$ and to $Q_\mu(t)+Q_\mu(1-t)=2m$ at all continuity points $t$ of $Q_\mu$, {i.e.}, the graphs of $F_\mu$ and $Q_\mu$ are point symmetric around $(m,\frac12)$ and $(\frac12,m)$, respectively.
    \end{remark}
   
\subsection{The set of \texorpdfstring{$\phi$}{ɸ}-unimodal distributions is not \texorpdfstring{{$p$-WB}\nocorr}{p-WB} convex on \texorpdfstring{$\R^d$}{ℝᵈ} for any \texorpdfstring{$d\geq2$}{d≥2} and \texorpdfstring{$p\in[1,\infty)$}{p∊[1,∞)}}
    \begin{example}[Counterexample to \autoref{theo: Wasserstein 2 barycenter convexity of M convex funcitons} for $d=2$]
        In \cite[Theorem 2.]{santambrogioConvexitySupportDisplacement2016}, the authors construct two convex sets $A,B$ in $\R^2$.
        If one takes the uniform measures $\mu_A=\frac{\1_A}{\Lebesgue(A)}\cdot\Lebesgue$, $\mu_B=\frac{\1_B}{\Lebesgue(B)}\cdot\Lebesgue$, then the support of their equal-weight {$2$-WB} $\mu^*$, depicted on the {r.h.s.\@} of \autoref{fig: counterexample unimodal preservation for p neq 2 and d>1}, is not convex. Thus, $\mu^*$ is not unimodal.
        Therefore, the set of $\infty$-unimodal distributions is not $2$-displacement convex and thus not {$2$-WB} convex.

    \end{example}

    For general $d>2$, we reuse these sets and take the McCann interpolant between $\mu_0=\mu_A\tensor\uniformmeasure^{d-2}$ and $\mu_1=\mu_B\tensor\uniformmeasure^{d-2}$. Both distributions are absolutely continuous and again ${\infty}$-unimodal.
    In general, $W_{\!p}(\mu\tensor\nu,\rho\tensor\nu)=W_{\!p}(\mu,\rho)$, $\mu,\rho\in\Wp(\R^d)$, $\nu\in\Wp(\R^k)$, with $\mu\tensor\nu,\rho\tensor\nu\in\Wp(\R^{d+k})$. The inequality \enquote{$\leq$} follows from a deterministic coupling in the last $k$ dimensions together with an optimal coupling of $\mu$ and $\rho$, while \enquote{$\geq$} holds because projection onto the first $d$ dimensions does not increase the $p$-Wasserstein distance. From this, it follows that $\mu^*\tensor\uniformmeasure^{d-2}$ is a barycenter of $\frac{1}{2}(\delta_{\mu_0}+\delta_{\mu_1})$, again by the projection inequality.
    The support of $\mu^*\tensor\uniformmeasure^{d-2}$ is still not convex, and the {$2$-WB} is thus not unimodal.

    For $p\neq2$ and $d\geq2$, we can do the same with the counterexample in \autoref{example: Counterexample Wasserstein p barycenter convexity for p neq 2}, again extended to dimension $d\geq2$ via independent couplings with $\uniformmeasure^{d-1}$.
    Thus, for $p\in[1,+\infty)\setminus\set{2}$, we have that $\mu_p=\mu^*\tensor \uniformmeasure^{d-1}$ is a {$p$-WB} of the extended mixtures $\sum_{j=1}^3\frac13\delta_{\mu_j\tensor\uniformmeasure^{d-1}}$ for $p\in[1,2)$, $\mu^*\defeq\mu_p^*$, and $\sum_{j=1}^3\frac13\delta_{\tilde{\mu}_j\tensor\uniformmeasure^{d-1}}$ for $p\in(2,+\infty)$, $\mu^*\defeq\tilde{\mu}_p^*$, respectively, based on \autoref{tab: uniform counterexmaples p neq2}.
    If $\mu_p$ were to admit a quasi-concave density $f_{\mu_p}$, then we would still have $f_{\mu_p}(x,y)=f_{\mu^*}(x)\1_{[0,1]^{d-1}}(y)$ for $\Lebesgue$-almost all $(x,y)\in\R^{1+d-1}$. Therefore, there exists a $y\in [0,1]^{d-1}$ such that $f_{\mu_p}(\cdot,y)=f_{\mu^*}$ $\Lebesgue$-{a.e}.
    This is a contradiction to \autoref{example: Counterexample Wasserstein p barycenter convexity for p neq 2} as $f_{\mu_p}(\cdot,y)$ would be a quasi-concave density of $\mu^*$ inheriting the quasi-concavity from $f_{\mu_p}$.

\section{Discussion}
\label{sec: Discussion}
 
\subsection{Different notions of barycenter convexity}
    \label{sec: Relation of different notions of barycenter convexity results}
    We give a short overview of the relations between the three notions of barycenter convexity mentioned in the introduction (see \autoref{tab: different concepts of barycenter convexity}), visualized in \autoref{tikz: Diagram realtions convexity}.
    It is shown in \cite[Proposition 7.6]{aguehBarycentersWassersteinSpace2011} that $2$-displacement convexity and convexity along $2$-barycenters coincide on $\Wtwo(\R)$.
    This does not hold for $p\neq 2$. \autoref{example: Counterexample Wasserstein p barycenter convexity for p neq 2} shows that the set of uniform distributions, which is displacement convex, is not convex along $p$-barycenters.
    Under the assumption that the set $A\subset\Wp(\R)$ is $W_{\!p}$-closed, convexity along $p$-barycenters and {$p$-WB} convexity are equivalent for $p\in(1,+\infty)$.
    We can approximate any $\popmeasure\in\Wp(\Wp(\R))$ supported on $A$ by $\popmeasure_n\in\mathfrak{M}_{<+\infty}$, {i.e.}, $\setcount{\support(\popmeasure_n)}<+\infty$,
    also supported on $A$, with $\popmeasure_n\emptyconvergencearrow{W_{\!p}}\popmeasure$ as $n\to\infty$.
    Then, we obtain from \cite[Theorem 3]{LeGouic.2017} that $\Barycenterfunc{p}{\Wp}({\popmeasure_n})\emptyconvergencearrow{W_{\!p}}\Barycenterfunc{p}{\Wp}(\popmeasure)$, as the barycenter of $\popmeasure$ is unique.
    Thus, if a set $A\in\B(\W_p(\R))$ is convex along $p$-barycenters and $W_{\!p}$-closed, then it is also {$p$-WB} convex. Extending the above argument to $\Wp(\R^d)$ would require a refinement of \cite[Theorem 3]{LeGouic.2017}, as one loses the uniqueness of {$p$-WBs} in general.
    Without closedness, convexity along $p$-barycenters need not imply {$p$-WB} convexity even on $\R$.
    \begin{figure}[ht]
        \centering

\begin{tikzcd}[
        every cell/.append style = {font = \normalsize},
        cells={nodes={execute at begin node=\displaystyle}}
    ]
    \text{$A$ {$p$-WB} convex}
    \arrow[r, Rightarrow, yshift = 3 pt]
    &
    \text{$A$ convex along $p$-barycenters}
    \arrow[l, "\text{ $\R$, $A$ $W_{\!p}$ closed, $p>1$}"{yshift=-5pt}, Rightarrow, yshift = -3 pt]
    \arrow[r, Rightarrow, yshift = 3 pt]
    &
    \text{$A$ $p$-displacement convex}
    \arrow[l, "\text{ $\R,p=2$}"{yshift=-5pt}, Rightarrow, yshift = -3 pt]
\end{tikzcd}

        \caption{Relations between the three notions of barycenter convexity from \autoref{tab: different concepts of barycenter convexity}.}
        \label{tikz: Diagram realtions convexity}
    \end{figure}
    \begin{example}[A set convex along $2$-barycenters that is not {$2$-WB} convex]
        \label{example: difference dispalcement convexity and barycenter convexity}
        Let $A$ denote the set of all finitely supported distributions on $\R$, {i.e.}, $A=\set*{\sum_{j=1}^n\convexindex_j\delta_{x_j}\given  \convexindex\in\weightings{n},(x_1,\dots,x_n)\in\R^n,n\in\N}$.\footnote{Here $\weightings{n}\defeq\set{(\convexindex_1,\dots,\convexindex_n)\in[0,1]^n\given\sum_{j=1}^n\convexindex_j=1}$.}
        Its {q.f.s} are given by
        \begin{equation*}
            \Functionimage{\Quantiletransform}{A}=\set[\Big]{x+\sum_{j=1}^n h_{j}\ind_{[W_j,1)}\given x\in\R,(h_1,\dots,h_n)\in\R_{\geq0}^{n},(W_1,\dots,W_n)\in[0,1]^n,n\in\N}.
        \end{equation*}
        Then $A\subset\Wtwo(\R)$ is $2$-displacement convex, as for $\mu_0,\mu_1\in A$, we have
        \begin{equation*}
            \convexindex Q_{\mu_1} + (1-\convexindex) Q_{\mu_0}
            = \convexindex x_1+(1-\convexindex)x_0+\sum_{j=1}^{n_0+n_1}\tilde{h}_{j}\ind_{[\tilde{W}_j,1)}\in \Quantiletransform{A},\quad\convexindex\in(0,1).
        \end{equation*}
        Therefore, it is also convex along $2$-barycenters.
        For $\popmeasure=\sum_{k=1}^\infty \frac1{2^k}\delta_{\mu_{Q_k}}$, $Q_k=\ind_{\left[\frac{1}{k+1},1\right)}\in \Quantiletransform(A)$, $k\in\N$, we have that $\popmeasure$ is supported on $A$, {i.e.}, $\support{\popmeasure}
            = \set{\mu_{Q_k}\given k\in\N}\cup\set{\delta_1} \subset A$, and $\popmeasure\in\Wtwo(\Wtwo(\R))$, as \begin{equation*}
                W_2(\popmeasure,\delta_{\delta_1})^2=\sum_{k=1}^{\infty}\frac{1}{2^k(k+1)}=2\log(2)-1<1.
            \end{equation*}
            However, $Q_{\Barycenterfunc{2}{\Wtwo}\,(\popmeasure)}=\sum_{k=1}^\infty\frac{1}{2^k}\ind_{\left[\frac{1}{k+1},1\right)} \notin \Quantiletransform(A)$.
        Even though $A$ is convex along $2$-barycenters, it is not {$2$-WB} convex.
    \end{example}
    The set of univariate log-concave distributions, as defined in \autoref{def: M concave general}, is weakly closed and thus, in particular, $W_2$-closed. Therefore, the $2$-displacement convexity of the set of univariate log-concave distributions already suffices for its {$2$-WB} convexity.
    The subset of absolutely continuous log-concave distributions is not weakly closed. As mentioned in the introduction, it can be written as the intersection of two {$2$-WB} convex subsets, namely, the set of all log-concave distributions and the set of all absolutely continuous distributions \cite[Theorem]{maAbsoluteContinuityWasserstein2025}, and is therefore itself {$2$-WB} convex.
\subsection{Generalization of absolutely continuous \texorpdfstring{$\phi$}{ɸ}-unimodal distributions}
    \label{sec: discussion generalization of phi concavity}
    Generalizations of absolutely continuous $\phi$-unimodal distributions, such as the one in \autoref{def: M concave general}, are common in the unimodality literature \cite{Dharmadhikari.1988}, as the set of absolutely continuous unimodal or log-concave distributions is not weakly closed. For the latter, the only non-absolutely continuous elements in the closure are Dirac distributions. For unimodal distributions, mixtures of Dirac distributions and absolutely continuous unimodal distributions are possible, as long as both share the same mode.
    The classical definition of general unimodal distributions uses the {c.d.f.\@}, while log-concave distributions are defined via an inequality involving the measure itself. A distribution $\mu\in\distributionon(\R)$ is \emph{log-concave} if, for all convex sets $A,B\subset\R$, the following holds \cite{prekopaLogarithmicConcaveMeasures1980}:
    \begin{equation}
        \label{eqn: Logconcave measure inequality}
        \mu(\convexindex A + (1-\convexindex) B)\geq M_{0,\convexindex}(\mu(A),\mu(B)),\quad \convexindex\in(0,1),
    \end{equation}
    where $\convexindex A + (1-\convexindex) B=\set{\convexindex x + (1-\convexindex)y \given x\in A , y\in B}$. It was then shown that if such a log-concave distribution is absolutely continuous, then it has a log-concave density \cite{prekopaLogarithmicConcaveMeasures1980}.

    This inequality and a similar characterization have been adapted to define $s$-concave distributions. A distribution $\mu$ is \emph{$s$-concave} if, for all convex $A,B$ such that $\mu(A),\mu(B)>0$, \eqref{eqn: Logconcave measure inequality} holds with $M_0$ replaced by $M_s$. For $s>0$, the condition on the masses of $A$ and $B$ is necessary to allow $s$-concave distributions $\mu$ with $\support{\mu}\neq\R$. For $s\leq0$, the inequality extends to all convex sets $A,B$.
    For $s\neq0$, an absolutely continuous $s$-concave distribution is not characterized by having an $s$-concave density. An absolutely continuous distribution is $t$-concave if and only if it has an $s$-concave density with $s=\frac{t}{1-t}$ \cite[Theorem 3.2]{borellConvexSetFunctions1975} or, alternatively, has an $M_r$-convex quantile density with $r= -\frac{1}{1-t}$, as we have shown. Even though $s$-concave densities cover all of $s\in\Rbar$, the only non-degenerate $t$-concave distributions lie in the range $-\infty\leq t\leq1$, with the boundary values $t\in\set{-\infty,1}$ covered by extending the relations continuously. This encompasses the $s$-concave densities for $-1\leq s<+\infty$ and the case of constant densities on intervals, {i.e.}, absolutely continuous $\infty$-unimodal distributions.

    Such a theory is missing for extending $\phi$-unimodal distributions beyond the absolutely continuous case. The authors conjecture that an investigation of $\phi$-unimodal measures defined in a manner analogous to \eqref{eqn: Logconcave measure inequality} will yield exactly the extension proposed in \autoref{def: M concave general}, {i.e.}, that the extension coincides with the definition of quantile $\monotonetransformation{\phi}$-unimodal distributions.

\section*{Acknowledgments}
We are grateful to Jan Victor Otte for various valuable comments. Generative AI was used for proofreading both the mathematical content and the grammar.

 \appendix
\section*{Appendix}
\label{sec supp}

The Appendix is twofold. Firstly, we reintroduce the setting of non-decreasing functions, their generalized inverses and the connections with their Stieltjes measures.
Secondly, we turn to proofs not given in the main text in \autoref{sec: Proofs of main paper}.
\section[Measurability and continuity on the set of non-decreasing functions]{Measurability and continuity on the set of non-decreasing\\functions}
\label{sec app: Non-decreasing functions}

In the following we give an introduction to the theory of non-decreasing functions on $\R$.
The focus lies on their connection to distributions on $\R$ via the {c.d.f.\@} and to locally finite measures on $(0,1)$ via the {q.f}. We establish measurability of the transformations between distributions, quantiles and inverse measures.
\subsection{Non-decreasing functions}

A natural setting for discussing real-valued non-decreasing functions is to start with extended-real-valued \emph{non-decreasing} functions $G:\Rbar\to\Rbar$, satisfying $G(x)\leq G(y)$ if $x\leq y$ for $x,y\in\Rbar$, where $\Rbar$ is the order completion of $\R$, {i.e.},
\begin{equation*}
-\infty = \inf{\R}\leq x \leq \sup{\R}=+\infty,\quad x\in\R.
\end{equation*}
Let $\increasingfunctions\defeq\set{G:\Rbar\to\Rbar\given \text{$G$ is non-decreasing}}$ denote the set of all such non-decreasing functions.
Then $\increasingfunctions\subset\Lp(\Rbar)$, the set of all measurable functions $G:\Rbar\to\Rbar$, as preimages of intervals under non-decreasing functions are intervals. Further, non-decreasing functions have left and right limits at each point and have only countably many discontinuities. For $G\in\increasingfunctions$, define the \emph{left- and right-continuous versions} $G_l$ and $G_r$, respectively, by
\begin{equation*}
    \begin{aligned}
        G_l(x) & \defeq G(x-)=\sup_{x>y\in \Rbar}{G(y)},
        \\
        G_r(x) & \defeq G(x+)=\inf_{x<y\in \Rbar}{G(y)},
    \end{aligned}
    \quad\forall x\in\Rbar
    .
\end{equation*}
By the convention $\sup\emptyset=-\infty$, $\inf\emptyset=+\infty$, we always have $G_l(-\infty)=-\infty$ and $G_r(+\infty)=+\infty$.
It is immediately clear that $G_l\leq G\leq G_r$ for each $G\in\increasingfunctions$ and $G_l=G= G_r$ $\Lebesgue$-{a.e.\@} for the Lebesgue measure extended to $(\Rbar,\B(\Rbar))$ with $\Lebesgue(\set{-\infty,+\infty})=0$. We always choose the largest right-continuous version as a representative and denote this collection by $\increasingfunctionclasses\defeq\set{G_r\given G\in\increasingfunctions}$.
Therefore, $\increasingfunctionclasses\subset\Lpmod(\Rbar)$, the set of all $\Lebesgue$-{a.e.\@} equivalence classes of measurable functions on $\Rbar$.
The real-valued right-continuous non-decreasing functions can be identified with elements of $\increasingfunctions$ .
Denote by $\increasingfunctions(I)\subset\Lp(I)$ the set of all non-decreasing real-valued functions $G:I\to\R$ for an open interval $I\subset\R$. We embed the real-valued non-decreasing functions $\increasingfunctions(I)$ into $\increasingfunctions$ by extending each element with $-\infty$ to the left of $I$ and with $+\infty$ to the right of $I$. Analogously, we denote by $\increasingfunctionclasses(I)\subset\Lpmod(I)$ the set of all right-continuous non-decreasing functions on $I$, which we embed into $\increasingfunctionclasses$ by first embedding them into $\increasingfunctions$ as before and then taking the maximal right-continuous version. In the following discussion, we will also consider extensions to the closure of $I$ instead of extensions to all of $\Rbar$, using the one-sided limits towards the boundary. While these functions will still be right-continuous on the closure of $I$, they might not be the restrictions to the closure of $I$ of the maximal right-continuous embeddings introduced above.

\subsection{Generalized inverses and Stieltjes measures}

For a non-decreasing function $G\in\increasingfunctionclasses$, we define its \emph{generalized inverse} as
\begin{equation*}
    \generalizedinverse(G)(t)\defeq \sup(\preimage{G}{[-\infty,t]}),\quad t\in\Rbar.
\end{equation*}
These definitions assemble into a well-defined function $\generalizedinverse:\increasingfunctionclasses\to\increasingfunctionclasses$, as the resulting function is again non-decreasing and right-continuous \cite{zobelmarkusQuantileCharacterizationUnivariate2026arxiv}. Furthermore, it is an involution.
To be able to use this symmetry, we use right-continuous {q.f.s} instead of the commonly used left-continuous ones.
Still, the generalized inverse need not be an actual inverse, {i.e.}, $\generalizedinverse(G)\circ G\neq \id$ in general, but if $G$ is bijective, and hence continuous and strictly increasing, then $\generalizedinverse{G}=G^{-1}$.
We denote by $\increasingfunctionclasses^{*}(I)\defeq\Functionimage{\generalizedinverse}{\increasingfunctionclasses(I)}$ the set of all generalized inverses of elements in $\increasingfunctionclasses(I)$. If $I=(a,b),a,b\in\R$, then one can easily check that
\begin{equation*}
    \increasingfunctionclasses^{*}(I)=\set{G\in\increasingfunctionclasses(\R)\given\lim_{x\to-\infty}G(x)=a,\lim_{x\to+\infty}G(x)=b},
\end{equation*}
where we identify the original $\increasingfunctionclasses^{*}(I)$ with its restrictions to $\R$, which are all real-valued for bounded $I$.
Thus, $\increasingfunctionclasses^{*}((0,1))$ denotes exactly the set of all {c.d.fs}.
In general, we can define $\Stieltjestransform:\increasingfunctionclasses(I)\to\locallyfinitemeasureon(I)$ by $\Stieltjestransform(G)((c,d])\defeq G(d)-G(c)$ for $c<d\in I$ \cite{natansonTheorieFunktionenReellen1961}, as $\set{(c,d]\given c,d\in I}$ generates $\B(I)$. For each $c$ in $I$, we can give a right-inverse $\CDFtransform^c:\locallyfinitemeasureon(I)\to\increasingfunctionclasses(I)$ as
\begin{equation*}
    \CDFtransform^c(\mu)(x)\defeq \begin{cases}
        -\mu((x,c]) & x\leq c \\
        \mu((c,x])  & c<x,
    \end{cases}\quad x\in I.
\end{equation*}
Therefore, $\Stieltjestransform$ is surjective but not injective, as $\Stieltjestransform(G)=\Stieltjestransform(G+x)$ for all $x\in\R$.
Assume now that the open interval $I=(a,b),a,b\in\R$ is bounded. Then the restriction $\Stieltjestransform:\increasingfunctionclasses^{*}(I)\to\locallyfinitemeasureon(\R)$ is injective, as the limiting conditions rule out constant shifts. Further, $\Functionimage{\Stieltjestransform}{\increasingfunctionclasses^{*}(I)}$ contains only finite measures with mass equal to $b-a$, the length of $I$. The composition $\Transportmeasure\defeq\Stieltjestransform\circ\generalizedinverse:\increasingfunctionclasses(I)\to\locallyfinitemeasureon(\R)$ also has a closed form as the pushforward measure of the Lebesgue measure restricted to bounded $I$. For this, we refer to \cite[Lemma 2.5]{zobelmarkusQuantileCharacterizationUnivariate2026arxiv}, where we prove that
\begin{align}
    \Transportmeasure(G)((c,d])
    &=\pushforward{G}{\restr{\Lebesgue}{I}}((c,d])=\generalizedinverse(G)(d)-\generalizedinverse(G)(c),&& G\in\increasingfunctionclasses(I),c<d\in\R.\nonumber
    \intertext{Restricting the more general notation there, we equivalently define $\Transportmeasure=\Stieltjestransform\circ\generalizedinverse:\increasingfunctionclasses^{*}{I}\to\locallyfinitemeasureon(I)$ satisfying}
    \label{eqn: transport formulation of inverse of Quantile transform}
    \Transportmeasure(G)((c,d])&=\pushforward{G}{\Lebesgue}((c,d]),&& G\in\increasingfunctionclasses^{*}(I),a<c<d<b.
\end{align}

\subsection{Convergence of non-decreasing functions is metric}

For the following discussion, fix an open interval $I\subset\R$.

Let $(G_n)_{n\in\N}\subset\increasingfunctionclasses{I}$, $G\in\increasingfunctionclasses{I}$, $\nu\in\locallyfinitemeasureon(I)$. Then we say that \emph{$G_n$ converges at continuity points to $G$}, $G_n\emptyconvergencearrow{\mathcal{C}}G$, if $G_n(x)\nto G(x)$ for all $x$ such that $G$ is continuous at $x$.
If we equip $I$ with a locally finite measure $\nu\in\locallyfinitemeasureon(I)$ which is further absolutely continuous $\nu\ll\Lebesgue$ then we get additional notions of convergence:

\begin{itemize}
    \item \emph{$G_n$ converges in $\nu$ to $G$}, $G_n\emptyconvergencearrow{\nu}G$, if $\nu(\abs{G_n-G}>\epsilon)\to0$ for each $\epsilon>0$,
    \item \emph{$G_n$ converges $\nu$-{a.e.\@} to $G$}, $G_n\emptyconvergencearrow{\text{$\nu$-{a.e.\@}\nocorr}}G$, if $\nu(\lim_{n\to\infty}\abs{G_n-G}\neq0)=0$.
\end{itemize}

The last notions of convergence we introduce arise from equipping $\increasingfunctionclasses(I)$ with one of the following $\Lpmod^p$-metrics. For each $\nu\in\locallyfinitemeasureon(I)$ with full support, we can define the \emph{$\Lpmod^p$-metric} $d_{\nu,p}$ for $p\in[1,\infty)$ as
\begin{align*}
    d_{\nu,p}(G,H)
    \defeq & \parenthesis*{\int_I d(G(t),H(t))^p\diff\nu(t)}^{\frac{1}{p}}=\parenthesis*{\int_I\abs{G(t)-H(t)}^p\diff\nu(t)}^{\frac{1}{p}} \\
    =      & \norm{d(G,H)}_{\Lpmod^p(\nu)}
\end{align*}
on the metric space $\Lpmod^p(\nu,d)=\set{G\in\Lpmod(I)\given d_{\nu,p}({G,0})<\infty}$ of $p$-integrable functions on $I$ {w.r.t.\@} the Euclidean distance $d$.
If we replace the Euclidean distance with $\bar{d}(x,y)=\min(d(x,y),1)$, then $\bar{d}_{\nu,p}$, defined analogously, gives us $\Lpmod^p(\nu,\bar{d})$.

Restricting to the non-decreasing functions, we obtain the metric spaces $\increasingfunctionclasses^{p}(\nu,\tilde{d})\defeq\increasingfunctionclasses(I)\cap\Lpmod^p(\nu,\tilde{d})$ for $\tilde{d}\in\set{d,\bar{d}}$.

Let $(G_n)_{n\in\N}\subset\increasingfunctionclasses^{p}{\nu,\tilde{d}}$ and $G\in\increasingfunctionclasses^{p}{\nu,\tilde{d}}$. Then $G_n\emptyconvergencearrow{\tilde{d}_{\nu,p}}G$ denotes the metric convergence of $G_n$ to $G$ in $\increasingfunctionclasses^{p}(\nu,\tilde{d})$.

On $\increasingfunctionclasses^{p}(\nu,d)$, metric convergence is the strongest of these notions. However, all the other notions of convergence coincide and are induced by the finite metric $\bar{d}_{\nu,1}$ when $\nu$ is a finite absolutely continuous measure with full support on $I$. For these choices, $\increasingfunctionclasses(I)=\increasingfunctionclasses^{1}(\nu,\bar{d})$.
\begin{lemma}
    \label{lem: convergence on non-decreasing functions}
    For a nonvoid open interval $I\subset\R$, $\nu\in\locallyfinitemeasureon(I)$ with full support, $p\in[1,\infty)$, and $\tilde{d}\in\set{d,\bar{d}}$, the following implications hold on $\increasingfunctionclasses^{p}(\nu,\tilde{d})$:
    \begin{equation*}
        G_n\emptyconvergencearrow{\tilde{d}_{\nu,p}}G\Rightarrow
        G_n\emptyconvergencearrow{\nu}G\Rightarrow G_n\emptyconvergencearrow{\mathcal{C}}G\Rightarrow
        G_n\emptyconvergencearrow{\text{$\Lebesgue$-{a.e.\@}\nocorr}}G.
    \end{equation*}
    If both $\tilde{d}=\bar{d}$ and $\nu\ll\restr{\Lebesgue}{I}$ with $\nu(I)<C_\nu$, then the reverse implications also hold.
\end{lemma}
\begin{proof}
\noindent
    \begin{itemize}
        \item[\enquote{$\emptyconvergencearrow{\tilde{d}_{\nu,p}}\Rightarrow\emptyconvergencearrow{\nu}$}:] Let $G_n\emptyconvergencearrow{\tilde{d}_{\nu,p}}G$, {i.e.}, $\tilde{d}_{\nu,p}(G_n,G)\nto0$, and take $\epsilon>0$. Then we have
              \begin{align*}
                  (\nu(\tilde{d}(G_n,G)>\epsilon))^{\frac1p}
                   & =\frac{1}{\epsilon}\parenthesis*{\int_{\set{\tilde{d}(G_n,G)>\epsilon}}\epsilon^p\diff\nu}^{\frac1p}
                  \leq\frac1{\epsilon}\parenthesis*{\int (\tilde{d}(G_n,G))^p\diff\nu}^{\frac1p}                          \\
                   & =\frac1\epsilon \tilde{d}_{\nu,p}(G_n,G)\nto0.
              \end{align*}
        \item[\enquote{$\emptyconvergencearrow{\nu}\Rightarrow\emptyconvergencearrow{\mathcal{C}}$}:] The following argument is adapted from Jonathan Hole's answer \cite{jonathanholehIfSequenceMonotone2020}.
              Let the sequence $(G_n)_{n\in\N}\subset\increasingfunctionclasses(I)$ converge to $G\in\increasingfunctionclasses(I)$ in $\nu$ {w.r.t.\@} $\tilde{d}$.
              Now take a continuity point $x\in I$ of $G$ and assume that $G_n(x)$ does not converge to $G(x)$ in $\tilde{d}$. Then there exist an $\epsilon>0$ and a subsequence $(G_{n_j})_{j\in\N}$ such that
              \begin{equation*}
                  \tilde{d}(G_{n_j}(x),G(x))>\epsilon
              \end{equation*}
              for each $j\in \N$.
              By continuity of $G$ at $x$, take a $\delta>0$ such that for all $y\in(x-\delta,x+\delta)$, we have
              \begin{equation*}
                  \tilde{d}({G(x),G(y)})\leq\frac\epsilon2.
              \end{equation*}
              For each $j\in\N$, we have either $G_{n_j}(x)>G(x)$ or $G_{n_j}(x)<G(x)$. This implies that either $G_{n_j}(x)>G(y)$, as otherwise
              \begin{equation*}
                  \tilde{d}(G_{n_j}(x),G(x))\leq \tilde{d}(G(x),G(y))\leq\frac\epsilon2,
              \end{equation*}
              or $G_{n_j}(x)<G(y)$, respectively, for each $y\in(x-\delta,x+\delta)$.
              In the first case, we thus have $G_{n_j}(y)\geq G_{n_j}(x)>G(y)$ for all $y\in(x,x+\delta)$, and thus
              \begin{equation*}
                  \tilde{d}(G_{n_j}(y),G(y))\geq \tilde{d}(G_{n_j}(x),G(y))\geq \tilde{d}(G_{n_j}(x),G(x))-\tilde{d}(G(x),G(y))>\frac{\epsilon}{2}.
              \end{equation*}
              In the other case, a similar calculation can be carried out for $y\in(x-\delta,x)$. Therefore,
              \begin{equation*}
                  \nu\parenthesis*{\tilde{d}(G_{n_j},G)>\frac{\epsilon}{2}}\geq\min(\nu((x,x+\delta)),\nu((x-\delta,x)))>0,\quad j\in\N.
              \end{equation*}
              This contradicts convergence in $\nu$. Thus, $\tilde{d}(G_n(x),G(x))\nto 0$ for all continuity points $x$ of $G$.
        \item[\enquote{$\emptyconvergencearrow{\mathcal{C}}\Rightarrow\emptyconvergencearrow{\text{$\Lebesgue$-{a.e.\@}\nocorr}}$}:]
              A non-decreasing function has only countably many discontinuities, and thus convergence at all continuity points implies $\Lebesgue$-{a.e.\@} convergence.

    \end{itemize}
    Assume now that $\nu\ll\restr{\Lebesgue}{I}$, $\nu(I)\leq C_\nu\in\R_{>0}$ and $\tilde{d}=\bar{d}$.
    \begin{itemize}
        \item[\enquote{$\emptyconvergencearrow{\text{$\Lebesgue$-{a.e.\@}\nocorr}}\Rightarrow\emptyconvergencearrow{\nu}$}:]
              This already holds for all measurable functions from $I$ to $\Rbar$ by \cite[Lemma 4.2.]{kallenbergFoundationsModernProbability2002}.
        \item[\enquote{$\emptyconvergencearrow{\nu}\Rightarrow\emptyconvergencearrow{\bar{d}_{\nu,p}}$}:]
              For each $\epsilon>0$, there exists a $\delta>0$ such that for sufficiently large $n$, the following holds:
              \begin{align*}
                  \bar{d}_{\nu,p}(G_n,G)^p
                   & =\int \bar{d}(G_n,G)^p\diff\nu
                  \leq \int_{\set{\bar{d}(G_n,G)\leq\delta}}\delta^p\diff\nu +\int_{\set{\bar{d}(G_n,G)>\delta}}1\diff\nu \\
                   & \leq \delta^p C_\nu+\nu(\bar{d}(G_n,G)>\delta)\leq\epsilon.
              \end{align*}
               Specifically, take $\delta\defeq\sqrt[p]{\frac{\epsilon}{2C_\nu}}$ and $N$ such that for all $n>N$, we have $\nu(\bar{d}(G_n,G)>\delta)\leq\frac\epsilon{2}$.\qedhere
    \end{itemize}
\end{proof}

\begin{remark}
    The remarkable part of this lemma is that convergence in measure implies $\Lebesgue$-{a.e.\@} convergence in this setting of non-decreasing functions and a full support reference measure.
\end{remark}

\subsection{The Borel structure of non-decreasing functions}
\label{sec: The Borel structure of non-decreasing functions}

We conclude from \autoref{lem: convergence on non-decreasing functions} that $\increasingfunctionclasses(I)=\increasingfunctionclasses^{1}(\rho,\bar{d})$ for $\rho=\restr{\gaussian(0,1)}{I}$ is a metric space with metric $\bar{d}_{\rho,1}$.
Denote the associated topology generated by open metric balls by
\begin{equation*}
    \bar{\tau}\defeq\tau(\bar{d}_{\rho,1})\defeq\tau({\set{B_{r,\bar{d}_{\rho,1}}(G)\given G\in \increasingfunctionclasses(I),r\in\R_{>0}}}).
\end{equation*}
This topology is sequential.
On the subspaces $\increasingfunctionclasses^{p}(\nu,\tilde{d})$, we have the metric topology $\tau(\tilde{d}_{\nu,p})$ as well as the subspace topology induced by $\bar{\tau}$, which we again denote by $\bar{\tau}$. Clearly, any $\bar{\tau}$-closed set is $\tau(\tilde{d}_{\nu,p})$-closed by \autoref{lem: convergence on non-decreasing functions}, and thus $\bar{\tau}\subset\tau(\tilde{d}_{\nu,p})$.
Based on these two topologies, we obtain two $\sigma$-algebras on $\increasingfunctionclasses^{p}(\nu,\tilde{d})$: the Borel $\sigma$-algebra $\B(\tilde{d}_{\nu,p})$ and $\sigma(\bar{\tau})$. We additionally consider the initial $\sigma$-algebra $\sigma(\iota_I)$ generated by the evaluation maps $\iota_x$ for $x\in I$. We demonstrate the equality of these $\sigma$-algebras for bounded $I$ and finite $\nu$.
\begin{lemma}
    \label{lem: equivalence of sigma algebras on non-decreasing functions}
    For a nonvoid, bounded, open interval $I=(a,b)$, $a<b\in\R$, $\restr{\Lebesgue}{I}\gg\nu\in\locallyfinitemeasureon(I)$ with full support and $\nu(I)\in\R$, $p\in[1,\infty)$, and $\tilde{d}\in\set{d,\bar{d}}$, it holds on $\increasingfunctionclasses^{p}(\nu,\tilde{d})$ that $\sigma(\iota_{I})=\sigma(\bar{\tau})=\B(\tilde{d}_{\nu,p})$.
\end{lemma}
\begin{proof}
    \noindent
    \begin{itemize}
        \item[\enquote{$\sigma(\iota_{I})\subset\sigma(\bar{\tau})$}:]
              \cite[Proof of Theorem 12.5 (ii)]{billingsleyConvergenceProbabilityMeasures1999}
              For each $x\in I$ and $k\in\N$ with $x+\frac1k\in I$, define $h_{k,x}:\increasingfunctionclasses^{p}(\nu,\tilde{d})\to\R$ by
              \begin{equation*}
                  h_{k,x}(G)=k\int_x^{x+\frac1k} G(t)\diff t.
              \end{equation*}
              Each $h_{k,x}:(\increasingfunctionclasses^{p}(\nu,\tilde{d}),\bar{\tau})\to(\R,d)$ is continuous. Let $G_n\emptyconvergencearrow{\bar{d}_{\rho,1}} G$. Then, in particular, $G_n\emptyconvergencearrow{\text{$\Lebesgue$-{a.e.\@}\nocorr}} G$. Therefore,
              \begin{equation*}
                  \lim_{n\to\infty} h_{k,x}(G_n)=k\int_x^{x+\frac1k} \lim_{n\to\infty}G_n(t)\diff t = h_{k,x}(G)
              \end{equation*}
              by dominated convergence: Choose continuity points $r<x<x+\frac1k<s$ in $I$. Then, since $G_n\emptyconvergencearrow{\mathcal{C}} G$, for sufficiently large $n$ and some $\epsilon>0$, we have
              \begin{equation*}
              G(r)-\epsilon\leq G_n(x)\leq G_n(y) \leq G_n(x+\frac1k) \leq G(s)+\epsilon,\quad y\in(x,x+\frac1k).
              \end{equation*}
              Thus, by the quasi-convexity of the absolute value, we have
              \begin{equation*}
                \abs{G_n(y)}\leq \max(\abs{G(r)-\epsilon},\abs{G(s)+\epsilon})\leq \epsilon + \max(\abs{G(r)},\abs{G(s)}),\quad y\in(x,x+\frac1k).
              \end{equation*}

              With this, we have in particular established the $\sigma(\bar{\tau}),\B(\R)$-measurability of $h_{k,x}$.
              We now show that $h_{k,x}$ converges pointwise to $\iota_x$ on $\increasingfunctionclasses^{p}(\nu,\tilde{d})$. Take an arbitrary $G\in\increasingfunctionclasses^{p}(\nu,\tilde{d})$. By the right-continuity of $G$, for each $\delta>0$, there exists $M_\delta$ such that for all $k>M_{\delta}$, we have $G(x+\frac1k)\leq G(x)+\delta$. Thus,
              \begin{align*}
                                                          && G(x)\leq  \lim_{k\to\infty} h_{k,x}(G) & \leq\lim_{k\to\infty,k>M_\delta}k\int_{x}^{x+\frac1k}G(x)+\delta\diff t=G(x)+\delta \\
                  \overset{\delta\searrow 0}{\Rightarrow} &&           \lim_{k\to\infty} h_{k,x}(G) & =G(x)=\iota_x(G).
              \end{align*}
              As the pointwise limit of measurable functions, $\iota_x$ is $\sigma(\bar{\tau}),\B(\R)$-measurable for each $x\in I$.
              Therefore, $\sigma(\bar{\tau})\supset\sigma(\iota_{I})$.
        \item[\enquote{$\sigma(\bar{\tau})\subset\B(\tilde{d}_{\nu,p})$}:]
              This holds as $\tau(\bar{d}_{\rho,1})\subset\tau(\tilde{d}_{\nu,p})$.
        \item[\enquote{$\B(\tilde{d}_{\nu,p})\subset\sigma(\iota_{I})$}:]
              We now prove that $\Id:(\increasingfunctionclasses^{p}(\nu,\tilde{d}),\sigma(\iota_{I}))\to(\increasingfunctionclasses^{p}(\nu,\tilde{d}),\B(\tilde{d}_{\nu,p}))$ is measurable to obtain the remaining inclusion.
              To this end, we show that we can approximate the identity by measurable projections onto simple non-decreasing functions.
              For $n\in\N$, denote by $x^{(n)}\in [a,b]^{2^n+1}$ the finite partition of $I=(a,b)$ with
              $x^{(n)}_0=a$, $x^{(n)}_{k+1}-x^{(n)}_k=\Delta_n$ and $x^{(n)}_{2^n}=b$, $\Delta_n=\frac{b-a}{2^n}$.
              Define $\pi_n:\increasingfunctionclasses^{p}(\nu,\tilde{d})\to\R^{2^n-1}$ by $(\pi_n(G))_k=\iota_{x^{(n)}_{k}}(G)$ for $k\in\set{1,\dots,2^n-1}$.
              By the definition of $\sigma(\iota_I)$, $\pi_n$ is $\sigma(\iota_{I}),\B(\R^{2^n-1})$-measurable for each $n\in\N$. Now define the subpartition
              \begin{equation*}
                  y_k^{(n)}=\begin{cases}
                      x_{k}^{(n)}   & k< 2^{n-1}    \\
                      x_{k+1}^{(n)} & k\geq 2^{n-1}
                  \end{cases},\quad k\in\set{0,\dots 2^n-1}
              \end{equation*}
              and $\Pi_n:\Functionimage{\pi_n}{\increasingfunctionclasses^{p}(\nu,\tilde{d})}\to\increasingfunctionclasses^{p}(\nu,\tilde{d})$ by
              \begin{equation*}
                  \Pi_n\parenthesis*{(z_k)_{k=1}^{2^n-1}}(x)
                  =\sum_{k=1}^{2^n-1}z_k\1_{[y_{k-1}^{(n)},y^{(n)}_{k})}(x),\quad x\in I.
              \end{equation*}
              For each $n\in\N$, we can define a metric $d^{(n)}$ on $\Functionimage{\pi_n}{\increasingfunctionclasses^{p}(\nu,\tilde{d})}$ that generates the subspace topology of the product topology on $\R^{2^n-1}$ by
              \begin{equation*}
                  d^{(n)}(v,z)=\parenthesis*{\sum_{k=1}^{2^n-1}(\tilde{d}(v_k,z_k))^p\nu([y^{(n)}_{k-1},y^{(n)}_{k}))}^{\frac{1}{p}}=\tilde{d}_{\nu,p}(\Pi_n(v),\Pi_n(z))
              \end{equation*}
              as $\nu$ has full support.
              We have $\sigma(d^{(n)})=\restr{\B(\R^{2^n-1})}{\Functionimage{\pi_n}{\increasingfunctionclasses^{p}(\nu,\tilde{d})}}$.
              This choice of metric turns $\Pi_n$ into an isometric embedding of $(\Functionimage{\pi_n}{\increasingfunctionclasses^{p}(\nu,\tilde{d})},d^{(n)})$ into $(\increasingfunctionclasses^{p}(\nu,\tilde{d}),\tilde{d}_{\nu,p})$, so it is continuous and thus $\sigma(d^{(n)}),\B(\tilde{d}_{\nu,p})$-measurable.
              Therefore, the functions $\xi_n\defeq\Pi_n\circ \pi_n$ are $\sigma(\iota_I),\B(\tilde{d}_{\nu,p})$-measurable for each $n$.
              Let $G\in\increasingfunctionclasses^{p}(\nu,\tilde{d})$ and let $t$ be an arbitrary continuity point.
              Define $t_n=x_{k'}^{(n)}$ such that $\xi_n(G)(t)\eqdef G(t_n)$.
              Then $t_n\nto t$ by construction of $x^{(n)}$, and therefore $\xi_n(G)(t)= G(t_n)\nto G(t)$.
              Furthermore, by construction, $\tilde{d}(\xi_n(G)(x),G(x))$ is monotone decreasing for each $x$, and thus
              \begin{equation*}
                \tilde{d}_{\nu,p}(\xi_n(G),G)^p\leq \tilde{d}_{\nu,p}(\xi_1(G),G)^p=\tilde{d}_{\nu,p}\parenthesis*{G\parenthesis*{\frac{a+b}{2}},G}^p
              \end{equation*}
              with the latter being $\nu$-integrable as $G\in\increasingfunctionclasses^{p}(\nu,\tilde{d})$ and $\nu$ is a finite measure. Thus, dominated convergence yields
              \begin{align*}
                  \lim_{n\to\infty}(\tilde{d}_{\nu,p}(\xi_n(G),G))^p
                   & =\lim_{n\to\infty}\int_{a}^b\tilde{d}(\xi_n(G)(t),G(t))^p\diff\nu(t)  \\
                   & = \int_{a}^b\lim_{n\to\infty}\tilde{d}(\xi_n(G)(t),G(t))^p\diff\nu(t)
                  =0.
              \end{align*}
              Thus, $\xi_n(G)\emptyconvergencearrow{\tilde{d}_{\nu,p}}G=\Id(G)$ for all $G\in\increasingfunctionclasses^{p}(\nu,\tilde{d})$, {i.e.}, $\xi_n\nto\Id$ on $(\increasingfunctionclasses^{p}(\nu,\tilde{d}),\tilde{d}_{\nu,p})$, and therefore $\Id$ is $\sigma(\iota_I),\B(\tilde{d}_{\nu,p})$-measurable as the limit of measurable functions in a metric space \cite[Lemma 1.10.]{kallenbergFoundationsModernProbability2002}.
              \qedhere
    \end{itemize}
\end{proof}

\subsection{Measurability and continuity on the set of locally finite measures}

    Denote by $\locallyfinitemeasureon(I)$ the set of locally finite measures on $I=(a,b)\subset\R$, as a subset of $\measureon(I)$, the set of all Borel measures on $I$.
    Define on it the \emph{vague topology} $\tau_{\vagueconvergencesymbol}$ as the initial topology that makes all integral evaluations $e_{C_c(I)}$ at continuous functions with compact support continuous.
    Here, for each $f\in C_c(I)$, we define $e_f(\nu)\defeq\int f \diff\nu$, $\nu\in\locallyfinitemeasureon(I)$.
    This topology $\tau_{\vagueconvergencesymbol}=\tau(e_{C_c(I)})$ is Polish \cite[4.2]{kallenbergFoundationsModernProbability2002} and thus sequential \cite[Proposition 3.1]{basrakNoteVagueConvergence2019}.
    A sequence of measures $(\mu_n)_{n\in\N}\subset\locallyfinitemeasureon(I)$ \emph{converges vaguely} to $\mu\in\locallyfinitemeasureon(I)$, $\mu_n\vagueto\mu$, if $e_f(\mu_n)\nto e_f(\mu)$ for all $f\in C_c(I)$.
    Furthermore, the Borel $\sigma$-algebra $\B(\locallyfinitemeasureon(I))\defeq\sigma(\tau_{\vagueconvergencesymbol})$ is equivalently generated by all evaluations, {i.e.}, $\B(\locallyfinitemeasureon(I))=\sigma(\iota_{\B(I)})$ \cite[Theorem 4.2]{kallenbergRandomMeasuresTheory2017}.
    Similarly, we define on $\distributionon(\R)$ the \emph{weak topology} $\tau_{\weakconvergencesymbol}=\tau(e_{C_b(\R)})$ as the initial topology that makes all integral evaluations $e_f$ at bounded continuous functions $f\in C_b(\R)$ continuous.
    For finite measures $\nu$, $e_{f}(\nu)=\int f\diff\nu$ is well-defined.
    Again, the weak topology is sequential, and a sequence $(\mu_n)_{n\in\N}\subset\distributionon(\R)$ \emph{converges weakly} to $\mu\in\distributionon(\R)$ if $e_f(\mu_n)\nto e_f(\mu)$ for all $f\in C_b(\R)$.
    Again, we have $\sigma(\tau_{\weakconvergencesymbol})=\sigma(\iota_{\B(\R)})$ by the same argument \cite[Theorem 4.2]{kallenbergRandomMeasuresTheory2017}, as weak convergence of probability distributions is the restriction of a generalization of vague convergence to finite measures based on a different localizing structure.
    If we instead restrict further to $\Wp(\R)$, then convergence in $W_{\!p}$ implies weak convergence \cite{Villani.2009} and thus $\tau_{\weakconvergencesymbol}\subset\tau(W_{\!p})$. Therefore, $\sigma(\iota_{\B(\R)})=\sigma(\tau_{\weakconvergencesymbol})\subset\B(\Wp(\R))$ for the subspace topology and $\sigma$-algebra, respectively.
    That the reverse inclusion holds is shown directly in \cite[Lemma 2.4.2]{Panaretos.2020} or in \cite[Lemma 5.8]{maAbsoluteContinuityWasserstein2025} using Souslin theory for the case $p=2$, with the proof directly applicable to arbitrary $p$.

\subsection{Continuity and measurability of the relations between measures and non-decreasing functions}
\label{sec: Continuity and measurability of the relations between measures and non-decreasing functions}

We can now turn to the continuity of the functions defined in \autoref{sec: non-decreasing functions mappings}. While we fix $I=(0,1)$ for this discussion, these results can be generalized directly to bounded $I$.
To do so, one has to replace $\distributionon(\R)$ with the appropriate set of finite measures.
\begin{lemma}\label{lem: quantile transform homeomorphism Stieltjes transform continuous}
    We have that
    \begin{enumerate}
        \item
              $\CDFtransform:(\distributionon(\R),\tau_{\weakconvergencesymbol})\to(\increasingfunctionclasses^{*}(0,1),\bar{\tau})$ is a
              homeomorphism,
        \item $\generalizedinverse:(\increasingfunctionclasses(0,1),\bar{\tau})\to(\increasingfunctionclasses^{*}(0,1),\bar{\tau})$ is a homeomorphism,
        \item
              $\Quantiletransform:(\distributionon(\R),\tau_{\weakconvergencesymbol})\to(\increasingfunctionclasses(0,1),\bar{\tau})$ is a
              homeomorphism, and
        \item
              $\Stieltjestransform:(\increasingfunctionclasses(0,1),\bar{\tau})\to(\locallyfinitemeasureon(0,1),\tau_{\vagueconvergencesymbol})$ is continuous.
    \end{enumerate}
\end{lemma}
\begin{delayedproof}{lem: quantile transform homeomorphism Stieltjes transform continuous}
    Bijectivity of the first three functions and surjectivity of the last function have already been discussed. It remains to show the continuity properties.
    \begin{enumerate}
        \item
              It is a classical result that a sequence of distributions $(\mu_n)_{n\in\N}\subset\distributionon(\R)$ converges weakly to a distribution $\mu\in\distributionon(\R)$ if and only if $F_{\mu_n}\emptyconvergencearrow{\mathcal{C}} F_\mu$, {i.e.}, at continuity points of $F_\mu$ \cite[13.23]{klenkeProbabilityTheoryComprehensive2008}. This completes the proof, as by \autoref{lem: convergence on non-decreasing functions}, convergence at continuity points is equivalent to convergence in the sequential topology $\bar{\tau}$ on $\increasingfunctionclasses(\R)$, and $\increasingfunctionclasses^{*}(0,1)\subset\increasingfunctionclasses(\R)$.
        \item
              \label{enum Proof utilizing elementary Skorokhod theorem}
              The elementary Skorokhod theorem \cite[Lemma A.5]{bobkovOnedimensionalEmpiricalMeasures2019} now gives us that this convergence at continuity points of $F_\mu$ implies $Q_{\mu_n}(t)\nto Q_\mu(t)$ for all continuity points $t$ of $Q_\mu$. The proof only uses the fact that $Q_{\mu}=\generalizedinverse(F_\mu)$ for all $\mu\in\distributionon(\R)$, and thus the reverse direction also holds \cite[Exercise 5.]{shorackEmpiricalProcessesApplications1986}. Our right-continuous {q.f.\@} convention causes no difficulties here, as left- and right-continuous versions coincide at their shared continuity points. Thus, the mapping $\generalizedinverse(F_\mu)=Q_\mu$ is a homeomorphism.
        \item
              As $\Quantiletransform=\generalizedinverse\circ\CDFtransform$, this composition is also a homeomorphism.
        \item
              We show that $\Stieltjestransform{Q_n}\emptyconvergencearrow{\vagueconvergencesymbol}\Stieltjestransform{Q}$ if $(Q_n)_{n\in\N}\subset\increasingfunctionclasses(0,1)$ converges to $Q\in\increasingfunctionclasses(0,1)$ at continuity points. We already know from \ref{enum Proof utilizing elementary Skorokhod theorem} that $F_n=\generalizedinverse(Q_n)\emptyconvergencearrow{\mathcal{C}} F=\generalizedinverse(Q)$.
              Let $\chi\in C_c(0,1)$. There exists a compact interval $[x,y]\subset(0,1)$ such that $\restr{\chi}{[x,y]^c}\equiv 0$ and $Q$ is continuous at $x,y$. By continuity, we can choose a compact interval $K\subset\R$ such that $[(Q_n)_l(x),Q_n(y)]\subset K$ for all sufficiently large $n$. Therefore, we have
              \begin{align*}
                  \lim_{n\to\infty} \int \chi\diff\Stieltjestransform{Q_n}
                                               & = \lim_{n\to\infty}\int_{[x,y]} \chi\diff\Transportmeasure{F_n}
                  \overset{\eqref{eqn: transport formulation of inverse of Quantile transform}}{=} \lim_{n\to\infty}\int_{[x,y]} \chi\diff\pushforward{{F_n}}{\Lebesgue}                                                                        \\
                                               & =
                  \lim_{n\to\infty}\int_{\preimage{F_n}{[x,y]}} \chi(F_n(u))\diff\Lebesgue(u)                                                                    \\
                                               & =                          \lim_{n\to\infty}\int\1_{{[(Q_n)_l(x),Q_n(y)]}}(u) \chi(F_n(u))\diff\Lebesgue(u) \\
                  \overset{\chi\in C_c(0,1)} &
                  =\int_K \lim_{n\to\infty}\1_{{[(Q_n)_l(x),Q_n(y)]}}(u)\chi(F_n(u))\diff\Lebesgue(u)                                                        \\
                                               & = \int \1_{[Q(x),Q(y)]}(u)\chi(F(u))\diff\Lebesgue(u)
                  = \int_{[x,y]} \chi\diff\pushforward{F}{\Lebesgue}                                                                                             \\
                                               & = \int_{[x,y]} \chi\diff\Transportmeasure(F)
                  = \int \chi\diff\Stieltjestransform(Q).\qedhere
              \end{align*}
    \end{enumerate}
\end{delayedproof}
Therefore, \autoref{lem: quantile transform homeomorphism Stieltjes transform continuous} also establishes measurability {w.r.t.\@} the $\sigma$-algebras initial to the corresponding evaluation maps. Thus, we can freely switch between representations of random elements on these spaces:
\begin{delayedproof}{lem: equivalence of kernels}
    By \autoref{lem: quantile transform homeomorphism Stieltjes transform continuous}, $\Quantiletransform:\distributionon(\R)\to\increasingfunctionclasses(0,1)$ is a homeomorphism with respect to weak convergence and convergence at continuity points, respectively, and $\Stieltjestransform:\increasingfunctionclasses(0,1)\to\locallyfinitemeasureon(0,1)$ is continuous {w.r.t.\@} convergence at continuity points and vague convergence, respectively.
    By \cite{kallenbergRandomMeasuresTheory2017}, we know that $\sigma(\tau_{\weakconvergencesymbol})=\sigma(\iota_{\B(\R)})$ on $\distributionon(\R)$ and $\sigma(\tau_{\vagueconvergencesymbol})=\sigma(\iota_{\B(0,1)})$ on $\locallyfinitemeasureon(0,1)$.
    Furthermore, we have seen that $\sigma(\bar{\tau})=\sigma(\iota_{(0,1)})$ by \autoref{lem: equivalence of sigma algebras on non-decreasing functions} applied to $\increasingfunctionclasses^{1}(\uniformmeasure,\bar{d})=\increasingfunctionclasses(0,1)$.
    Therefore, $\Quantiletransform$ is $\sigma(\iota_{\B(\R)}),\sigma(\iota_{(0,1)})$-bimeasurable.
    Thus, $\randommeasure:(\Omega,\A)\to(\distributionon(\R),\sigma(\iota_{\B(\R)}))$ is measurable if and only if $\Quantiletransform\circ\randommeasure$ is measurable.
    Moreover, $\InverseMeasure\circ\randommeasure=\Stieltjestransform\circ\Quantiletransform\circ\randommeasure$ is $\A,\sigma(\iota_{\B(0,1)})$-measurable if $\randommeasure$ is $\A,\sigma(\iota_{\B(\R)})$-measurable.
\end{delayedproof}

\subsection{The Borel \texorpdfstring{$\sigma$}{σ}-algebra on \texorpdfstring{$\Wp(\R)$}{𝒲ₚ(ℝ)}}
    \label{proof: thm: measurable structure on Wasserstein space}
    Analogously to \autoref{lem: quantile transform homeomorphism Stieltjes transform continuous}, we have that
    \begin{enumerate}
        \item $\Quantiletransform:(\Wp(\R),\tau(W_{\!p}))\to(\increasingfunctionclasses^{p}(0,1),\tau(d_{\uniformmeasure,p}))$ is a homeomorphism, as it is an isometry,
        \item and $\Stieltjestransform:(\increasingfunctionclasses^{p}(0,1),\tau(d_{\uniformmeasure,p}))\to(\locallyfinitemeasureon(0,1),\tau_{\vagueconvergencesymbol})$ is still continuous by \autoref{lem: quantile transform homeomorphism Stieltjes transform continuous} and the inclusion $\bar{\tau}\subset\tau(d_{\uniformmeasure,p})$.
    \end{enumerate}

    We now use \autoref{lem: equivalence of sigma algebras on non-decreasing functions}, {i.e.}, $\sigma(\bar{\tau})=\B(d_{\uniformmeasure,p})=\sigma(\iota_{(0,1)})$, to characterize the Borel structure on $\Wp(\R)$. As mentioned before, this special case also follows directly from \cite[Lemma 2.4.2]{Panaretos.2020} and \cite[Theorem 4.2]{kallenbergRandomMeasuresTheory2017}.
    \begin{lemma}
        \label{Borel sigma algebra on Wasserstein p on R}
        We have $\B(\Wp(\R))=\sigma(\iota_{\B(\R)})$ on $\Wp(\R)$ for all $p\geq1$.
    \end{lemma}
    \begin{proof}
        We have $\sigma(\iota_{\B(\R)})=\sigma(\tau_{\weakconvergencesymbol})$ by \cite{kallenbergRandomMeasuresTheory2017}.
        It is clear that convergence in $W_{\!p}$ implies weak convergence \cite[Theorem 6.9]{Villani.2009}. Therefore, $\tau_{\weakconvergencesymbol}\subset\tau(W_{\!p})$ for the subspace $\sigma$-algebra.
        We thus have
        \begin{equation*}
            \sigma(\iota_{\B(\R)})=\sigma(\tau_{\weakconvergencesymbol})\subset\sigma(\tau({W_{\!p}}))=\B(\Wp(\R)).
        \end{equation*}
        It remains to show that $\Id:(\Wp(\R),\sigma(\tau_\weakconvergencesymbol))\to(\Wp(\R),\B(\Wp(\R)))$ is measurable.
        We conclude from \autoref{lem: quantile transform homeomorphism Stieltjes transform continuous} that $\Quantiletransform:\Wp(\R)\to\increasingfunctionclasses^{p}(0,1)$ is $\sigma(\tau_\weakconvergencesymbol),\sigma(\bar{\tau})$-measurable.
        By \autoref{lem: equivalence of sigma algebras on non-decreasing functions}, we have that $\sigma(\bar{\tau})=\B(\increasingfunctionclasses^{p}(0,1))$ and thus $\Quantiletransform$ is $\sigma(\tau_\weakconvergencesymbol),\B(\increasingfunctionclasses^{p}(0,1))$-measurable.
        Furthermore, since $\Transportmeasure:\increasingfunctionclasses^{p}(0,1)\to\Wp(\R)$ is $\B(\increasingfunctionclasses^{p}(0,1)),\B(\Wp(\R))$-measurable, we conclude that $\Id_{\Wp(\R)}=\Transportmeasure\circ\Quantiletransform$ is $\sigma(\tau_\weakconvergencesymbol),\B(\Wp(\R))$-measurable as a composition of measurable functions. Thus, we obtain the desired equality $\sigma(\iota_{\B(\R)})=\sigma(\tau_\weakconvergencesymbol)=\B(\Wp(\R))$.
    \end{proof}
    \begin{delayedproof}{thm: measurable structure on Wasserstein space}
        Take $p=2$ in \autoref{Borel sigma algebra on Wasserstein p on R}.
    \end{delayedproof}

\section{Further proofs}
\label{sec: Proofs of main paper}

\subsection{Proof of \autoref{thm: equivilance between M concave and quantile M convex quasi arithmetic}}
    \label{sec: Proof of general generalized inverse pair theorem on equivalenve of derivatives}

    We call a function $f:I\to\R$ \emph{locally absolutely continuous} if, for all compact intervals $[a,b]\subset I$, the restriction $\restr{f}{(a,b)}$ is absolutely continuous in the sense of Banach, {i.e.,} for each $\epsilon>0$ there exists a $\delta>0$ such that for any disjoint sequence of open subintervals $((a_j,b_j))_{j\in N}$ of $[a,b]$ with a possibly countable index set $N\subseteq\N$, it holds that
    \begin{equation*}
        \sum_{j\in N} b_j-a_j<\delta\Rightarrow\sum_{j\in N} \abs{f(b_j)-f(a_j)}<\epsilon.
    \end{equation*}
    A function $f':I\to\Rbar$ is called the \emph{Radon-Nikodym derivative} of a locally absolutely continuous function $f:I\to\R$ if $f'$ coincides with the local derivative of $f$ wherever $f$ is differentiable, which holds $\Lebesgue$-{a.e.\@} on $I$.
    In the case of non-decreasing functions, these definitions coincide with those given in the main text.
    We will use the following theorems to deal with the absolute continuity of compositions.
    \begin{corollary}[{One-dimensional version of \cite[Theorem 4.2]{helselAbsolutelyContinuousProduct1952}}]
        \label{theo: composition of absolutely continuous functions}
        A composition $T_2\circ T_1$ of absolutely continuous functions $T_1:I_1\to\R$, $T_2:I_2\to\R$ defined on open bounded intervals $I_1,I_2\subset\R$ is absolutely continuous if and only if $\Functionimage{T_1}{I_1}\subset I_2$ and, for two Radon-Nikodym derivatives $\RadonNikodym{T_1}{\Lebesgue},\RadonNikodym{T_2}{\Lebesgue}$, the function $\RadonNikodym{T_2}{\Lebesgue}\circ T_1\cdot\RadonNikodym{T_1}{\Lebesgue}:I_1\to\Rbar$ is integrable.\footnote{For products of Radon-Nikodym derivatives we use the convention $0\cdot(\pm\infty)=0$.} If this composition is absolutely continuous, then we always have $\RadonNikodym{T_2\circ T_1}{\Lebesgue}=\RadonNikodym{T_2}{\Lebesgue}\circ T_1\cdot\RadonNikodym{T_1}{\Lebesgue}$ $\restr{\Lebesgue}{I_1}$-{a.e}.
    \end{corollary}
    \begin{theorem}[{\cite[IX 1. Satz 3]{natansonTheorieFunktionenReellen1961}}]
        \label{theo: composition of absolute continuous and absolute continuous non-decreasing function}
        A well-defined composition $T_2\circ T_1$ of absolutely continuous functions $T_1:I_1\to\R$, $T_2:I_2\to\R$ defined on open bounded intervals $I_1,I_2\subset\R$ is absolutely continuous if $T_1$ is additionally strictly increasing.
    \end{theorem}
    We call a function $f:I\to\R$ on an open interval $I$ \emph{locally Lipschitz} if for all compacta $[a,b]\subset I$ we have that $\restr{f}{[a,b]}$ is Lipschitz, {i.e.}, there exists a constant $L_{[a,b]}$ such that
            \begin{equation*}
                \abs{f(x)-f(y)}\leq L_{[a,b]}\abs{x-y},\quad x,y\in [a,b].
            \end{equation*}
    We can now formulate the central lemma for the proof of \autoref{thm: equivilance between M concave and quantile M convex quasi arithmetic}.
    \begin{lemma}
        \label{lem: derivative property of generalized inverse and monotone transformation}
        Let $G,H$ be a generalized inverse pair such that $G$ and $H$ are locally absolutely continuous on the open intervals $I_G$, $I_H$, respectively, satisfying $I_G=\Functionimage{H}{I_H}$, $I_H=\Functionimage{G}{I_G}$. For every generator $\phi$ and Radon-Nikodym derivative $g$ of $G$, define $h\defeq\frac{1}{g\circ H}$, which is a Radon-Nikodym derivative of $H$. If $\phi\circ g$ is locally Lipschitz on $I_G$, then $\monotonetransformation{\phi}\circ h$ is locally Lipschitz on $I_H$ and
        \begin{equation*}
            \RadonNikodym{\monotonetransformation{\phi}\circ h}{\Lebesgue}  =-\RadonNikodym{\phi\circ g}{\Lebesgue}\circ H\quad\text{$\restr{\Lebesgue}{I_H}$-{a.e.\@}\nocorr}
        \end{equation*}
    \end{lemma}
    \begin{remark}
        The existence of $h=\frac{1}{g\circ H}$ on $I_H$ follows directly since $H$ is bi-absolutely continuous on $I_H$
        \cite[Proposition A.18]{bobkovOnedimensionalEmpiricalMeasures2019}.
    \end{remark}
    \begin{lemma}
        \label{lem: Welldefinedness of extension of integrated transform for ragular transformed density}
        Under the assumptions of \autoref{lem: derivative property of generalized inverse and monotone transformation}, there exists an extension of $\integratedgeneralizedinverse{\phi}:\Functionimage{\phi}{\R_{>0}}\to\R$ to an open interval $I_{\integratedgeneralizedinverse{\phi}}$ that is locally absolutely continuous and non-increasing, such that, using this extension, $\monotonetransformation{\phi}\circ h=\integratedgeneralizedinverse{\phi}\circ \phi\circ g \circ H$ holds and the right-hand side is real-valued and well-defined.
    \end{lemma}
    \begin{proof}
        A locally absolutely continuous and monotone function on a bounded interval is absolutely continuous if finite one-sided limits exist at the boundaries. If there is a finite one-sided limit, then we can extend the function in a locally absolutely continuous manner beyond that boundary by extending it constantly with the finite limit.
        For our extension of $\integratedgeneralizedinverse{\phi}$, we consider the values of $\integratedgeneralizedinverse{\phi}$ at the boundaries $\phi(0),\phi(+\infty)$, namely,
        $\integratedgeneralizedinverse{\phi}(\phi(0))\defeq\lim_{x\searrow \phi(0)}\integratedgeneralizedinverse{\phi}(x)$ and $\integratedgeneralizedinverse{\phi}(\phi(+\infty))\defeq\lim_{x\nearrow \phi(+\infty)}\integratedgeneralizedinverse{\phi}(x)$, or equivalently, the values $\monotonetransformation{\phi}(0)=\integratedgeneralizedinverse{\phi}(\phi(+\infty))$, $\monotonetransformation{\phi}(+\infty)=\integratedgeneralizedinverse{\phi}(\phi(0))$.
        We define $I_{\integratedgeneralizedinverse{\phi}}\eqdef(a_\phi,b_\phi)$ with
        \begin{align*}
            a_{\phi} & \defeq\begin{cases}
                                \phi(0) & \integratedgeneralizedinverse{\phi}(\phi(0))=+\infty  \\
                                -\infty & \integratedgeneralizedinverse{\phi}(\phi(0))<+\infty,
                            \end{cases}      \\
            b_{\phi} & \defeq\begin{cases}
                                \phi(+\infty) & \integratedgeneralizedinverse{\phi}(\phi(+\infty))=-\infty \\
                                +\infty       & \integratedgeneralizedinverse{\phi}(\phi(+\infty))>-\infty
                            \end{cases}
        \end{align*}
        and set $\integratedgeneralizedinverse{\phi}(x)=\integratedgeneralizedinverse{\phi}(\phi(0))$ for all $x\in(a_\phi,\phi(0)]$ and $\integratedgeneralizedinverse{\phi}(x)=\integratedgeneralizedinverse{\phi}(\phi(+\infty))$ for all $x\in[\phi(+\infty),b_\phi)$. Therefore, $\integratedgeneralizedinverse{\phi}:I_{\integratedgeneralizedinverse{\phi}}\to\R$ is real-valued and still non-increasing. It is also locally absolutely continuous since it is continuous, locally of bounded variation, and maps null sets to null sets.
        Recall that, by the assumptions of \autoref{lem: derivative property of generalized inverse and monotone transformation}, $\mathfrak{g}\defeq \phi\circ g$ is locally Lipschitz.
        To show that, with this extension, $\integratedgeneralizedinverse{\phi}\circ \phi\circ g\circ H=\integratedgeneralizedinverse{\phi}\circ\mathfrak{g}\circ H$ is well-defined, we need to check that $\Functionimage{\mathfrak{g}}{I_G}\subset I_{\integratedgeneralizedinverse{\phi}}$.
        \begin{itemize}
            \item[\enquote{$\Functionimage{\phi\circ g}{I_G}\subset(a_{\phi},+\infty)$}:]
                If $g$ does not attain a global minimum on $I_G$ or if that minimum is greater than $0$, then we have that $a_\phi\leq \phi(0)<\Functionimage{\mathfrak{g}}{I_G}$.
                Thus, assume that there exists $x\in I_G$ such that $g$ attains its minimum $0$ at $x$. Then necessarily $\phi(0)\neq-\infty$, since otherwise $\abs{\mathfrak{g}(x)}=+\infty$ on $I_G$, contradicting the local Lipschitz continuity of $\mathfrak{g}$. We can calculate
                \begin{equation*}
                    \integratedgeneralizedinverse{\phi}(\phi(0))=\lim_{x\searrow \phi(0)}\int_x^{\phi(\frac1c)}\phi^{-1}\diff\Lebesgue+\phi(c)\leq\int_{\phi(0)}^{\phi(\frac1c)}\frac1c\diff\Lebesgue+\phi(c)=\frac{\phi(\frac1c)-\phi(0)}{c} +\phi(c)<+\infty.
                \end{equation*}
                Therefore, $a_\phi=-\infty<\phi(0)\leq\Functionimage{\mathfrak{g}}{I_G}$.
            \item[\enquote{$\Functionimage{\phi\circ g}{I_G}\subset(-\infty,b_\phi)$}:]
                Analogously, if $g$ does not attain a global maximum on $I_G$ or if that maximum is less than $+\infty$, then we have that $b_\phi\geq \phi(+\infty)>\Functionimage{\mathfrak{g}}{I_G}$.
                Thus, assume that there exists $x\in I_G$ such that $g$ attains its maximum $+\infty$ at $x$. Then again $\phi(+\infty)<+\infty$, since $\mathfrak{g}$ is locally Lipschitz. As $\phi^{-1}$ is unbounded as $x\nearrow\phi(+\infty)$, we need a different approach than before.
                There exists $\epsilon>0$ such that $[x-\epsilon,x+\epsilon]\subset I_G$.
                Since $\mathfrak{g}$ is locally Lipschitz, there exists $L\defeq L_{[x-\epsilon,x+\epsilon]}\in\R_{>0}$ such that
                \begin{equation*}
                    \mathfrak{g}(x+z)\geq \mathfrak{g}(x)- L\abs{z}=\phi(+\infty)- L\abs{z},\quad z\in[-\epsilon,\epsilon].
                \end{equation*}
                Take now $0<\delta<\epsilon$ such that $\mathfrak{g}(x)-L\delta>\phi(0)$.
                Because $g=\phi^{-1}\circ \mathfrak{g}$ is a density of $G$, we obtain, after applying $\phi^{-1}$ to both sides and integrating, that
                \begin{align*}
                    +\infty & >\int_0^{\delta}g(x+z)\diff\Lebesgue(z)=\int_{0}^{\delta}\phi^{-1}\circ\mathfrak{g}(x+y)\diff\Lebesgue(y)
                    \geq\int_{0}^\delta\phi^{-1}(\phi(+\infty)- L\abs{y})\diff\Lebesgue(y)                                                                                                       \\
                            & =\int^{\phi(+\infty)}_{\phi(+\infty)-L\delta}\frac{\phi^{-1}(u)}{L}\diff\Lebesgue(u)
                    =\frac{1}{L}\parenthesis*{\int^{\phi(+\infty)}_{\phi(\frac1c)}\phi^{-1}(u)\diff\Lebesgue(u)-\int_{\phi(\frac{1}{c})}^{\phi(+\infty)-L\delta}\phi^{-1}(u)\diff{\Lebesgue}(u)} \\
                            & =-\frac{1}{L}\integratedgeneralizedinverse{\phi}(\phi(+\infty))-C_2.
                \end{align*}
                Therefore, $\integratedgeneralizedinverse{\phi}(\phi(+\infty))>-\infty$ and $b_\phi=+\infty>\phi(+\infty)\geq\Functionimage{\mathfrak{g}}{I_G}$.
        \end{itemize}
        Thus, in all cases, $I_{\integratedgeneralizedinverse{\phi}}\supset\Functionimage{\mathfrak{g}}{I_G}=\Functionimage{\mathfrak{g}\circ H}{I_H}$ and, in particular, for all compact intervals $[a,b]\subset I_H$, we have $ I_{\integratedgeneralizedinverse{\phi}}\supset\closure{\Functionimage{\mathfrak{g}\circ H}{[a,b]}}$.
        Finally, we have that
        \begin{equation*}
            \mathfrak{h}
            \defeq\monotonetransformation{\phi}\circ h
            =\monotonetransformation{\phi}\circ \frac{1}{\Id} \circ g \circ H
            =\integratedgeneralizedinverse{\phi}\circ\phi\circ g \circ H
            =\integratedgeneralizedinverse{\phi}\circ \mathfrak{g} \circ H
        \end{equation*}
        is therefore real-valued on $I_H$.
    \end{proof}
    \begin{delayedproof}{lem: derivative property of generalized inverse and monotone transformation}
        Let $\phi$ be a generator and $G,H,g,h$ as above. In the following, we will only use the restrictions of $G,g$ to  $I_G$ and of $H,h$ to $I_H$.
        Denote $\mathfrak{g}\defeq\phi\circ g$ and
        \begin{equation}
            \label{eqn: def associated absolutely continuous function phi}
            \mathfrak{h}
            \defeq\monotonetransformation{\phi}\circ h
            =\integratedgeneralizedinverse{\phi}\circ \mathfrak{g} \circ H
        \end{equation}
        by \autoref{lem: Welldefinedness of extension of integrated transform for ragular transformed density}.
        Therefore, $\mathfrak{h}$ is a composition of three locally absolutely continuous functions.
        The function $\integratedgeneralizedinverse{\phi}:I_{\integratedgeneralizedinverse{\phi}}\to\R$ is locally absolutely continuous by \autoref{lem: Welldefinedness of extension of integrated transform for ragular transformed density}.
        The local absolute continuity of $\mathfrak{g}$ and $H$ holds by assumption.
        We now prove the local absolute continuity of $\mathfrak{h}$ using \autoref{theo: composition of absolutely continuous functions}.
        Take an arbitrary $[a,b]\subset I_H$. It holds that $\Functionimage{H}{[a,b]}\eqdef[c,d]$ is again a compact subset of $I_G$. Thus, $\mathfrak{g}\circ H:(a,b)\to\R$ is a composition of absolutely continuous functions, and since $H$ is strictly increasing, $\mathfrak{g}\circ H$ is itself absolutely continuous by \autoref{theo: composition of absolute continuous and absolute continuous non-decreasing function}.
        By the resulting continuity, we have that $\Functionimage{\mathfrak{g}\circ H}{[a,b]}\subset I_{\integratedgeneralizedinverse{\phi}}$ is a compact interval. Therefore, there exists a second compact subinterval of $I_{\integratedgeneralizedinverse{\phi}}$ with interior $I_{a,b}$ such that $\Functionimage{\mathfrak{g}\circ H}{(a,b)}\subset I_{a,b}$ and $\integratedgeneralizedinverse{\phi}$ is absolutely continuous when restricted to $I_{a,b}$ by its local absolute continuity.
        We have thus shown the first condition in \autoref{theo: composition of absolutely continuous functions}.
        For the second condition, we start by calculating $\RadonNikodym{\mathfrak{g}\circ H}{\Lebesgue}$ on $(a,b)$.
        We denote a Radon-Nikodym derivative of $\mathfrak{g}$ by $\mathfrak{g}'$. By \cite{zobelmarkusQuantileCharacterizationUnivariate2026arxiv}, $h\defeq\frac{1}{g\circ H}$ is a Radon-Nikodym derivative of $H$, since $G$ is also locally absolutely continuous and $\frac{h}{h}=1$ $\restr{\Lebesgue}{I_H}$-{a.e}.
        Using the second part of \autoref{theo: composition of absolutely continuous functions}, we calculate
        \begin{equation*}
            \RadonNikodym{\mathfrak{g}\circ H}{\Lebesgue}=\RadonNikodym{\mathfrak{g}}{\Lebesgue}\circ H\cdot\RadonNikodym{H}{\Lebesgue}= (\mathfrak{g}'\circ H)\cdot h\quad \text{$\restr{\Lebesgue}{(a,b)}$-{a.e.\@}\nocorr}
        \end{equation*}
        On $I_{\integratedgeneralizedinverse{\phi}}$, we have the Radon-Nikodym derivative
        \begin{equation*}
        \RadonNikodym{\integratedgeneralizedinverse{\phi}}{\Lebesgue}(x)=\begin{cases}
            -\phi^{-1}(x) & x\in\Functionimage{\phi}{\R_>0}\\
            0 & x\notin\Functionimage{\phi}{\R_>0}
        \end{cases},\quad x\in I_{\integratedgeneralizedinverse{\phi}}.
        \end{equation*}
        Since $g$ is positive and finite $\restr{\Lebesgue}{I_G}$-{a.e.\@}, we have that $\restr{\Lebesgue}{a,b}(\phi\circ g\circ H\notin\Functionimage{\phi}{\R_>0})=0$, as $G$ is absolutely continuous on $[c,d]=\Functionimage{H}{[a,b]}$. Therefore,
        \begin{equation*}
        \RadonNikodym{\integratedgeneralizedinverse{\phi}}{\Lebesgue}\circ\mathfrak{g}\circ H=-\phi^{-1}\circ\mathfrak{g}\circ H\quad \text{$\restr{\Lebesgue}{a,b}$-{a.e.\@}\nocorr}
        \end{equation*}
        To finish the proof of the absolute continuity of $\mathfrak{h}$ on $(a,b)$, it remains to check the integrability over $(a,b)$ of
        \begin{align*}
            \mathfrak{h}'\defeq & \left(\RadonNikodym{\integratedgeneralizedinverse{\phi}}{\Lebesgue}\circ(\mathfrak{g}\circ H)\right)\cdot\left(\RadonNikodym{\mathfrak{g}\circ H}{\Lebesgue}\right) = -(\phi^{-1}\circ\phi\circ g\circ H)\cdot(\mathfrak{g}'\circ H)\cdot h \\
            =                   & -\frac{1}{h}\cdot \mathfrak{g}'\circ H\cdot h
            =-(\mathfrak{g}'\circ H)\quad \text{$\restr{\Lebesgue}{(a,b)}$-{a.e.\@}\nocorr}
        \end{align*}
        For $[c,d]$, there exists a Lipschitz constant $L_{[c,d]}$ for $\mathfrak{g}$. Therefore, $\abs{\mathfrak{g}'}\leq L_{[c,d]}$ $\restr{\Lebesgue}{(c,d)}$-{a.e}. Again, by the bi-absolute continuity of $H$ on $(a,b)$, we also have $\abs{\mathfrak{g}'\circ H}\leq L_{[c,d]}$ $\restr{\Lebesgue}{(a,b)}$-{a.e.\@}, and $\mathfrak{h}'$ is integrable over $(a,b)$. Thus, we obtain from \autoref{theo: composition of absolutely continuous functions} that $\mathfrak{h}$ is absolutely continuous on $(a,b)$ and $-(\mathfrak{g}'\circ H)$ is a Radon-Nikodym derivative of $\mathfrak{h}$ there. Therefore, $\mathfrak{h}$ is Lipschitz on $[a,b]$, because
        \begin{equation*}
            \abs{\mathfrak{h}(x)-\mathfrak{h}(y)}=\abs*{\int_y^x\mathfrak{h}'\diff\Lebesgue}\leq L_{[c,d]}\abs{x-y},\quad x<y\in[a,b].
        \end{equation*}
        Since $(a,b)\subset I_H$ was arbitrary, this concludes the proof of the local Lipschitz continuity of $\mathfrak{h}$ on $I_H$ and the equation $\RadonNikodym{\mathfrak{h}}{\Lebesgue}=-(\mathfrak{g}'\circ H)$ $\restr{\Lebesgue}{I_H}$-{a.e}.
    \end{delayedproof}
    \begin{delayedproof}{thm: equivilance between M concave and quantile M convex quasi arithmetic}
        In both cases, $F_\mu$ and $Q_\mu$ are locally absolutely continuous. The former is absolutely continuous by assumption and the latter by either \autoref{Prop: Characterization of unimodal distributions} or \autoref{def: quantile phi-convex distribution} \ref{enum: quantile phi unimodal}. Therefore, $\interior{\support{\mu}}=\Functionimage{Q_\mu}{(0,1)}\neq\emptyset$ is convex, and we set $I_{F_\mu}=\interior{\support{\mu}}$ and $I_{Q_{\mu}}=(0,1)$.
        \begin{itemize}
            \item[\enquote{$\Rightarrow$}:]
                  By assumption, some density $f_{\mu}:\interior{\support{\mu}}\to\Rbar_{>0}$ of $\mu$ is $M_{\phi}$-concave, so the composition $\phi\circ f_\mu:\interior{\support{\mu}}\to\R$ is concave and therefore locally Lipschitz \cite{robertsAnotherProofThat1974}, thus locally absolutely continuous, and there exists a Radon-Nikodym derivative $\RadonNikodym{\phi\circ f_\mu}{\Lebesgue}$ that is non-increasing.
                  By \autoref{lem: derivative property of generalized inverse and monotone transformation} applied to $G=F_\mu$, $H=Q_\mu$, it holds for $q_\mu\defeq\frac{1}{f_\mu\circ Q_\mu}$ that $\monotonetransformation{\phi}\circ q_\mu$ is locally Lipschitz and $\RadonNikodym{\monotonetransformation{\phi}\circ q_\mu}{\Lebesgue}=-\RadonNikodym{\phi\circ f_\mu}{\Lebesgue}\circ Q_\mu$ is non-decreasing. In particular, $\monotonetransformation{\phi}\circ q_\mu$ is real-valued by \autoref{lem: Welldefinedness of extension of integrated transform for ragular transformed density}.
                  Therefore, $\monotonetransformation{\phi}\circ q_\mu:(0,1)\to\R$ is convex, since it is locally absolutely continuous with a non-decreasing Radon-Nikodym derivative \cite[12 A]{robertsConvexFunctions1973}. By \autoref{cor: Aczél correspondence}, $q_\mu$ is $M_{\monotonetransformation{\phi}}$-convex.
            \item[\enquote{$\Leftarrow$}:]
                  Let $\mu$ be quantile $\monotonetransformation{\phi}$-unimodal with $M_{\monotonetransformation{\phi}}$-convex quantile density $q_\mu$ and density $f_\mu=\frac{1}{q_\mu\circ F_\mu}$.
                  Then we conclude from \autoref{lem: derivative property of generalized inverse and monotone transformation} applied to $G=Q_\mu$, $ H=F_\mu$ that $\RadonNikodym{\phi\circ {f_{\mu}}}{\Lebesgue}=-\RadonNikodym{\monotonetransformation{\phi}\circ q_\mu}{\Lebesgue}\circ F_\mu$ is non-increasing on $\interior{\support(\mu)}$.
                  This again gives the $M_{\phi}$-concavity of $f_\mu:\interior{\support{\mu}}\to\Rbar_{>0}$.  \qedhere
        \end{itemize}
    \end{delayedproof}

\subsection{Proof of \autoref{lem: existence of quantile density selections}}
\label{sec: appendix density selection}
    
Let $I$ be an open non-void real interval.
    For two arbitrary measures $\mu,\nu\in\locallyfinitemeasureon(I)$, there exists a set $A\subset I$ such that $\mu_a(B)\defeq\mu(B\cap A^c)$ and $\mu_s(B)\defeq\mu(B\cap A)$ decompose $\mu=\mu_a+\mu_s$ and satisfy $\mu_a\ll\nu$ and $\mu_s\perp\nu$ with $\mu_a,\mu_s\in\locallyfinitemeasureon(I)$. \cite[Corollary 1.29]{kallenbergRandomMeasuresTheory2017} shows that there exists a measurable function $\mathfrak{d}_{I}:(\locallyfinitemeasureon(I)\times\locallyfinitemeasureon(I)\times I,\sigma(\iota_{\B(I)})\otimes\sigma(\iota_{\B(I)})\otimes\B(I))\to(\Rbar_{\geq0},\B(\Rbar_{\geq0}))$ such that
    \begin{equation*}
        \mu=\mathfrak{d}_I(\mu,\nu,\cdot)\cdot\nu+\mu_{s}=\mu_a+\mu_s.
    \end{equation*}
    Thus, $\mathfrak{d}_I(\mu,\nu,\cdot)=\RadonNikodym{\mu_a}{\nu}$ $\nu$-{a.e.\@} for all $\mu,\nu\in\locallyfinitemeasureon(I)$. By construction using \cite[Theorem 1.28]{kallenbergRandomMeasuresTheory2017}, we also have that the choice $A=A_{\mu,\nu}\defeq\set{\mathfrak{d}_I(\mu,\nu,\cdot)=\infty}$ gives the Lebesgue decomposition $\mu_a,\mu_s$ as above. Here, we denote the product $\sigma$-algebra by $\mathcal{A}\otimes\BC$, {i.e.}, $\mathcal{A}\otimes\BC=\sigma(\set{A\times B\given A\in\mathcal{A},B\in\BC})$.
    \begin{lemma}[{\cite[Corollary 1.29]{kallenbergRandomMeasuresTheory2017}}]
        \label{lem: existence of density selections}
        For each open nonvoid interval $I\subset\R$, there exists a function $\DensSelec_I:\locallyfinitemeasureon(I)\times I\to\Rbar_{\geq0}$ called a \emph{density selection} on $I$ which satisfies
        \begin{enumerate}
            \item $\DensSelec_I$ is $\sigma(\iota_{\B(I)})\otimes\B(I),\B(\Rbar_{\geq0})$-measurable and
            \item $\DensSelec_I(\mu,\cdot)=\RadonNikodym{\mu}{\restr{\Lebesgue}{I}}$ for all absolutely continuous $\mu\in\locallyfinitemeasureon(I)$.
        \end{enumerate}
    \end{lemma}
\begin{proof}
    \label{proof: lem: existence of density selections}
    Take $\mathfrak{d}_{I}$, which exists by \cite[Corollary 1.29]{kallenbergRandomMeasuresTheory2017}, and define $\DensSelec_I(\mu,t)=\mathfrak{d}_{I}(\mu,\restr{\Lebesgue}{I},t)$ for each $(\mu,t)\in\locallyfinitemeasureon(I)\times I$. If $\mu$ is absolutely continuous, then $\mu_a=\mu$, and thus $\DensSelec_I(\mu,\cdot)=\mathfrak{d}_{I}(\mu,\restr{\Lebesgue}{I},\cdot)=\RadonNikodym{\mu}{\restr{\Lebesgue}{I}}$. Further, as $\mathfrak{d}_{I}$ is measurable, $\mathfrak{d}_{I}\circ\iota_{\cdot,\restr{\Lebesgue}{I},\cdot}=\DensSelec_I$ is a measurable function from $\locallyfinitemeasureon(I)\times I$ to $\Rbar_{\geq0}$, {i.e.}, $\DensSelec_I$ is $(\sigma(\iota_{\B(I)})\otimes\B(I)),\B(\Rbar_{\geq0})$-measurable.
\end{proof}

\begin{delayedproof}{lem: existence of quantile density selections}
    Based on \autoref{lem: existence of density selections}, there exists a density selection $\DensSelec_{(0,1)}$ on $(0,1)$. We define $\QuantDensSelec\defeq \DensSelec_{(0,1)}\circ(\InverseMeasure\times\Id_{(0,1)})$.
    This mapping is measurable as a composition of measurable functions. If $\mu$ is quantile absolutely continuous, then $\QuantDensSelec(\mu,\cdot)=\RadonNikodym{\mu^{-1}}{\uniformmeasure}=q_\mu$ $\uniformmeasure$-{a.e.\@}, and hence  $\QuantDensSelec(\mu,\cdot)$ is a quantile density.
\end{delayedproof}

\subsection{Proof of \autoref{lem: aczel correspondence}}
\label{sec: Proof of aczel correspondence}

\begin{lemma}
    \label{lem: non-constant M_phi shaped functions}
    Let $\phi$ be a generator on $\R_{>0}$ and $I\subset\R$ an open interval. If $g:I\to\Rbar_{\geq0}$ is an $M_\phi$-convex (-concave) function and $\phi\circ g$ attains the value $-\infty$ ($+\infty$), then $g$ is constant.
\end{lemma}
\begin{proof}
    Let $g:I\to\Rbar_{\geq0}$ be $M_\phi$-convex and let $\phi\circ g$ take the value $-\infty$ at $x\in I$. Thus, necessarily $g(x)=\phi^{-1}(-\infty)=0$. We have that for all $y\in I$, there are $z\in I$, $\convexindex\in(0,1)$ such that $A_{\convexindex}(x,z)=y$. As $g$ is $M_\phi$-well-defined, we have that $\phi\circ g(z)\neq+\infty$ and thus $\convexindex\phi(g(x))+(1-\convexindex)\phi(g(z))=-\infty$. Thus, we obtain for each $x\neq y\in I$
    \begin{equation*}
        0\leq g(y)\leq M_{\phi,\convexindex}(g(x),g(z))=\phi^{-1}(-\infty)=0.
    \end{equation*}
    Therefore, $g\equiv0$, so $g$ is constant.
    Analogously, if $g$ is $M_\phi$-concave and $\phi\circ g$ takes the value $+\infty$ at $x\in I$, then again $g$ is constant, this time with value $+\infty$, since for all $y\in I$ there exists $z\in I$ such that
    \begin{equation*}
        +\infty\geq g(y)\geq M_{\phi,\convexindex}(g(x),g(z))=\phi^{-1}(+\infty +(1-\convexindex)\phi(g(z)))=\phi^{-1}(+\infty)=+\infty
    \end{equation*}
    as this time $\phi(g(z))\neq-\infty$ for all $z\in I$.
\end{proof}
\begin{delayedproof}{lem: aczel correspondence}
    \label{proof: lem: aczel correspondence}
    Let $g$ be $M_\phi$-convex (-concave) and non-constant. By \autoref{lem: non-constant M_phi shaped functions}, we have that $\phi\circ{g}:I\to\Rbar$ does not attain the value $-\infty$ ($+\infty$). Conversely, let $\phi\circ g: I\to\Rbar$ be convex (concave). Then it does not attain $-\infty$ ($+\infty$) and $g$ is $M_\phi$-well-defined.
    Thus, we only need to match the inequalities. For $x,y\in I$, $\convexindex\in(0,1)$, the following holds:
    \begin{align*}
        (\phi\circ g)(A_{\convexindex}(x,y))\overset{(\geq)}     & {\leq}\convexindex(\phi\circ g)(x)+(1-\convexindex)(\phi\circ g)(y)        \\
        \Leftrightarrow g(A_{\convexindex}(x,y))\overset{(\geq)} & {\leq}\phi^{-1}(A_{\convexindex}((\phi\circ g)(x),(\phi\circ g)(y)))=M_{\phi,\convexindex}(g(x),g(y)).\qedhere
    \end{align*}
\end{delayedproof}
We record the fact that the converse direction does not require the non-constant assumption.
\begin{corollary}
    \label{cor: Aczél correspondence}
    Let $I$ be an open nonvoid interval and $\phi$ a generator. If $\phi\circ g:I\to\Rbar$ is convex (concave) for $g:I\to\Rbar_{\geq0}$, then $g$ is $M_{\phi}$-convex (-concave).
\end{corollary}

\subsection{Proof of \autoref{lem: Bijectivity of Integrated inverse Transform}}
\label{proof: lem: Bijectivity of Integrated inverse Transform}

    \begin{delayedproof}{lem: Bijectivity of Integrated inverse Transform}
        Let $\phi:\R_{>0}\to\R$ be a generator and fix $c=\frac12$. That is, $\phi$ is continuous and strictly increasing. We have that
        \begin{equation*}
            \monotonetransformation{\phi}(c)
            =-\int_{\phi(\frac1c)}^{\phi(\frac1c)}\phi^{-1}\diff \Lebesgue+\phi(c)
            =\phi(c),
        \end{equation*}
        and thus $\monotonetransformation{(\monotonetransformation{\phi})}(c)=\phi(c)$. Therefore, it suffices to show that $\phi(b)-\phi(a)=\monotonetransformation{(\monotonetransformation{\phi})}(b)-\monotonetransformation{(\monotonetransformation{\phi})}(a)$ for all $a,b\in \R_{>0}$ to obtain the equality.
        This is equivalent to showing that $\Stieltjestransform{\phi}=\Stieltjestransform{\monotonetransformation{(\monotonetransformation{\phi})}}$.
        We thus begin by investigating $\Stieltjestransform{\monotonetransformation{\phi}}$.
        Take an arbitrary $[a,b]\subset\R_{>0}$. Then
        \begin{align*}
            \Stieltjestransform(\monotonetransformation{\phi})((a,b])
             & =
            \monotonetransformation{\phi}(b)-\monotonetransformation{\phi}(a)
            = -\int_{\phi(\frac1c)}^{\phi(\frac1b)}\phi^{-1}\dd\Lebesgue +\int_{\phi(\frac1c)}^{\phi(\frac1a)}\phi^{-1}\dd\Lebesgue  \\
             & =\int_{\phi(\frac1b)}^{\phi(\frac1a)}\phi^{-1}\dd\Lebesgue
            =\int_{[\phi(\frac1b),\phi(\frac1a)]}\phi^{-1}\dd\Lebesgue=\int_{\Functionimage{\phi}{[\frac1b,\frac1a)}}\phi^{-1}\dd\Lebesgue \\
             & =\int_{[\frac1b,\frac1a)}\Id\dd\pushforward{\phi^{-1}}{\Lebesgue}
            =\int_{[\frac1b,\frac1a)}\Id\dd\Stieltjestransform{\phi}
            =\int_{[\frac1b,\frac1a)}\frac{1}{\frac1y}\diff\Stieltjestransform{\phi}(y)                              \\
             & =\int_{(a,b]}\frac1y\diff\pushforward{\frac1\Id}{\Stieltjestransform{\phi}}(y)
            =\frac{1}{\Id}\cdot\pushforward{\frac1\Id}{\Stieltjestransform(\phi)}((a,b]).
        \end{align*}
        Therefore, $\Stieltjestransform(\monotonetransformation{\phi})=\frac{1}{\Id}\cdot\pushforward{\frac1\Id}{\Stieltjestransform(\phi)}$.
        
        Taking $[a,b]\subset\R_{>0}$ again arbitrary, we obtain
        \begin{align*}
            \monotonetransformation{(\monotonetransformation{\phi})}(b) -\monotonetransformation{(\monotonetransformation{\phi})}(a)
             & =\Stieltjestransform(\monotonetransformation{(\monotonetransformation{\phi})})((a,b])
            = \int_{[\frac1b,\frac1a)}\Id\diff\Stieltjestransform(\monotonetransformation{\phi})                \\
             & = \int_{[\frac1b,\frac1a)}\frac{\Id}{\Id}\diff\pushforward{\frac1\Id}{\Stieltjestransform(\phi)}
            =\int_{\preimage{\frac{1}{\Id}}{(a,b]}}1\diff\pushforward{\frac1\Id}{\Stieltjestransform(\phi)}     \\
             & =\int_{(a,b]}1\diff\Stieltjestransform(\phi)
            =\phi(b)-\phi(a).
        \end{align*}
        Therefore, $\monotonetransformation{(\monotonetransformation{\phi})}(x)- \phi(x)=\monotonetransformation{(\monotonetransformation{\phi})}(c)- \phi(c)=0$ for all $x\in\R_{>0}$.
    \end{delayedproof}

\subsection{Proof of \autoref{thm: characterization of not absolutely continuous M concave distribution}}
\label{sec: proof of chararacterization of phi unimodality}

    Recall the characterization of non-absolutely continuous and non-degenerate unimodal distributions.
    \begin{lemma}[{\cite[Lemma 4.4]{zobelmarkusQuantileCharacterizationUnivariate2026arxiv}}]
        \label{lem: characterization of unimodal distributions via their Lebesgue decompsition}
        A distribution $\mu$ is unimodal if and only if $\mu=\convexindex\mu_{a}+(1-\convexindex)\delta_{m}$ for some $\convexindex\in[0,1]$, $m\in\R$, where $\mu_a\ll\Lebesgue$ is absolutely continuous and unimodal, and $\mu_a$ and $\delta_{m}$ share their mode, {i.e.}, there exists a quasi-concave density $f_a$ of $\mu_a$ such that $m$ maximizes $f_a$.
    \end{lemma}

    \begin{lemma}
        \label{lem: Calculation relation M_phi convexity of quantile density and quantile density of absollutely continuous part}
        Let $\mu=\convexindex\mu_{a}+(1-\convexindex)\delta_{m}$ be unimodal with $\convexindex\in(0,1)$, $\mu_a\ll\Lebesgue$, $m\in\R$. Then
        $\mu$ is quantile $\phi$-unimodal if and only if $\mu_a$ is quantile $\tilde{\phi}$-unimodal with $\tilde{\phi}^{-1}=\convexindex\phi^{-1}$ and $\lim_{t\to F_{\mu_a}(m),t\in(0,1)}q_{\mu_a}(t)=0$ for an $M_{\tilde{\phi}}$-convex quantile density $q_{\mu_a}$ of $\mu_a$.
    \end{lemma}
    \begin{delayedproof}{thm: characterization of not absolutely continuous M concave distribution}
        By \autoref{lem: characterization of unimodal distributions via their Lebesgue decompsition}, we know that in both cases $\mu$ is unimodal with decomposition $\mu=\convexindex\mu_a +(1-\convexindex)\delta_{m}$ for $\convexindex\in[0,1]$.
        For $\convexindex=1$, the statement reduces to \autoref{thm: equivilance between M concave and quantile M convex quasi arithmetic}. For $\convexindex=0$, we have that $\mu=\delta_m$, so $\mu$ is $\phi$-unimodal. Thus, we can take $\mu_a\defeq\restr{\Lebesgue}{[m-\frac12,m+\frac12]}$ and $f_a\defeq\1_{[m-\frac12,m+\frac12]}$. Then, $\convexindex f_a$ is $M_\phi$-concave and $m$ is a mode. Assume now that $\convexindex\in(0,1)$.
        \begin{itemize}
            \item[\enquote{$\Rightarrow$}:]
                If $\mu$ is $\phi$-unimodal, then there exists a quantile density $q_\mu$ which is $M_{\monotonetransformation{\phi}}$-convex.
                We apply \autoref{lem: Calculation relation M_phi convexity of quantile density and quantile density of absollutely continuous part} to obtain that there exists a quantile density $q_{\mu_a}$ of $\mu_a$ which is $M_{\tilde{\monotonetransformation{\phi}}}$-convex and $\lim_{t\to F_{\mu_a}( m),t\in(0,1)}q_{\mu_a}(t)=0$.
                As $\mu_a$ is absolutely continuous, it follows from \autoref{thm: equivilance between M concave and quantile M convex quasi arithmetic} that $\mu_a$ is ${\monotonetransformation{(\tilde{\monotonetransformation{\phi}})}}$-unimodal with an $M_{\monotonetransformation{(\tilde{\monotonetransformation{\phi}})}}$-concave density $f_{a}\defeq\frac{1}{q_{\mu_a}\circ F_{\mu_a}}$ on $\interior{\support{\mu_a}}$. Because $q_{\mu_a}\not\equiv0$, the limit condition $\lim_{t\to F_{\mu_a}( m),t\in(0,1)}q_{\mu_a}(t)=0$ implies that $q_{\mu_a}$ and thus $f_{a}$ are non-constant.
                Since $\tilde{\monotonetransformation{\phi}}(x)=\monotonetransformation{\phi}(\frac{1}{\convexindex}x)$, we have that
                \begin{align*}
                    \monotonetransformation{(\tilde{\monotonetransformation{\phi}})}(x)
                    & =-\int_{\tilde{\monotonetransformation{\phi}}(\frac{1}{c})}^{\tilde{\monotonetransformation{\phi}}(\frac1x)}\tilde{\monotonetransformation{\phi}}^{-1}\diff\Lebesgue +\tilde{\monotonetransformation{\phi}}(c)
                    =-\int_{\monotonetransformation{\phi}(\frac{1}{\convexindex c})}^{\monotonetransformation{\phi}(\frac1{\convexindex x})}\convexindex\monotonetransformation{\phi}^{-1}\diff\Lebesgue +\tilde{\monotonetransformation{\phi}}(c)                                                                                                                                                             \\
                    & =-\convexindex\int_{\monotonetransformation{\phi}(\frac{1}{c})}^{\monotonetransformation{\phi}(\frac1{\convexindex x})}\monotonetransformation{\phi}^{-1}\diff\Lebesgue +\int_{\monotonetransformation{\phi}(\frac{1}{c})}^{\monotonetransformation{\phi}(\frac1{\convexindex c})}\convexindex\monotonetransformation{\phi}^{-1}\diff\Lebesgue+\tilde{\monotonetransformation{\phi}}(c) \\
                    & =\convexindex\monotonetransformation{(\monotonetransformation{\phi})}(\convexindex \cdot x) + C
                    =\convexindex \phi(\convexindex\cdot x) + C \eqdef \convexindex \breve{\phi}(x)+ C.
                \end{align*}
                Therefore, we have that $f_a$ is $M_{\monotonetransformation{(\tilde{\monotonetransformation{\phi}})}}=M_{\convexindex\breve{\phi}+C}=M_{\breve{\phi}}$-concave and non-constant. Thus, by \autoref{lem: aczel correspondence},
                \begin{equation*}
                    \breve{\phi}\circ f_a=\phi\circ(\convexindex \cdot f_a)
                \end{equation*}
                is concave. Hence, $\convexindex f_a$ is $M_\phi$-concave.
                Further, by the definition of $f_a$ on $\interior{\support{\mu_a}}$, we have that
                \begin{equation*}
                    \lim_{\substack{x\to m\\x\in\interior{\support{\mu_a}}}} f_a(x)
                    =\lim_{\substack{x\to m\\x\in\interior{\support{\mu_a}}}} \frac{1}{q_{\mu_a}\circ {F_{\mu_a}}(x)}
                    =\lim_{\substack{t\to F_{\mu_a}( m)\\t\in(0,1)}}\frac{1}{q_{\mu_a}(t)}=+\infty.
                \end{equation*}
            \item[\enquote{$\Leftarrow$}:]
                Since $\convexindex f_a$ is $M_\phi$-concave, we obtain that $f_a$ is again $M_{\breve{\phi}}$-concave. Thus, by \autoref{thm: equivilance between M concave and quantile M convex quasi arithmetic}, $q_{\mu_a}\defeq\frac{1}{f_a\circ Q_{\mu_a}}$ is $M_{\monotonetransformation{\breve{\phi}}}=M_{\tilde{\monotonetransformation{\phi}}}$-convex. Furthermore, from $\lim_{\substack{x\to m,x\in\interior{\support{\mu_a}}}} f_a(x)=+\infty$, we calculate that
                \begin{equation*}
                    \lim_{\substack{t\to F_{\mu_a}( m)\\t\in(0,1)}}q_{\mu_a}(t)= \lim_{\substack{x\to m\\x\in\interior{\support{\mu_a}}}}\frac{1}{f_{a}(x)}=0.
                \end{equation*}
                Applying \autoref{lem: Calculation relation M_phi convexity of quantile density and quantile density of absollutely continuous part}, we obtain that $\mu$ is quantile $\monotonetransformation{\phi}$-unimodal and thus $\mu$ is $\phi$-unimodal by \autoref{def: M concave general}.\qedhere
        \end{itemize}
    \end{delayedproof}
    \begin{delayedproof}{lem: Calculation relation M_phi convexity of quantile density and quantile density of absollutely continuous part}
        By assumption, we know that both $\mu$ and $\mu_a$ are unimodal and thus both $Q_\mu$ and $Q_{\mu_{a}}$ are locally absolutely continuous on $(0,1)$ with quantile densities $q_\mu,q_{\mu_a}$, respectively. We formally extend $Q_\mu,Q_{\mu_a},q_\mu,q_{\mu_a}$ to $\set{0,1}$ via one-sided limits, {e.g.\@}, $Q_\mu(0)\defeq\lim_{t\searrow0}Q_\mu(t)\in\Rbar$ and $q_{\mu_a}(1)\defeq\essliminf_{t\nearrow1}q_{\mu_a}(t)$.
        Additionally, we set $Q_\mu,Q_{\mu_a}$ to $-\infty$ on $(-\infty,0)$ and $+\infty$ on $(1,+\infty)$.
        We first relate $Q_\mu$ and $Q_{\mu_a}$ on $(0,1)$. The relation between the {c.d.f.s} is $F_\mu= \convexindex F_{\mu_a} + (1-\convexindex) \ind_{[ m,\infty)}$. We calculate for $t\in(0,1)$
        \begin{align*}
            Q_\mu(t)
            & =\sup(\preimage{F_\mu}{[0,t]})
            =\sup\set{x\in\R\given F_\mu(x)\leq t}\\
            &=\sup\set{x\in\R\given \convexindex F_{\mu_a}(x) + (1-\convexindex) \ind_{[ m,\infty)}(x)\leq t}                                                                                                                           \\
            & =\sup\set*{x\in\R\given  F_{\mu_a}(x) \leq \frac{t}{\convexindex} - \frac{1-\convexindex}{\convexindex}\ind_{[ m,\infty)}(x)}                                                                                            \\
            & =\max\set*{\sup\set*{x\in\R\given  F_{\mu_a}(x) \leq \frac{t}{\convexindex} - \frac{1-\convexindex}{\convexindex}\ind_{[ m,\infty)}(x)},\sup{\set*{ x\in\R\given F_{\mu_a}(x)\leq \frac{t-1+\convexindex}{\convexindex}}}} \\
            & =\max\set*{\sup\set*{x<m\given F_{\mu_a}(x) \leq \frac{t}{\convexindex}},\sup{\set*{ x\in\R\given F_{\mu_a}(x)\leq \frac{t-1+\convexindex}{\convexindex}}}}                                                                \\
            & =\max{\set*{
            \min\parenthesis*{Q_{\mu_a}\parenthesis*{\frac{t}{\convexindex}}, m},Q_{\mu_a}\parenthesis*{\frac{t-1+\convexindex}{\convexindex}}}}                                                                                       \\
            & =\begin{cases}
                    Q_{\mu_a}(\frac{t}{\convexindex})                             & t< F_\mu( m-)               \\
                    m                                                             & F_\mu( m-)\leq t< F_\mu( m) \\
                    Q_{\mu_a}\parenthesis*{\frac{t-1+\convexindex}{\convexindex}} & F_\mu( m)\leq t.
                \end{cases}
        \end{align*}
        Analogously,
        \begin{equation*}
            Q_{\mu_a}(t)
            =\begin{cases}
                Q_\mu(\convexindex t)                  & t <F_{\mu_a}( m)     \\
                Q_\mu(\convexindex t+(1-\convexindex)) & F_{\mu_a}( m)\leq t.
            \end{cases}
        \end{equation*}
        We can define two piecewise affine deformations on $(0,1)$ as follows. For $t\in(0,1)$, let
        \begin{align*}
            W(t)
            \defeq\begin{cases}
                    \frac{t}{\convexindex}                & t< F_\mu( m-)               \\
                    F_{\mu_a}( m)                         & F_\mu( m-)\leq t< F_\mu( m) \\
                    \frac{t-1+\convexindex}{\convexindex} & F_\mu( m)\leq t
                \end{cases}\quad
            W_a(t)
            \defeq\begin{cases}
                    \convexindex t                   & t<F_{\mu_a}( m)     \\
                    \convexindex t +(1-\convexindex) & F_{\mu_a}( m)\leq t
                \end{cases}
        \end{align*}
        such that $Q_\mu=Q_{\mu_a}\circ W$ and $Q_{\mu_a}= Q_\mu\circ W_a$, using the one-sided limit extensions of both to $[0,1]$.
        Thus, we obtain the following equations relating the quantile densities of $\mu$ and $\mu_a$:
        \begin{align}
            \label{eq:def: q_mu from q_mu_a}
            q_\mu     & =
            \1_{(0,F_\mu(m-))\cup(F_\mu(m),1)}
            \frac{1}{\convexindex}q_{\mu_a}\circ W
                    \quad \text{$\uniformmeasure$-{a.e.\@},}
            \\
            \label{eq: def: q_mu_a from q_mu}
            q_{\mu_a} & =\convexindex q_\mu\circ W_a
                    \quad \text{$\uniformmeasure$-{a.e.\@}\nocorr}
        \end{align}
        We again use the extensions of $q_\mu,q_{\mu_a}$ to $[0,1]$ and interpret $0\cdot(+\infty)=0$ in \eqref{eq:def: q_mu from q_mu_a}.
        If $q_{\mu_a}(F_{\mu_a}(m))=0$, then \eqref{eq:def: q_mu from q_mu_a} simplifies to $ q_\mu =\frac{1}{\convexindex}q_{\mu_a}\circ W$ without the need for that convention.

        For the following proof, denote $\beta\defeq F_{\mu_a}(m)$.
        \begin{itemize}
            \item[\enquote{$\Rightarrow$}:]
                Let $q_\mu$ be the $M_{\phi}$-convex quantile density of $\mu$. As $\convexindex>0$, we have that $q_\mu$ is not constant, and thus, by \autoref{lem: aczel correspondence}, we have that $\phi\circ q_\mu$ is convex.
                Take $q_{\mu_a}$ as defined in \eqref{eq: def: q_mu_a from q_mu}. Then
                \begin{equation*}
                    \tilde{\phi}\circ q_{\mu_a}=\tilde{\phi}(\convexindex\cdot q_\mu\circ W_a)=\phi\circ q_\mu\circ W_a.
                \end{equation*}
                As $\convexindex<1$, we have that $F_\mu(m-)\neq F_\mu(m)$. Therefore, $Q_\mu$ is constant on $[F_\mu(m-),F_\mu(m)]$ and $\restr{q_\mu}{[F_\mu(m-),F_\mu(m)]}\equiv 0$.
                Thus, $\tilde{\phi}\circ q_{\mu_a}$ is convex and non-increasing on $(0,\beta)$ and convex and non-decreasing on $(\beta,1)$, since $W_a$ is affine on these open sets. For $\beta\in\set{0,1}$, we have already shown the convexity of $\tilde{\phi}\circ q_{\mu_a}$.
                By the continuity of $q_\mu=\phi^{-1}\circ(\phi\circ q_\mu)$, we have that
                \begin{equation*}
                    \lim_{\substack{t\nearrow F_\mu(m-)\\t\in(0,1)}}q_\mu(t)=\lim_{\substack{t\nearrow \beta\\t\in(0,1)}}q_\mu(W_a(t))=\lim_{\substack{t\searrow F_\mu(m)\\t\in(0,1)}}q_\mu(t)=\lim_{\substack{t\searrow \beta\\t\in(0,1)}}q_\mu(W_a(t))=0.
                \end{equation*}
                When $\beta\in\set{0,1}$, we can only establish these limits one-sidedly.
                Nevertheless, this suffices to obtain
                \begin{equation*}
                    \lim_{\substack{t\to \beta\\t\in(0,1)}}q_{\mu_a}(t)=\lim_{\substack{t\to \beta\\t\in(0,1)}}\convexindex q_\mu\circ W_a(t)=0=\convexindex q_{\mu}(F_{\mu}(m))=\convexindex q_{\mu}(W_a(\beta))=q_{\mu_a}(\beta).
                \end{equation*}
                Therefore, $q_{\mu_a}$ is continuous at $\beta$. This implies the continuity of $\tilde{\phi}\circ q_{\mu_a}$ at $\beta$.
                Joining a non-increasing convex function and a non-decreasing convex function continuously yields a convex function. Thus, we have shown both the $M_{\tilde{\phi}}$-convexity of $q_{\mu_a}$ and the limit condition $\lim_{t\to F_{\mu_a}(m),t\in(0,1)}q_{\mu_a}(t)=0$.
            \item[\enquote{$\Leftarrow$}:]
                Assume that $q_{\mu_a}$ is $M_{\tilde{\phi}}$-convex and $\lim_{t\to \beta,t\in(0,1)}q_{\mu_a}(t)=0$. Therefore, $q_{\mu_a}$ is non-constant, since otherwise the limit condition would force it to vanish identically, contradicting $\mu_a\ll\Lebesgue$. By \autoref{lem: aczel correspondence}, $\tilde{\phi}\circ q_{\mu_a}$ is convex.
                Therefore, $q_{\mu_a}$ is continuous and $q_{\mu_a}(\beta)=0$.
                Then, \eqref{eq:def: q_mu from q_mu_a} simplifies to $q_\mu\defeq \frac1\convexindex q_{\mu_a}\circ W$, giving a quantile density of $\mu$, and
                \begin{equation*}
                    \phi\circ q_\mu=\phi\circ \parenthesis *{\frac{1}{\convexindex}\cdot (q_{\mu_a}\circ W)}=\tilde{\phi}\circ q_{\mu_a}\circ W.
                \end{equation*}
                As $\lim_{t\to \beta,t\in(0,1)}q_{\mu_a}(t)=0$, we have that $\phi\circ q_\mu$ is strictly decreasing and convex to the left of $F_\mu( m-)$ and strictly increasing and convex to the right of $F_\mu( m)$, as $W$ is affine on these sets.
                We further know that $\tilde{\phi}\circ q_{\mu_a}(\beta)\neq-\infty$, as $\tilde{\phi}\circ q_{\mu_a}$ is bounded from below as a convex function on $(0,1)$.
                As $\beta$ is the minimizer of $\tilde{\phi}\circ q_{\mu_a}$ and $\tilde{\phi}\circ q_{\mu_a}$ is continuous there, we know that the constant region $[F_{\mu}( m-),F_\mu( m)]$ of $\phi\circ q_\mu$ fits between the other two, and $\phi\circ q_\mu$ is convex.
                More formally, by the continuity of $\tilde{\phi}\circ q_{\mu_a}\circ W$, we have that $\phi\circ q_\mu$ is locally Lipschitz, as $\restr{\phi\circ q_\mu}{[F_{\mu}(m-),F_\mu(m)]}$ is constant and finite and a Radon-Nikodym derivative is non-decreasing and negative, constant $0$, non-decreasing and positive on $(0,F_\mu( m-)),[F_{\mu}( m-),F_\mu( m)],(F_\mu( m),1)$, respectively, and thus non-decreasing on $(0,1)$. Therefore, $q_\mu$ is $M_\phi$-convex.\qedhere
        \end{itemize}
    \end{delayedproof}

\subsection{Proof of \autoref{lem: existence of lower semicontinuous modification of quantile density}}
\label{sec: proof of density selection for quantile densities}
\begin{delayedproof}{lem: existence of lower semicontinuous modification of quantile density}
    By \autoref{lem: existence of quantile density selections}, there exists a quantile density selection $\QuantDensSelec:\distributionon(\R)\times(0,1)\to\Rbar_{\geq0}$.
    Define $\QuantDensSelec^*:\distributionon(\R)\times(0,1)\to\Rbar_{\geq0}$ as
    \begin{equation*}
        \QuantDensSelec^*(\mu,t)=\essliminf_{y\to t}\QuantDensSelec(\mu,y),\quad \mu\in\distributionon(\R),t\in(0,1).
    \end{equation*}
    By \cite{piccioniContinuousVersionsRadon1983}, $\QuantDensSelec^*(\mu)=\QuantDensSelec^*(\mu,\cdot)$ is lower semicontinuous and the greatest lower semicontinuous {a.e.\@} minorant of $\QuantDensSelec(\mu)$, {i.e.,} $\QuantDensSelec^*(\mu)\leq\QuantDensSelec(\mu)$ $\Lebesgue$-{a.e.\@} for all $\mu\in\distributionon(\R)$. Furthermore, it holds that $\QuantDensSelec^*(\mu,t)=\essliminf_{y\to t}q_\mu(y)$ for any quantile density $q_\mu$ of a quantile absolutely continuous $\mu$.
    The following measurability argument adapts an idea from gerw's answer \cite{gerwAnswerMeasurabilityEssential2021}.
    We first extend $\QuantDensSelec$ to all of $\R$ by defining $\overline{\QuantDensSelec}:\distributionon(\R)\times \R\to[0,+\infty]$ as
    \begin{equation*}
        \overline{\QuantDensSelec}(\mu,t)=\begin{cases}
            \QuantDensSelec(\mu,t) & t\in(0,1)\\
            0 & t\in(0,1)^c.
        \end{cases}
    \end{equation*}
    Now, for each $r\in\R_{>0}$, we define the translated function
    $\overline{\QuantDensSelec}_r:\distributionon(\R)\times(0,1)\times(-r,r)\to[0,+\infty]$ as $\overline{\QuantDensSelec}_r(\mu,t,\delta)=\overline{\QuantDensSelec}(\mu,t+\delta)$. This function is still product measurable as the composition of product measurable functions. With this, we can now define the local essential lower bounds of $\QuantDensSelec$ as
    $\QuantDensSelec_r:\distributionon(\R)\times (0,1)\to[0,+\infty]$ with $\QuantDensSelec_r(\mu,t)=\essinf_{\delta\in(-r,r)}\overline {\QuantDensSelec}_r(\mu,t,\delta)$. This function is measurable:
    Define, for each $R\in[0,+\infty]$,
    \begin{equation*}
        h_{R}(\mu,t)
        \defeq\Lebesgue(\set{\delta\in(-r,r)\given \overline{\QuantDensSelec}_r(\mu,t,\delta)<R})
        =\int_{-r}^r
        \1_{\set{\overline{\QuantDensSelec}_r<R}}(\mu,t,\delta)
        \dd\Lebesgue(\delta).
    \end{equation*}
    Then $h_R:\distributionon(\R)\times(0,1)\to[0,+\infty]$ is measurable by Tonelli's theorem, as $\1_{\set{\overline{\QuantDensSelec}_r<R}}$ is product measurable and nonnegative.
    Now we have that
    \begin{align*}
        \set{\QuantDensSelec_r(\mu,t)<R}&=\set{(\mu,t)\in\distributionon(\R)\times(0,1)\given \Lebesgue(\set{\delta\in(-r,r)\given \overline{\QuantDensSelec}_r(\mu,t,\delta)<R})>0}\\
        &=\set{h_R(\mu,t)>0}\in\sigma(\iota_{\B(\R)})\tensor\B((0,1))
    \end{align*}
    for each $R\in[0,+\infty]$, and thus $\QuantDensSelec_r$ is measurable. Now,
    \begin{equation*}
        \sup_{n\in\N}\QuantDensSelec_{\frac{1}{n}}(\mu,t)=\sup_{r>0}\essinf_{\delta\in(-r,r)}\overline{\QuantDensSelec}(\mu,t+\delta)=\essliminf_{y\to t}\QuantDensSelec(\mu,y)=\QuantDensSelec^*(\mu,t).
    \end{equation*}
    As the $(\QuantDensSelec_r)_{r\in\R_{>0}}$ are non-decreasing as $r\to0$, the countable supremum coincides with the supremum. Therefore, $\QuantDensSelec^*$ is measurable.

    Define $f(\mu)=\int_{(0,1)}\abs*{\frac{\QuantDensSelec^*(\mu,t)}{1+\QuantDensSelec^*(\mu,t)}-\frac{\QuantDensSelec(\mu,t)}{1+\QuantDensSelec(\mu,t)}}\diff(t)$. Then this function is measurable by Tonelli's theorem, as $\abs*{\frac{\QuantDensSelec^*}{1+\QuantDensSelec^*}-\frac{\QuantDensSelec}{1+\QuantDensSelec}}$ is product measurable and nonnegative, using the continuous extension $\frac{+\infty}{1+\infty}=1$.
    If $f(\mu)=0$, then $\QuantDensSelec^*(\mu)=\QuantDensSelec(\mu)$ $\Lebesgue$-{a.e}. Thus, $\QuantDensSelec(\mu)$ has a lower semicontinuous version. On the other hand, if $\QuantDensSelec(\mu)$ has a lower semicontinuous version $q_\mu$, then, as $\QuantDensSelec^*(\mu)$ is the {a.e.\@} lower semicontinuous minorant of $\QuantDensSelec(\mu)$, we have
    \begin{equation*}
        \QuantDensSelec(\mu)\geq\QuantDensSelec^*(\mu)\geq q_\mu=\QuantDensSelec(\mu)\quad \text{$\Lebesgue$-{a.e}.}
    \end{equation*}
    Therefore, $f(\mu)=0$. By this argument, the set of all measures $\mu$ such that there exists a lower semicontinuous version of $\QuantDensSelec(\mu,\cdot)$ is $\preimage{f}{\set{0}}$ and therefore measurable.
    Thus,
    \begin{equation*}
        \tilde{\QuantDensSelec}(\mu,t)\defeq\begin{cases}
            \QuantDensSelec^*(\mu,t) & f(\mu)=0\\
            \QuantDensSelec(\mu,t) & f(\mu)\neq 0
        \end{cases},\quad \mu\in\distributionon(\R),t\in(0,1)
    \end{equation*}
    defines a quantile density selection. By definition, $\tilde{\QuantDensSelec}(\mu)=\QuantDensSelec^*(\mu)$ holds if and only if $f(\mu)=0$, and above we have shown that $f(\mu)=0$ is equivalent to $\QuantDensSelec(\mu)$ having a lower semicontinuous version, which is the same as $\tilde{\QuantDensSelec}(\mu)$ having one.

\end{delayedproof}

\subsection{Proof of \autoref{theo:generalized Jensen inequality}}
\label{sec: proof of Jensens inequality}

\begin{delayedproof}{theo:generalized Jensen inequality}
    Let $(b_n)_{n\in\N}\subset\R_{>0}$ with $b_n\nearrow +\infty$ and let $h=(h^{(1)},h^{(2)}):\Omega\to\Rbar_{\geq0}^2$ be measurable and $M_\phi$-well-defined with $\int h\diff\P\in C_\phi$. Define
    \begin{equation*}
        h_n\defeq  \begin{pmatrix}
            \min(h^{(1)},b_n)\\
            \min(h^{(2)},b_n)
        \end{pmatrix}.
    \end{equation*}
    Then $h_n^{(j)}\in\Lp^1(\P)$, $j\in\set{1,2}$, and $(h_n)_{n\in\N}$ is a non-decreasing sequence of $M_\phi$-well-defined functions to $\R_{\geq0}^2$ that converges pointwise to $h$. By this monotonicity, we have that $\int h_n\diff\P\nearrow \int h\diff\P$ as $n\to\infty$, and therefore also $\int h_n\diff\P\in C_\phi$ for all $n\in\N$. Let $\convexindex\in(0,1)$ be arbitrary.
    For a convex quasi-arithmetic mean $M_\phi:\R_{\geq0}^2\to\R_{\geq0}$, we can apply the multivariate Jensen inequality (see, for example, \cite[2.8 Lemma 1]{ferguson1967mathematical}) to obtain, for each $n\in\N$ and the integrable $h_n$, that
    \begin{equation*}
        M_{\phi,\convexindex}\parenthesis*{\int_\Omega h_n\diff\P}
        \leq \int_\Omega M_{\phi,\convexindex}\parenthesis*{h_n}\diff\P.
    \end{equation*}
    Now, as $M_{\phi,\convexindex}$ is non-decreasing, we have that $M_{\phi,\convexindex}\circ h_n$ is a non-decreasing sequence of functions. Then monotone convergence and continuity of $M_{\phi,\convexindex}:C_\phi\to\Rbar_{\geq0}$ finish the proof:
    \begin{align*}
        M_{\phi,\convexindex}\parenthesis*{\int_\Omega h\diff\P}
        &= M_{\phi,\convexindex}\parenthesis*{\lim_{n\to\infty}\int_\Omega h_n\diff\P}
        = \lim_{n\to\infty}M_{\phi,\convexindex}\parenthesis*{\int_\Omega h_n\diff\P}\\
        &\leq \lim_{n\to\infty}\int_\Omega M_{\phi,\convexindex}\parenthesis*{h_n}\diff\P
        =\int_\Omega \lim_{n\to\infty}M_{\phi,\convexindex}\parenthesis*{h_n}\diff\P
        =\int_\Omega M_{\phi,\convexindex}\parenthesis*{h}\diff\P.
    \end{align*}
    For concave $M_\phi$, all inequalities are reversed.
\end{delayedproof}

\subsection{Details for \autoref{example: Counterexample Wasserstein p barycenter convexity for p neq 2}}
\label{sec: Details counterexample}

In \autoref{example: Counterexample Wasserstein p barycenter convexity for p neq 2}, we use the mixture distributions $\popmeasure=\frac13\sum_{j=1}^3\delta_{\mu_j}$ and $\tilde{\popmeasure}=\frac13\sum_{j=1}^3\delta_{\tilde{\mu}_j}$ as given in \autoref{tab: uniform counterexmaples p neq2}.
As $\popmeasure$ is finitely supported, we can use \cite[Theorem 5.1]{brizziWassersteinBarycenters2025} ($1<p<+\infty$) and \cite[Proposition 4.1.]{carlierWassersteinMediansRobustness2024} ($p=1$) to obtain
\begin{equation}
    \label{eqn: Wasserstein p barycenter characterization quantile functions}
    \iota_t(\Quantiletransform(\Barycenterfunc{p}{\Wp}(\popmeasure)))=\iota_t(\Lambda_p^{\Lpmod^p(0,1)}(\pushforward{\Quantiletransform}{\popmeasure}))=\Lambda_p^{\R}(\pushforward{(\iota_t\circ\Quantiletransform)}{\popmeasure}),\quad t\in(0,1),p\in[1,\infty),
\end{equation}
{i.e.}, the {q.f.\@} of the \WB{p} is the pointwise $p$-barycenter on $\R$ of the evaluated {q.f.s}{}. For $p=1$, this equation is set-valued in general. In our example, this will not be the case, as the median $\Lambda_{1}^{\R}$ of three points is unique.
\begin{remark}%
    Equation \eqref{eqn: Wasserstein p barycenter characterization quantile functions} is the analogue of \eqref{eqn: characterization of Qunaitle fucntionas mean} for $p\neq2$ with the restriction to $\popmeasure\in\Wp(\Wp(\R))$ satisfying $\setcount{\support{\popmeasure}}\in\N$.
    The extension of this result to arbitrary $\popmeasure\in\Wp(\Wp(\R))$ has not yet been reported in the literature but is straightforward for $p>1$ with additional care needed for $p=1$. 
\end{remark}
Fixing a distribution $\rho=\frac13\sum_{j=1}^3\delta_{x_j}$ on $\R$, we have that $\Lambda_p^{\R}(\rho)$ is monotone in $p$ with extremal values
\begin{align*}
    \Lambda_1^{\R}(\rho)      & =\median{x_1,x_2,x_3}=x_{(2)},                     \\
    \Lambda_{\infty}^\R(\rho) & =\midrange(x_1,x_2,x_3)=\frac{x_{(1)}+x_{(3)}}{2},
\end{align*}
where $x_{(1)}\leq x_{(2)}\leq x_{(3)}$ are the order statistics of $(x_1,x_2,x_3)\in\R^3$. The function $\Lambda_{p}^\R(\rho)$ is strictly monotone in $p\in[1,+\infty]$ if and only if $\Lambda_1^{\R}(\rho)\neq\Lambda_{\infty}^\R(\rho)$.
If $\Lambda_{p}^\R(\rho)$ is not strictly monotone, then it is constant in $p$ and given by $\median{x_1,x_2,x_3}=\Lambda_p^{\R}(\rho)=\midrange(x_1,x_2,x_3)$.
Applied to $\popmeasure$, we have that
\begin{align*}
    Q_1(t)       \defeq{} & \median(\set{Q_{{\mu}_1}(t),Q_{{\mu}_2}(t),Q_{{\mu}_3}(t)})< Q_p(t)                                                  \\
    <{}                   & Q_\infty(t)  \defeq\midrange(\set{Q_{{\mu}_1}(t),Q_{{\mu}_2}(t),Q_{{\mu}_3}(t)}),\quad t\in\parenthesis*{0,\frac13},
\end{align*}
which gives an upper and a lower bound for the {q.f.\@} $Q_p$ of the \WB{p} of $\popmeasure$ for each $p\in(1,\infty)$, with the roles switching on $(\frac{1}{3},\frac{2}{3})$ and again on $(\frac{2}{3},1)$. Similarly, for $\tilde{\popmeasure}$, we have $\tilde{Q}_1>\tilde{Q}_p>\tilde{Q}_\infty$ on $(0,\frac13)\cup (\frac23,1)$ and $\tilde{Q}_1<\tilde{Q}_p<\tilde{Q}_\infty$ on $(\frac13,\frac23)$.
By construction, we have that $Q_1=\tilde{Q}_{\infty}$. Furthermore, the \WB{2} $\mu^*_2$ and $\tilde{\mu}_2^*$ also coincide, as depicted on the {l.h.s.\@} of \autoref{fig: counterexample unimodal preservation for p neq 2 and d>1}.
Thus, for each $p\in(1,\infty)\setminus{\set{2}}$, there exists a distribution supported on uniform distributions, namely either $\popmeasure$ for $p<2$ or $\tilde{\popmeasure}$ for $p>2$, such that the {q.f.\@} $Q_{\mu_p^*}$ of its \WB{p} lies strictly inside the region bounded by $Q_1$ and $Q_{2}$, except at the crossing points $\frac13,\frac23$.

\end{document}